\documentclass{article}%
\usepackage{amsmath}
\usepackage{amsfonts}
\usepackage{amssymb}
\usepackage{tikz}
\usepackage{graphicx}

\usepackage{fancyhdr}
\usepackage{hyperref}

\newtheorem{theorem}{Theorem}

\newtheorem{corollary}[theorem]{Corollary}

\newtheorem{definition}[theorem]{Definition}

\newtheorem{lemma}[theorem]{Lemma}

\newtheorem{proposition}[theorem]{Proposition}
\newtheorem{remark}[theorem]{Remark}

\newenvironment{proof}[1][Proof]{\noindent\textbf{#1.} }{\ \rule{0.5em}{0.5em}}

\begin{document}

\title{Computing the exact number of periodic orbits for planar flows}

\author{Daniel S.~Gra\c{c}a\\Universidade do Algarve, C. Gambelas, 8005-139 Faro, Portugal\\ and Instituto de Telecomunica\c{c}\~{o}es, Portugal \and
Ning Zhong\\DMS, University of Cincinnati, Cincinnati, OH 45221-0025, U.S.A.}

\maketitle

\begin{abstract}
In this paper, we consider the problem of determining the \emph{exact} number
of periodic orbits for polynomial planar flows. This problem is a variant of
Hilbert's 16th problem. Using a natural definition of computability, we show
that the problem is noncomputable on the one hand and, on the other hand,
computable uniformly on the set of all structurally stable systems defined on
the unit disk. We also prove that there is a family of polynomial planar
systems which does not have a computable sharp upper bound on the number of
its periodic orbits.
\end{abstract}

\section{Introduction\label{Sec:Introduction}}

\thispagestyle{fancy}
In his famous lecture of the 1900 International Congress of Mathematicians,
David Hilbert stated a list of 23 problems. There has been intensive research
on these problems ever since. The second part of Hilbert's 16th problem asks
for the maximum number and relative positions of periodic orbits of planar
polynomial (real) vector fields of a given degree%
\begin{equation}
x^{\prime}=p(x), \label{ODE_poly}%
\end{equation}
where the components of $p:\mathbb{R}^{2}\rightarrow\mathbb{R}^{2}$ are
polynomials of degree $n$. More than a century later, and despite extensive
work on this topic (see \cite{Ily02} for an overview of its rich history),
this problem remains open even for the simplest non-linear systems where $p$
consists of quadratic polynomials. In this paper, we investigate the following
related problem:\bigskip

\textbf{Problem.} Is there some general procedure that, given as input a
function $p:\mathbb{R}^{2}\rightarrow\mathbb{R}^{2}$ with polynomial
components of a given degree, yields as output the number and relative
positions of periodic orbits of (\ref{ODE_poly})?\bigskip

A \textquotedblleft general procedure\textquotedblright\ is referred to a formula or an algorithm. We show that the answer to
this problem is negative:\bigskip

\textbf{Theorem A.} The operator $\Theta$ which maps a function $p:A\subseteq\mathbb{R}%
^{2}\rightarrow\mathbb{R}^{2}$ with polynomial components to the number
of periodic orbits of (\ref{ODE_poly}) is noncomputable when:
\begin{enumerate}
    \item $A$ is the unit ball;
    \item $A=\mathbb{R}^2$. In this case there exists a family of polynomials $\{p_n\}_{n\in\mathbb{N}}$ and a value $\delta>0$ such that $\|p_n-p_m\|\geq\delta$ whenever $n\ne m$, with the property that the operator $\Theta$ is still noncomputable over the set $\{p_n:n\in\mathbb{N}\}$.
\end{enumerate}\bigskip

In the second item of the theorem, we notice that the existence of a family $\{p_n\}_{n\in\mathbb{N}}$ of polynomials which are not close to each other but on which $\Theta$ is noncomputable shows that noncomputability can arise even if continuity problems are avoided (it is well known that, over real numbers, discontinuous functions are also noncomputable).

On the other hand, there is an algorithm that computes the number and
depicts the positions of periodic orbits - the portraits can be made with
arbitrarily high precision - for any structurally stable vector field defined
on the closed unit disk; moreover, the computation is uniform on the set of
all such vector fields. This is our second main result, Theorem B. Recall that
the density theorem of Peixoto \cite[Theorem 2]{Pei62} shows that, on
two-dimensional compact manifolds, structurally stable systems are
\textquotedblleft typical\textquotedblright\ in the sense that such systems
form a dense open subset in the set of all $C^{1}$ systems%
\begin{equation}
x^{\prime}=f(x)\text{.} \label{ODE_Main}%
\end{equation}
Moreover,  structurally stable system can only have a finite number of  equilibrium points and of periodic orbits.

\bigskip

\textbf{Theorem B.} Let $\mathbb{D}\subseteq\mathbb{R}^{2}$ be the closed unit
disk and let $SS_{2}$ be the subset of $\mathcal{X}(\mathbb{D})$ consisting of all $C^{1}$ structurally stable vector
fields $f:\mathbb{D}\rightarrow\mathbb{R}^{2}$ (the definition of $\mathcal{X}(\mathbb{D})$ is given in subsection 3.1). Then the operator which maps
$f\in SS_{2}$ to the number of periodic orbits of (\ref{ODE_Main}) is (uniformly)
computable. Meanwhile, the algorithm which produces the computation can depict
the periodic orbits with arbitrarily high precision.\bigskip

Another related problem is to find a \textit{sharp} upper bound (see Section
\ref{Sec:ThmAProof} for definition) for the number of periodic orbits that a
polynomial system (\ref{ODE_poly}) of degree $n$ can have. We show that this
problem is in general not computable.\bigskip

\textbf{Theorem C.} There is a family of polynomial systems (\ref{ODE_poly}), namely the family $\{p_n\}_{n\in\mathbb{N}}$ of Theorem A,
which does not have a computable sharp upper bound on the number of its
periodic orbits.\bigskip

The structure of the paper is as follows. In Section \ref{Sec:Computability},
we discuss classical computability theory and computability over the real numbers.
In Section \ref{Sec:StructuralStability},\ we review some notions about
structurally stable systems. For completeness of the paper, several
\textquotedblleft folklore\textquotedblright\ results are presented in this
section which are not, up to our knowledge, coherently presented elsewhere in
the literature (Appendix \ref{Sec:AppendixDynSys} includes some proofs of these results). In Section \ref{Sec:ThmAProof},  we prove Theorems A and C.
Section \ref{Sec:ThmBProof} presents an outline for the proof of Theorem B,
while the remaining sections present all the technical parts of the argument. Section \ref{Sec:Hilbert16} discusses possible connections of the present work with Hilbert's 16th problem and proves that Hilbert's 16th problem is computable, relative to the Halting problem, over a dense and open subset, both for vector fields defined over the unit ball $\mathbb{D}$ or over the whole plane $\mathbb{R}^2$. The conclusion reviews the main results of the paper and discusses some open problem.

\section{Introduction to computability theory\label{Sec:Computability}}

In this section we briefly outline the classical computability theory and
computability over the reals. The presentation draws heavily from
\cite[Section 2.1]{BY06}.

\textbf{Classical Computability.} Computability theory allows one to classify
problems as algorithmically solvable (computable) or algorithmically
unsolvable (noncomputable). For example, most common computational tasks such
as performing arithmetic operations with integers, finding whether a graph is
connected, etc. are all computable.

A major contribution of computability theory is to show the existence of
noncomputable problems. The best known examples of noncomputable problems (see
e.g. \cite{Sip12}, \cite{Mat93})\ are the \textit{Halting problem} and \textit{the solvability of diophantine equations} (Hilbert's 10th problem).

In the setting of formal computability theory, computations are performed on
Turing machines (TM for short), which were introduced by Alan Turing in 1936
\cite{Tur36}; the notion of Turing machines has since become a universally
accepted formal model of computation. A function $u:\mathbb{N}^{k}%
\rightarrow\mathbb{N}^{l}$ is computable if there exists a TM that takes $a$
as an input and outputs the value of $u(a)$. There are only countably many Turing
machines, which can be enumerated in a natural way. (See, e.g., \cite{Sip12} and references therein for more details.)

Since Turing machines solve exactly the same problems which are
solvable by digital computers, it suffices to regard a TM as a
computer program written in any programming language. We will often
take this approach in the paper.

This notion of computability can be naturally extended to the
rational numbers $\mathbb{Q}$, some countable subsets of the real
numbers, or any domain that can be ``effectively encoded" in
$\mathbb{N}$. For example, mathematical expressions consisting of
finitely many symbols can be symbolically manipulated as performed
by computer algebra systems (there are only countably many such
expressions, viewed as strings of finite length). Another example is
that any finite union of balls or rectangles in $\mathbb{R}^n$
having rational radii and centers with rational coordinates or
having corners with rational coordinates are computable objects. On
the other hand, it is clear that $\mathbb{R}$ - the set of all real
numbers - is too big to be encoded in $\mathbb{N}$. Computability of
the real numbers and real functions is the subject studied in
computable analysis.

\textbf{Computability of real functions and sets.} Computable
analysis was originated from the work of Banach and Mazur
\cite{BM37, Maz63}. There are several equivalent modern approaches;
for example, the axiomatization approach \cite{PR89}, the type two
theory of effectivity or representation approach \cite{BHW08, Wei00}
(the approach used in this paper), and the oracle Turing machine
approach \cite{Ko91, BY09}. These approaches provide a common
framework for combining approximation, computation, computational
complexity, and implementation. Roughly speaking, in this model of
computation, an object is computable if it can be approximated by
computer-generated approximations with an arbitrarily high
precision.

Formalizing this idea to carry out computations on infinite objects,
those objects are encoded as infinite sequences of finite-sized
approximations with arbitrary precision,
using representations (see \cite{BHW08, Wei00} for a complete
development). Let $\Sigma$ be a countable set of symbols. A
represented space is a pair $(X; \delta)$, where $X$ is a set,
$\delta: \Sigma^{\mathbb{N}}\to X$ is an onto map, and
$\mbox{dom}(\delta)\subseteq\Sigma^{\mathbb{N}}$. Every
$q\in\mbox{dom}(\delta)$ such that $\delta(q)=x$ is called a name
(or a $\delta$-name) of $x$. Note that $q$ is an infinite sequence
in $\Sigma$. In majority of cases, $\delta(q)$ is an infinite
sequence in $X$ with finite-sized members that converges to $x$ in
$X$ with a prescribed rate. For instance, a name of a $C^k$ function
$f: \mathbb{R} \to \mathbb{R}$, $k\geq 0$, can be taken as an
infinite sequence of finite-sized functions such as polynomials
$P_{l}$ with rational coefficients that satisfies $\| f(x) -
P_l(x)\|_{k}\leq 2^{-l}$ for all $l\in \mathbb{N}$,  where $\|
\cdot\|_k$ is a $C^k$-norm.  Naturally, an element $x\in X$ is
computable if it has a computable name in $\Sigma^{\mathbb{N}}$;
namely, it has a name that is classically computable. For example, a
popular name for a real number $x$ is a sequence $\{ r_{l}\}$ of
rationals satisfying $|x-r_{l}| \leq 2^{-l}$. Thus, $x$ is
computable if there is a Turing machine (or a computer program or an
algorithm) that outputs a rational $r_{l}$ on input $l$ such that
$|r_{l}-x|\leq 2^{-l}$ for all $l\in \mathbb{N}$.

The notion of computable maps between represented spaces now arises naturally.
A map $\Phi: (X;\delta_{X})\to(Y;\delta_{Y})$ between two represented spaces
is computable if there is a (classically) computable map $\phi:\subseteq\Sigma^{\mathbb{N}%
}\to\Sigma^{\mathbb{N}}$ such that
$\Phi\circ\delta_{X}=\delta_{Y}\circ\phi$. Informally speaking, this
means that there is a computer program that outputs a name of $\Phi
(x)$ when given a name of $x$ as input. Thus, for example, a $C^k$
function $f: \mathbb{R} \to \mathbb{R}$, $k\geq 0$, is computable if
there exists a machine capable of computing an approximation $p_l$
of $f(x)$ satisfying $\| f(x) - p_l\|_k \leq 2^{-l}$  when given as
input of $l\in \mathbb{N}$ (accuracy) and (a name of) $x\in
\mathbb{R}$. The precise definition of a computable $C^k$ function
is given below. Since only the planar vector fields are considered
in this paper, the definition is given to planar functions defined
on $\mathbb{R}^2$ only. The definitions A and B are equivalent.
\begin{definition} \label{computable-function} Let $f: \mathbb{R}^2 \to \mathbb{R}^2$ be a $C^k$ function, and let $K=\{ x\in \mathbb{R}^2: \, \|x\|\leq r \}$, $r$ is a rational number.
\begin{itemize}
\item[A.]  $f$ is said to be ($C^k$-) computable on $K$ if there is a Turing machine  that, on input $l$ (accuracy), outputs the rational coefficients of a polynomial $P_l$ such that $d^k(f, P_l)\leq 2^{-l}$, where
\[ d^k(f, P_l)=\max_{0\leq j\leq k}\max_{x\in K} \|D^{j}f(x)-D^{j}P_{l}(x)\| \]
\item[B.] $f$ is said to be ($C^k$-) computable on $K$ if there is an oracle Turing machine such that for any input $l\in \mathbb{N}$ (accuracy) and any name of $x\in K$ given as an oracle, the machine will output the rational vectors $q_0, q_1, \ldots, q_k$ in $\mathbb{R}^2$ such that $\|q_j - D^{j}f(x)\|\leq 2^{-l}$ for all $0\leq j\leq k$ (see e.g. \cite{Ko91, BY06}).
\end{itemize}
\end{definition}
In practice, an oracle can be conveniently treated as an interface
to a program computing $f$: for every $l\in\mathbb{N}$, the oracle
supplies a good enough rational approximation $p$ of $x$ to the
program, the program then performs computations based on inputs $l$
and $p$, and returns rational vectors $q_j$, $0\leq j\leq k$, in the
end such that $\|q_j-D^{j}f(x)\|<2^{-l}$. In other words, with an
access to arbitrarily good approximations for $x$, the machine
should be able to produce arbitrarily good approximations for $f(x),
Df(x), \ldots, D^k f(x)$. This is often termed as $f(x)$ is computable
from (a name of) $x$.

\begin{definition} \label{computable-compact-set} Let $C$ be a compact subset of $\mathbb{R}^2$. Then $C$ is said to be computable if there is a Turing machine that, on input $l\in \mathbb{N}$ (accuracy), outputs finite sequences $r_j\in \mathbb{Q}$ and $c_j\in \mathbb{Q}^2$, $1\leq j\leq j_l$, such that $d_{H}(C, S_l)\leq 2^{-l}$, where $S_l = \cup_{j=1}^{j_l}B(c_j, r_j)$ is the finite union of the closed balls $B(c_j, r_j)=\{ x\in \mathbb{R}^2: \, \| x - c_j\|\leq r_j\}$ and $d_{H}(\cdot, \, \cdot)$ denotes the Hausdorff distance between two compact subsets of $\mathbb{R}^2$.
\end{definition}
Thus, if we imagine those rational balls as pixels, then $C$ is computable provided it can be drawn on a computer screen with arbitrarily high precision.

We mention in passing  that although the computation can only exploit approximations up to some finite precision in a finite time, nevertheless it is always possible to continue the computation with better approximations of the input. In other words, the computation can be performed to achieve arbitrary precision and obtain results which are guaranteed correct.

We also note that many standard functions like arithmetic operations ($+, \times, \ldots$), polynomials with computable coefficients, the usual trigonometric functions $\sin$, $\cos$, the exponential $e^x$, their composition, etc.~are all computable \cite{BHW08}. Many other standard operations are also well-known to be computable. In particular, the operator which yields the solution of some initial-value problem (\ref{ODE_Main}) with $x(t_0)=x_0$ is also computable \cite{GZB07}, \cite{CG08}, \cite{CG09} and one can also often determine bounds on the computational resources needed to compute it (see e.g.~\cite{Mul87}, \cite{MM93} \cite{Kaw10}, \cite{BGP12}, \cite{KORZ14}, \cite{PG16}, \cite{GZ21}).

\section{Structural stability\label{Sec:StructuralStability}}

\subsection{Classical theory}

We recall the definition of structurally stable systems for the case of flows
(see e.g. \cite[pp. 317-318]{Per01}). Let $K$ be a compact set in
$\mathbb{R}^{n}$ with non-empty interior and smooth $(n-1)$-dimensional
boundary. Consider the space $\mathcal{X}(K)$ consisting of restrictions to
$K$ of $C^{1}$ vector fields on $\mathbb{R}^{n}$ that are transversal to the
boundary of $K$ and inward oriented. $\mathcal{X}(K)$ is equipped with the
norm%
\[
\left\Vert f\right\Vert _{1}=\max_{x\in K}\left\Vert f(x)\right\Vert
+\max_{x\in K}\left\Vert Df(x)\right\Vert
\]
where $\| \cdot\|$ is either the max-norm $\| x\|=\max\{ |x_{1}|, \ldots,
|x_{n}|\}$ or the $l^{2}$-norm $\|x\|=\sqrt{x_{1}^{2}+\cdots+ x_{n}^{2}}$;
these two norms are equivalent.

\begin{definition}
\label{Def:StructuralStability}The system (\ref{ODE_Main}), where
$f\in\mathcal{X}(K)$, is structurally stable if there exists some
$\varepsilon>0$ such that for all $g\in C^{1}(K)$ satisfying $\left\Vert
f-g\right\Vert _{1}\leq\varepsilon$, the trajectories (orbits) of%
\begin{equation}
y^{\prime}=g(y) \label{Eq:PertubedDynSys}%
\end{equation}
are homeomorphic to the trajectories of (\ref{ODE_Main}), i.e. there exists
some homeomorphism $h$ such that if $\gamma$ is a trajectory of
(\ref{ODE_Main}), then $h(\gamma)$ is a trajectory of (\ref{Eq:PertubedDynSys}%
). Moreover, the homeomorphism $h$ preserves the orientation of trajectories
with time.
\end{definition}

Intuitively, (\ref{ODE_Main}) is structurally stable if the shape of its
dynamics is robust to small perturbations. We now recall the notion of
non-wandering set.

\begin{definition}
The non-wandering set $NW(f)$ of (\ref{ODE_Main}) is the set of all points $x$
with the following property: for any neighborhood $U$ of $x$, given some
arbitrary $T>0$, there exists $t\geq T$ such that $\phi_{t}(U)\cap
U\neq\varnothing$, where $\phi_{t}(U)=\{\phi_{t}(y):y\in U\}$, and $\phi
_{t}(y)$\ is the solution of (\ref{ODE_Main}) with the initial condition
$x(0)=y$.
\end{definition}

For a structurally stable planar vector field $f$, the set $NW(f)$ consists of
equilibrium points and periodic orbits. A point $x_{0} \in K$ is called an
equilibrium (point) of the system (\ref{ODE_Main}) if $f(x)=0$. Accordingly
any trajectory starting at an equilibrium stays there for all $t\in\mathbb{R}%
$. An equilibrium $x_{0}$ is called a hyperbolic equilibrium if the
eigenvalues of $Df(x_{0})$ have non-zero real parts. If both eigenvalues of
$Df(x_{0})$ have negative real parts, $x_{0}$ is called a sink - it attracts
nearby trajectories; if both eigenvalues have positive real parts, $x_{0}$ is
called a source - it repels nearby trajectories; if the real parts of the
eigenvalues have opposite signs, $x_{0}$ is called a saddle. A system is
locally robust near a hyperbolic equilibrium. Indeed, it follows from the
Hartman-Grobman theorem that the nonlinear vector field $f(x)$ is conjugate to
its linearization $Df(x_{0})$ in a neighborhood of $x_{0}$ provided that
$x_{0}$ is a hyperbolic equilibrium.

A closed curve $\gamma$ in $NW(f)$ is called a \textit{periodic orbit} if there is some $T>0$ such that for any $x\in\gamma$ one
has $\phi_{T}(x)=x$. Periodic orbits can also be hyperbolic, with similar
properties as of hyperbolic equilibria. However, there are only attracting periodic orbits and repelling periodic orbits, both are pictured in
Fig.\ \ref{fig:nested_orbits}. There is no equivalent of a saddle point for
periodic orbits in dimension two (one dimension is \textquotedblleft used
up\textquotedblright\ by the flow of the periodic orbit. The remaining
direction can only be attracting or repelling). See \cite[p. 225]{Per01} for
more details.

The following well-known theorem proved by Peixoto in 1962 \cite{Pei62} is a
refinement of the Poincar\'{e}-Bendixson theorem.

\begin{theorem}
[Peixoto]\label{Th:GeneralPeixoto}Let $f$ be a $C^{1}$ vector field defined on
a compact two-dimensional differentiable manifold $K\subseteq\mathbb{R}^{2}$.
Then $f$ is structurally stable on $K$ if and only if:

\begin{enumerate}
\item The number of equilibria (i.e.\ zeros of $f$) and of periodic orbits is
finite and each is hyperbolic;

\item There are no trajectories connecting saddle points, i.e.\ there are no
\emph{saddle connections};

\item The non-wandering set $NW(f)$ consists only of equilibria and periodic orbits.
\end{enumerate}

Moreover, if $K$ is orientable, the set of structurally stable vector fields
in $C^{1}(K)$ is an open, dense subset of $C^{1}(K)$. Similar results hold for
$\mathbb{D}=\{x\in\mathbb{R}^{2}:\left\Vert x\right\Vert \leq1\}$, assuming
that the vector fields always point inwards on the boundary of $\mathbb{D}$
(see \cite[Theorem 3 of p. 325]{Per01}, \cite{Pei59}, and the references
therein), with the difference that condition 3 is not needed (it follows from the Poincaré-Bendixson theorem). It is mentioned in \cite[p.~200]{Pei59} that the results remain true
if $\mathbb{D}$ is replaced by any region bounded by a $C^{1}$ Jordan curve.
\end{theorem}

Due to the above result, we will always assume that the vector field points
inwards along the boundary of $\mathbb{D}$.

\subsection{Results about structurally stable
systems\label{Sec:StructuralResults}}

In this section, we present several results which are needed for proving
Theorem B. Most ideas and results of this section are implicitly present in
the literature or are \textquotedblleft folklores\textquotedblright. However,
as details are usually important when working with computability theory, and
because, up to our knowledge, these results are not presented coherently
elsewhere, we decided to include a section with these results. Readers
familiar with structural stability and related fields may skip this section.
For completeness, we have included a proof for each result if no reference is
provided in Appendix \ref{Sec:AppendixDynSys}.

The following theorems can be found in \cite[Theorem 1 of p. 130 and Theorem 3
of p.\ 226]{Per01}. They state that there is a neighborhood, called the basin
of attraction, around an attracting hyperbolic equilibrium point or hyperbolic
periodic orbit (the attractor) such that the convergence to the attractor is
exponentially fast in this neighborhood. Let $\mathcal{N}_{\bar{\epsilon}}(A)$
denote the $\bar{\epsilon}$-neighborhood of $A$:%
\[
\mathcal{N}_{\bar{\epsilon}}(A)=%
{\displaystyle\bigcup\limits_{x\in A}}
\mathring{B}(x,\bar{\epsilon})=%
{\displaystyle\bigcup\limits_{x\in A}}
\{y\in\mathbb{D}:\left\Vert x-y\right\Vert <\bar{\epsilon}\}.
\]

\begin{proposition}
\label{Prop:ConvEquilibrium}Let $x_{0}$ be a sink of
(\ref{ODE_Main}) such that $\operatorname{Re}(\lambda_{i})\leq-\alpha<0$ for
all eigenvalues $\lambda_{i}$ of $Df(x_{0})$. Then given $\varepsilon>0$,
there exists $\delta>0$ such that for all $x\in\mathcal{N}_{\delta}(x_{0})$,
the flow $\phi_{t}(x)$ of (\ref{ODE_Main}) satisfies%
\[
\left\vert \phi_{t}(x)-x_{0}\right\vert \leq\varepsilon e^{-\alpha t}%
\]
for all $t\geq0$.
\end{proposition}

\begin{proposition}
\label{Prop:ConvPeriodic}Let $\gamma=\gamma(t)$ be an attracting
periodic orbit of (\ref{ODE_Main}) with period $T$. Then there exist some
$\alpha>0$, $\delta>0$ and $\varepsilon>0$ such that for any $x\in
\mathcal{N}_{\delta}(\gamma)$, there is an asymptotic phase $t_{0}$ such that
for all $t\geq0$%
\[
\left\vert \phi_{t}(x)-\gamma(t-t_{0})\right\vert \leq\varepsilon e^{-\alpha
t/T}.
\]

\end{proposition}

\begin{theorem}
\label{Cor:uniform_convergence_full}Let (\ref{ODE_Main})\ be structurally
stable and defined on the compact set $\mathbb{D}\subseteq\mathbb{R}^2$, and suppose it does not have any
saddle equilibrium point. Then for every $\epsilon>0$ there exists some $T>0$
such that for any $x\in\mathbb{D}$ satisfying $d(x,NW(f))\geq\epsilon$ one has
$d(\phi_{t}(x),NW(f))\leq\epsilon$ for every $t\geq T$.
\end{theorem}

\begin{remark}
\label{Re:no_saddle} Theorem \ref{Cor:uniform_convergence_full} is relevant to prove Theorem B in a simplified case where (\ref{ODE_Main}) has no saddle points, since it shows that for any given
precision $n$, there always exists some time $T_{n}>0$ such that for every
$x\in\mathbb{D}$, either $x$ is already in $\mathcal{N}_{1/n}(NW(f))$ or
$\phi_{t}(x)$ will enter $\mathcal{N}_{1/n}(NW(f))$ no later than $T_{n}$ and
stay in $\mathcal{N}_{1/n}(NW(f))$ thereafter. This \textquotedblleft
uniform\textquotedblright\ time bound is an essential ingredient in
constructing the algorithm presented in section 8.1.

When (\ref{ODE_Main}) has saddle points, the situation becomes more
complicated because a saddle does not have a basin of attraction: no matter
how small a neighborhood of a saddle is, almost all trajectories entering the
neighborhood will leave it after some time; even worse, there is no
\textquotedblleft uniform\textquotedblright\ bound on time sufficient for the
trajectories to leave the neighborhood.
\end{remark}

The following results provide some useful tools to handle the saddles.

\begin{lemma}
(See e.g. \cite[Theorem 3, p. 177]{BR89}) \label{Lemma:unstable_convergence}%
Let (\ref{ODE_Main}) define a structurally stable system over the compact set
$\mathbb{D}\subseteq\mathbb{R}^{2}$, and let $x_{0}$ be a saddle point.
Then there are at most two attractors $\Omega_{1}(x_{0}),\Omega_{2}(x_{0})$,
with $\Omega_{1}(x_{0})\cap\Omega_{2}(x_{0})=\varnothing$, such that any
trajectory starting on $U_{x_{0}}-\{x_{0}\}$, where $U_{x_{0}}$ is a local
unstable manifold of $x_{0},$ will converge to one of the attractors
$\Omega_{1}(x_{0}),\Omega_{2}(x_{0})$.

\end{lemma}

\begin{theorem}
\label{Thm:No_saddle}Let (\ref{ODE_Main}) define a structurally stable system
over the compact set $\mathbb{D}\subseteq\mathbb{R}^{2}$ with saddle points
$x_{1},\ldots,x_{m}$. Then there exists some $\varepsilon>0$ such that for every $1\leq i\leq m$,  any trajectory starting in
$B(x_{i},\varepsilon)$ will never intersect $B(x_{j},\varepsilon)$, where
$j=1,\ldots,i-1,i+1,\ldots,m.$
\end{theorem}

We end this section with two more lemmas.

\begin{lemma}
(\cite{GH83}) \label{Lemma:equilibrium_point_cycle}Given a periodic orbit
$\gamma$ of a structurally stable dynamical system (\ref{ODE_Main}) defined on
a compact subset of $\mathbb{R}^{2}$, there is always an equilibrium point in
the interior of the region delimitated by $\gamma$.
\end{lemma}

\begin{corollary}
\label{corollary:Nonzero_zero} The previous lemma implies that for every $f\in SS_{2}$,
the system (\ref{ODE_Main}) has at least one equilibrium point.
\end{corollary}

A sketch of proof: Suppose otherwise. The Poincar\'{e}-Bendixson Theorem
implies that if a trajectory does not leave a closed and bounded
region of phase space which contains no equilibria, then the trajectory must
approach a periodic orbit as $t\to\infty$ (see, e.g. \cite{GH83}). Thus there
is at least one periodic orbit inside $\mathbb{D}$. Then it follows from Lemma
\ref{Lemma:equilibrium_point_cycle} that there is at least one equilibrium
inside the region delimitated by the periodic orbit. We arrive at a contradiction.

The next result can be found e.g.~in \cite{BR89}.
\begin{lemma}
\label{Lemma:Perturbation}Assume that $x$ and $y$ are, respectively, the
solutions of the ODEs (\ref{ODE_Main}) and (\ref{Eq:PertubedDynSys}) with
initial conditions set at $t_{0}=a$, which are defined on a domain
$D\subseteq\mathbb{R}^{n}$, where $f$ and $g$ are continuous and satisfy
$\left\Vert f(x)-g(x)\right\Vert \leq\varepsilon$ for every $x\in D$. Suppose $f$ also satisfies a Lipschitz condition in $D$ with Lipschitz constant
$L$. Then%
\[
\left\Vert x(t)-y(t)\right\Vert \leq\left\Vert x(a)-y(a)\right\Vert
e^{L\left\vert t-a\right\vert }+\frac{\varepsilon}{L}\left(  e^{L\left\vert
t-a\right\vert }-1\right)  .
\]

\end{lemma}

Let $x(\cdot,x_0)$ be the solution of (\ref{ODE_Main}) with initial condition $x(0)=x_0$. Note that the previous lemma implies that given any $x_0\in\mathbb{D}, T\geq0$, and $\delta>0$, one has that $\|x_0-y_0\|\leq\delta e^{-LT}$ yields $\|x(T,x_0)-x(T,y_0)\|\leq\delta$. This fact implies the following corollary.

\begin{corollary}
\label{Cor:Continuity}Let $T>0$. The maps $\phi_{T}:\mathbb{D}\rightarrow
\mathbb{D}$ and $\phi_{-T}:\phi_{T}(\mathbb{D})\rightarrow\mathbb{D}$ are continuous, where $\phi_t(x_0)=x(t,x_0)$ and $x(\cdot,x_0)$ is the solution of (\ref{ODE_Main}) with initial condition $x(0)=x_0$.
\end{corollary}

\section{Proof of Theorems A and C\label{Sec:ThmAProof}}

We begin this section by reviewing a variant of the halting problem and its
proof. Recall from Section \ref{Sec:Computability} that it is possible to
enumerate all the Turing machines in a natural (computable) manner:\ $TM_{1}%
,TM_{2},\ldots$.

\begin{lemma}
\label{Lem:NonComputable}The function
\[
h(k)=\left\{
\begin{array}
[c]{ll}%
1\text{ \ \ } & \text{if }TM_{k}\text{ halts on }k\\
0 & \text{if }TM_{k}\text{ doesn't halt on }k
\end{array}
\right.
\]
is not computable.
\end{lemma}

\begin{proof}
Suppose otherwise $h$ was computable. Then there is a Turing machine $TM_{j}$
computing it. Let $TM_{i}$ be the following machine: given some input $k$ it
simulates $TM_{j}$ with this input. If $TM_{j}$ outputs 1, then $TM_{i}$
enters an infinite loop and never halts; if $TM_{j}$ outputs 0, then $TM_{i}$
halts and outputs 1. In other words:
\[
TM_{i}(k)=\left\{
\begin{array}
[c]{ll}%
\text{enters an infinite loop \ \ } & \text{if }TM_{k}\text{ halts on }k\\
\text{halts and outputs 1} & \text{if }TM_{k}\text{ doesn't halt on }k.
\end{array}
\right.
\]
What is the output of running $TM_{i}(i)$? If $TM_{i}$ halts on input $i$,
then $TM_{i}(i)$ enters an infinite loop;\ if $TM_{i}(i)$ does not halt, then
it halts and outputs 1. We arrive at a contradiction.
\end{proof}

We now prove Theorem A of Section \ref{Sec:Introduction}. We start with the case of the unit ball and we argue by way of a
contradiction. Suppose the operator $\Theta$ which maps $p$ to the number of
periodic orbits of (\ref{ODE_poly}) was computable.

Let $g: \mathbb{N}^{2} \to\mathbb{N}$ be the function defined as follows:%

\[
g(k,i)=\left\{
\begin{array}
[c]{ll}%
1 & \mbox{if $TM_k$ halts in $\leq i$ steps on input $k$}\\
0 & \mbox{otherwise}
\end{array}
\right.
\]
This function is computable since its output can be computed in finite time.
Let $G:\mathbb{N}\rightarrow\lbrack0,1/2]$, $G(k)=\sum_{i=1}^{\infty
}g(k,i)/2^{i+1}$. Then $G$ is also a computable function with $\sum_{i=1}%
^{n}g(k,i)$ being a rational approximation of $G(k)$ with accuracy $2^{-n}$.
Moreover, $0<G(k)\leq1/2$ if $TM_{k}$ halts on $k$ and $G(k)=0$ if $TM_{k}$
doesn't halt on $k$.

We now define a family of polynomial systems with parameters $G(k)$:
$x^{\prime}=p_{k}(x)$, where $x=(x_{1},x_{2})\in\mathbb{R}^{2}$, $p_{k}%
(x_{1},x_{2})=(p_{k,1}(x_{1},x_{2}),p_{k,2}(x_{1},x_{2}))$,%
\[
p_{k,1}(x_{1},x_{2})=-x_{2}+x_{1}(x_{1}^{2}+x_{2}^{2}-G(k))
\]
and%
\[
p_{k,2}(x_{1},x_{2})=x_{1}+x_{2}(x_{1}^{2}+x_{2}^{2}-G(k)).
\]
Since $G:\mathbb{N}\to[0, 1]$ is computable, so is the function $P:
\mathbb{N}\to\mathcal{P}$, $k\mapsto p_{k}$, where $\mathcal{P}$ is the set of
functions $p: \mathbb{R}^{2} \to\mathbb{R}^{2}$ with polynomial components. By
assumption that $\Theta$ is computable, it follows that the composition
$\Theta \circ P: \mathbb{N}\to\mathbb{N}$ is a computable function.

In the polar coordinates, the system is converted to the following form: let
$\theta=t$,%
\[
dr/dt=r(r^{2}-G(k)),d(\theta)/dt=1
\]
Thus, if $TM_{k}$ doesn't halt on $k$, then $dr/dt=r^{3}$, and there is only
one equilibrium point at the origin and no periodic orbit; if $TM_{k}$ does
halt on $k$, then $dr/dt=r(r^{2}-G(k))$ and there is one periodic orbit and
one equilibrium point. In other words,
\[
\Theta\circ P(k)=\left\{
\begin{array}
[c]{ll}%
1\text{ \ \ } & \text{if }TM_{k}\text{ halts on }k\\
0 & \text{if }TM_{k}\text{ doesn't halt on }k.
\end{array}
\right.
\]
We arrive at a contradiction because it follows from Lemma
\ref{Lem:NonComputable} that $\Theta\circ P$ cannot be a computable function.

Notice that the above argument also proves the second item of Theorem A, with the exception of the existence of the family of polynomials $\{p_n\}_{n\in\mathbb{N}}$, since the same argument works over the whole plane $\mathbb{R}^2$. Before we show the existence of such a family $\{p_n\}_{n\in\mathbb{N}}$, we first proceed to the proof of Theorem C. The missing part of the second item of Theorem A will follow as a corollary of the proof of Theorem C.

\begin{definition}
Let $\{p_{i}\}_{i\in I}$ be a family of functions $p_{i}:\mathbb{R}%
^{2}\rightarrow\mathbb{R}^{2}$ whose components are polynomials. A \emph{sharp
upper bound} for the number of periodic orbits for the family of polynomial
ODEs of the form (\ref{ODE_poly}) is a function $f:\mathbb{N}\rightarrow
\mathbb{N}$ with the following properties:

\begin{enumerate}
\item If $p$ is a function from $\{p_{i}\}_{i\in I}$
with components of degree at most $n$, then the number of periodic orbits of
(\ref{ODE_poly}) is $\leq f(n)$;

\item There is a function $p$ in $\{p_{i}\}_{i\in
I}$ with components of degree at most $n$ such that (\ref{ODE_poly}) has
exactly $f(n)$ periodic orbits.
\end{enumerate}
\end{definition}

\noindent \textbf{Theorem C.}  There is a family
$\{p_{k}\}_{k\in\mathbb{N}}$ of polynomial systems on
$\mathbb{R}^{2}$ for which there is no computable sharp upper bound
for the number of periodic orbits of (\ref{ODE_poly}).

\begin{proof}
Let $\{p_{k}\}_{k\in\mathbb{N}}$ be the family of polynomial systems
with
parameters $G(k)$: $dx/dt=p_{k}(x,y)$, where $p_{k}(x,y)=(p_{k,1}%
(x,y),p_{k,2}(x,y))$,
\begin{align*}
p_{k,1}(x,y)  &  =-y+x\Pi_{j=1}^{k}(x^{2}+y^{2}-\Sigma_{i=1}^{j}iG(i))\\
&  =-y+x(x^{2}+y^{2}-G(1))\cdots(x^{2}+y^{2}-(G(1)+2G(2)+\cdots+kG(k)))
\end{align*}
and
\[
p_{k,2}(x,y)=x+y\Pi_{j=1}^{k}(x^{2}+y^{2}-\Sigma_{i=1}^{j}iG(i))
\]
In the polar coordinates, the system is converted to the following one: let
$\theta= t$,
\[
dr/dt=r\Pi_{j=1}^{k}(r^{2}-\Sigma_{i=1}^{j}iG(i)),\text{ and }d(\theta)/dt=1.
\]
Again we argue by way of a contradiction. Suppose otherwise there was a
computable sharp upper bound $f$ for the number of periodic orbits for this
family of polynomial systems. We note first that $f$ is defined for every $n$.
It follows from the definition that the components of $p_{k}$ have degree
$1+2k$, and so $f(3)$ would yield a sharp upper bound for the number of
periodic orbits of (\ref{ODE_poly}) when $p=p_{1}$. In particular, if
$f(3)=0$, then $G(1)=0$, which implies that $TM_{1}$ does not halt on input
$1$ and thus $h(1)=0$; if $f(3)=1$, then $G(1)>0$ and thus $h(1)=1$. Now let
$k>1$. We observe that if $h(k)=0$, then $G(k)=0$ and thus
$f(1+2(k-1))=f(1+2k)$; on the other hand, if $h(k)=1$, then $G(k)>0$ and
$f(1+2k)=f(1+2(k-1))+1$. This observation together with the assumption that
$f$ is computable generates the following algorithm for computing $h(k)$ of
Lemma \ref{Lem:NonComputable}%
\[
h(k)=\left\{
\begin{array}
[c]{ll}%
1\text{ \ \ } & \text{if }f(1+2k)=f(1+2(k-1))+1\\
0 & \text{if }f(1+2k)=f(1+2(k-1))
\end{array}
\right.
\]
for $k>1$. But the function $h$ cannot be computable according to Lemma
\ref{Lem:NonComputable}. We arrive at a contradiction.
\end{proof}

We now remark that the proof of Theorem C can be used to prove the missing part of the second item of Theorem A.
Indeed, the family $\{p_{k}\}_{k\in\mathbb{N}}$ constructed in the
proof for Theorem C  has, at most, one polynomial vector field of
degree $n$ for every $n\in\mathbb{N}$. Hence, finding a computable
sharp upper bound for the number of periodic orbits of the family
$\{p_{k}\}_{k\in\mathbb{N}}$ is equivalent to finding an algorithm
that computes the exact number of periodic orbits of each $p_k$,
uniformly in $k$. Moreover, the construction also indicates that the
non-computability results are not simply the consequences of
discontinuity.


\section{Proof of Theorem B\label{Sec:ThmBProof}}

The following theorem is needed in order to prove Theorem B. Recall that
$\mathbb{N}$ denotes the set of all positive integers; $\mathbb{D}%
\subseteq\mathbb{R}^{2}$ the closed unit disk; $SS_{2}$ the set of all $C^{1}$
structurally stable planar vector fields defined on $\mathbb{D}$; and
$\mathcal{K}(\mathbb{D})$ the set of all non-empty compact subsets contained
in $\mathbb{D}$. Then $SS_{2}$ is a subspace of $C^{1}(\mathbb{D};
\mathbb{R}^{2})$ and $\mathcal{K}(\mathbb{D})$ is a metric space with the
Hausdorff metric. Recall that if $f\in SS_{2}$, then the trajectories of
(\ref{ODE_Main}) are transversal to the boundary of $\mathbb{D}$ and are
inward oriented.

\begin{theorem}
\label{Th:NW(f)} The operator $\Psi: SS_{2} \to\mathcal{K}(\mathbb{D})$,
$f\mapsto NW(f)$ of (\ref{ODE_Main}), is computable.
\end{theorem}

Note that it follows from Theorem \ref{Th:GeneralPeixoto} that $NW(f)$
consists of equilibria and periodic orbits only, and $NW(f)\neq\emptyset$
according to Corollary \ref{corollary:Nonzero_zero}.

To prove Theorem \ref{Th:NW(f)} it suffices to construct an algorithm that
takes as input $(k, f)$ and returns a set in $\mathcal{K}(\mathbb{D})$ such
that the Hausdorff distance between the output set and $NW(f)$ is less than
$1/k$, for every $k\in\mathbb{N}$ and every $f\in SS_{2}$. The construction
concept is intuitive and not entirely new (see, for example, \cite{DFJ01}); it
can be outlined as follows: first, cover the compact set $\mathbb{D}$ with a
finite number of square \textquotedblleft pixels;" second, use a rigorous
numerical method to compute the (flow) images of all pixels after some time
$T$, and take the union $\Omega_{T}$ of all images of pixels as a candidate
for an approximation to $NW(f)$; third, test whether $\Omega_{T}$ is an
over-approximation of $NW(f)$ within the desired accuracy. If the test fails,
increase $T$, and use a finer lattice of square pixels when numerically
approximating the flow of (\ref{ODE_Main}) after time $T$. Similar simulations
using time $-T$ are run in parallel to find repellers.

The novel and intricate components of the algorithm are where the saddle
points are dealt with, and the search for a time $T$ such that the Hausdorff
distance between $\Omega_{T}$ and $NW(f)$ is less than $1/k$ with $(k, f)$
being the input to the algorithm. Two comments seem in order. The problem with
a saddle point is that it may take an arbitrarily long time for the flow
starting at some point near but not on the stable manifold of the saddle to
eventually move away from the saddle. This undesirable behavior is dealt with
by transforming the original flow near a saddle to a linear flow using a
computable version of Hartman-Grobman's theorem (\cite{GZD12}). The time
needed for the linear flow to go through a small neighborhood can be
explicitly calculated (see section 8.2 for details). Another key feature of
the algorithm is that for every $k\in\mathbb{N}$, the algorithm computes a
uniform time bound $T$ such that $d_{H}(\Omega_{T}, NW(f))< 1/k$, where each
connected component of $\Omega_{T}$ is in a donut shape containing at least
one periodic orbit if it doesn't contain any equilibrium. This feature is
crucial for finding the number and positions of the periodic orbits of the
system (\ref{ODE_Main}). A coloring program is constructed for checking
whether $T$ is a good enough time (see section 8.3 for details).

We proceed to prove Theorem B once Theorem \ref{Th:NW(f)} is proved. The idea
is to use the coloring algorithm to find a cross-section for each connected
component of $\Omega_{T}$ and then compute the Poincar\'{e} maps. By finding
the number of fixed points of each Poincar\'{e} map, via a zero-finding algorithm, we will be able to count the
total number of periodic orbits of (\ref{ODE_Main}).

The remaining sections are devoted to the proof of Theorem B. Section
\ref{Sec:Zeros} explains how the algorithm for computing the number of zeros of a
function works. Section \ref{Sec:Simulation} explains how we can numerically compute the solution of the ODE (\ref{ODE_Main}) at
a time $T$ in a rigorous manner. This result is, in a sense, not surprising, but we present the
details so that we can adapt them later to the more subtle case when saddle
points are present in (\ref{ODE_Main}). In Section \ref{Sec:Proof} we prove
Theorem \ref{Th:NW(f)}; we start with the simpler case where (\ref{ODE_Main})
has no saddle points, and then proceed to the general case. In Section
\ref{Sec:Limit_cycles}, we explain in more details the coloring algorithm used
to determine the accuracy of approximations to periodic orbits as well to
define cross-sections. In Section \ref{Sec:Poincare} we show how to compute
the Poincar\'{e} maps and their derivatives defined on those cross-sections,
and finally prove Theorem B in Section \ref{Sec:ThmBProofFinal}.

\section{Computing the number of zeros of a function\label{Sec:Zeros}}

For every $f\in SS_{2}$, let $Zero(f)=\{ x\in\mathbb{D}: f(x)=0\}$, and let
$\#(f)=|Zero(f)|$ be the number of zeros of $f$. Then $\#(f)\neq0$ according
to Corollary \ref{corollary:Nonzero_zero}. It is proved in \cite{GZ21a} that
there exists an algorithm taking as input $(k, f)$, $k\geq1$, and outputting
$\#(f)$ and a non-empty set $C$ in $\mathcal{K}(\mathbb{D})$ such that
$Zero(f)\subset C$ and $d_{H}(C, Zero(f))<\frac{1}{k}$. The construction of
the algorithm relies on the fact that the equilibrium point(s) of the system
$x^{\prime}=f(x)$ are all hyperbolic for every $f\in SS_{2}$, and thus $f$ is
invertible in some neighborhood of each equilibrium point.

The following is a brief sketch of the algorithm: Start from $l=k$ and cover
$\mathbb{D}$ with side-length $1/l$ square pixels. For each pixel $s$, compute
$d(0, f(s))$ and $\min_{x\in s}\| Df(x)\|$, increase $l$ if necessary until
either $d(0, f(s))> 2^{-l}$ or $\min_{x\in s}\|Df(x)\| > 2^{-l}$ after
finitely many increments. If $d(0, f(s))> 2^{-l}$, then there is no
equilibrium inside $s$; if $\min_{x\in s}\|Df(x)\| > 2^{-l}$, then use the
elements from the proof of the inverse function theorem to either locate the
squares contained in $s$ such that each hosts a unique equilibrium or else
increase $l$ and repeat the process on smaller pixels contained in $s$. The
full details of the construction can be found in \cite{GZ21a}.

The algorithm also works for any function defined on a computable compact
subset of $\mathbb{R}^{2}$ that has finitely many zeros if all of them are invertible.

We mention in passing that the algorithm is fully automated in comparison with
most familiar root-finding numerical algorithms - Newton's method, Secant
method, etc - in the sense that no extra \textit{ad hoc} information or
analysis - such as a good initial guess or a priori knowledge on existence of
zeros or requiring further properties of the given function in order to
distinguish nearby zeros - is needed.

\section{Discrete simulation of planar dynamics\label{Sec:Simulation}}

We begin with some preliminary notions.

\begin{definition}
The $\delta$-grid in $\mathbb{R}^{d}$ is the set $G=(\mathbb{Z}\delta)^{d}$,
where $\delta\mathbb{Z}=\{\delta z\in\mathbb{R}:z\in\mathbb{Z}\}=\{\ldots
,-2\delta,-\delta,0,\delta,2\delta,\ldots\}$.
\end{definition}

An initial-value problem $x^{\prime}=f(x)$, $x(t_{0})=x_{0}$, is often solved
numerically using, for example, Euler's method, which can be described as
follows: select an initial point $y_{0}\simeq x_{0}$ (ideally $y_{0}=x_{0}$),
choose $h>0$ (the stepsize), and set%
\begin{equation}
y_{n+1}=y_{n}+hf(y_{n})\text{, \quad}n=0,1,2,\ldots\label{Eq:Euler}%
\end{equation}
If $f$ has bounded partial derivatives on $\mathbb{D}$ and hence a bounded
jacobian, then $f$ satisfies a Lipschitz condition%
\[
\left\Vert f(x)-f(y)\right\Vert \leq L\left\Vert x-y\right\Vert
\]
for all $x,y\in\mathbb{D}$, where $L$ is a Lipschitz constant, which can be
chosen as $L=\max_{x\in\mathbb{D}}\left\Vert Df(x)\right\Vert $ (see
e.g.~\cite[p. 26]{BR89}). Assume that $x(t)$ is defined for all $t\in\lbrack
t_{0},b]$ for some $b>t_{0}$. Let $x_{n}=x(t_{0}+nh)$, $n=0,1,2,\ldots$, and
let $\rho$ be a rounding error bound when computing $y_{n+1}$ using
(\ref{Eq:Euler}). Then the (global truncation) error of Euler's method is
bounded by the following formula (see \cite[p. 350]{Atk89} for one-dimensional
case and Appendix \ref{Sec:AppendixEuler} for the two-dimensional case):%
\begin{equation}
\left\Vert y_{n}-x_{n}\right\Vert \leq e^{(b-t_{0})L}\left\Vert y_{0}%
-x_{0}\right\Vert +\left(  \frac{e^{(b-t_{0})L}-1}{L}\right)  \left(  \frac
{h}{2}E+\frac{\rho}{h}\right)  \label{Eq:ErrorEuler1}%
\end{equation}
for all $n=0,1,2,\ldots$ satisfying $t_{0}\leq t_{0}+nh\leq b$, where $E$ is a
bound for $\left\Vert x^{\prime\prime}(t)\right\Vert $. Since
\[
\left\Vert x^{\prime\prime}(t)\right\Vert =\left\Vert \left(  f(x(t))\right)
^{\prime}\right\Vert =\left\Vert \left(  \frac{\partial f_{1}}{\partial x_{1}%
}f_{1}(x)+\frac{\partial f_{1}}{\partial x_{2}}f_{2}(x),\frac{\partial f_{2}%
}{\partial x_{1}}f_{1}(x)+\frac{\partial f_{2}}{\partial x_{2}}f_{2}%
(x)\right)  \right\Vert
\]
it follows that (\ref{Eq:ErrorEuler1}) can be rewritten as%
\begin{equation}
\left\Vert y_{n}-x_{n}\right\Vert \leq e^{(b-t_{0})L}\left\Vert y_{0}%
-x_{0}\right\Vert +\left(  \frac{e^{(b-t_{0})L}-1}{L}\right)  \left(
hM^{2}+\frac{\rho}{h}\right)  \label{Eq:ErrorEuler}%
\end{equation}
where $M=\max_{x\in\mathbb{D}}(\left\Vert f(x)\right\Vert ,\left\Vert
Df(x)\right\Vert )$, which is computable from (a $C^{1}$-name of) $f$.

The global truncation error bound (\ref{Eq:ErrorEuler}) depends on $h$, $\rho
$, $b$, $L$, and $M$, where $L$ and $M$ are computable from $f$. We may assume
that $L, M\geq1$. It is possible to make the error smaller than any given
$\epsilon>0$ over any time interval $[t_{0}, T]$, as long as the solution
$x(t)$ is defined, by selecting appropriate values for $h$ and $\rho$. In
other words, we can choose $h$ and $\rho$ such that $\left\Vert y_{n}%
-x_{n}\right\Vert \leq\epsilon$ for all $n=0,1,2,\ldots$ satisfying $t_{0}\leq
t_{0}+nh\leq T$ (we assume that $y_{0}$ is obtained from $x_{0}$ by rounding
it with error bounded by $\rho$). This can be achieved, for example, by
requiring that each term of the sum in the right-hand side of
(\ref{Eq:ErrorEuler}) is bounded by $\epsilon/2$. For the first term,
\[
e^{TL}\rho\leq\frac{\epsilon}{2}\text{\quad}\Longrightarrow\text{\quad}%
\rho<\frac{\epsilon}{2e^{TL}}.
\]
For the second term,
\begin{align}
\left(  \frac{e^{TL}-1}{L}\right)  \left(  hM^{2}+\frac{\rho}{h}\right)   &
\leq e^{TL}\left(  hM^{2}+\frac{\rho}{h}\right) \nonumber\\
&  =e^{TL}hM^{2}+e^{TL}\frac{\rho}{h}. \label{Eq:bound_1}%
\end{align}
If we require both terms on the right hand side of (\ref{Eq:bound_1}) are
bounded by $\epsilon/4$, we obtain the desired estimate. For the first term of
(\ref{Eq:bound_1}),%
\[
e^{TL}hM^{2}\leq\frac{\epsilon}{4}\text{\quad}\Longrightarrow\text{\quad}%
h\leq\frac{\epsilon}{4e^{TL}M^{2}}.
\]
Now we fix some $h$ satisfying the inequality above. We can then derive a
desirable value of $\rho$ by bounding the second term of (\ref{Eq:bound_1})
with $\epsilon/4$
\[
e^{TL}\frac{\rho}{h}\leq\frac{\epsilon}{4}\text{\quad}\Longrightarrow
\text{\quad}\rho\leq\frac{\epsilon h}{4e^{TL}}.
\]

The above analysis is summarized in the following theorem:

\begin{theorem}
Let $\epsilon>0$, $T>0$ be given and assume that the solution $x(t)$ of
(\ref{ODE_Main}) with initial condition $x(0)=x_{0}$ is defined in the
interior of $\mathbb{D}$ for all $t\in\lbrack0,T]$. By selecting a stepsize $h$ satisfying%
\begin{equation}
h\leq\frac{\epsilon}{4e^{TL}M^{2}} \label{Eq:h}%
\end{equation}
and then using a rounding error $\rho$ bounded by
\begin{equation}
\rho\leq\min\left(  \frac{\epsilon h}{4e^{TL}},\frac{\epsilon}{2e^{TL}%
}\right)  \label{Eq:rho}%
\end{equation}
the approximations $y_n$ generated by Euler's mathod have the property that $\left\Vert y_{n}-x_{n}\right\Vert \leq\epsilon$ for all
$n=0,1,2,\ldots$ satisfying $0\leq nh\leq T$.
\end{theorem}

We now define a computable function \texttt{ChooseParameters} as follows: it
takes as input $(f,\epsilon,T)$, where $\epsilon>0$ and $T>0$ are rational
numbers, and returns in finite time as output $(h,\rho,n_{T})$, where $h$ and
$\rho$ are rational numbers satisfying (\ref{Eq:h}) and (\ref{Eq:rho}),
$n_{T}\in\mathbb{N}$ and $T=n_{T}h$. Then, using the function
\texttt{ChooseParameters} as a subroutine and Euler's method, we can devise a
new computable function \texttt{TimeEvolution} that receives as input some
compact set $D\subseteq\mathbb{D}$, $f$ and some rational numbers
$0<\epsilon<1$ and $T>0$, with the following properties (usually we omit the
explicit dependence on $f$ and write \texttt{TimeEvolution}$(D,\epsilon,T)$
instead of \texttt{TimeEvolution}$(f,D,\epsilon,T)$ if that is clear from the context):

\begin{itemize}
\item $\phi_{T}(D)\subseteq$ \texttt{TimeEvolution}$(D,\epsilon,T)\subseteq
\phi_{T}(D)+\epsilon B(0,1)$. In particular this implies that $d_{H}(\phi
_{T}(D),$ \texttt{TimeEvolution}$(D,\epsilon,T))\leq\epsilon$.

\item The computation of \texttt{TimeEvolution} with input $(D,\epsilon,T)$
halts in finite time$.$
\end{itemize}

The algorithm that computes \texttt{TimeEvolution} is designed as follows:

\begin{enumerate}
\item Compute \texttt{ChooseParameters}$(\epsilon/4,T)$, and obtain the
corresponding values $h,\rho,n_{T}$. Next create a $(\rho/2)$-grid over
$\mathbb{D}$ (technically, create a $(\rho/2)$-grid over $\mathbb{R}^{2}$ and
then intersect it with $\mathbb{D}$).

\item For each (rational) point $p$ of the $\frac{\rho}{2}$-grid over
$\mathbb{D}$, decide whether $d(p,D)\leq\rho/2$ or $d(p,D)\geq\rho/3$. Let
$D_{1}$ be the set of all $\frac{\rho}{2}$-grid points $p$ with the test
result $d(p,D)\leq\rho/2$.

\item Apply Euler's method to all points in $D_{1}$, using $h$ as the timestep
and $\rho$ as the rounding error. Let $\tilde{p}_{n}$ be the $n$th iterate
obtained by applying Euler's method for the IVP (\ref{ODE_Main}), $x(0)=p$,
with these parameters. Then output%
\[%
{\displaystyle\bigcup\limits_{p\in D_{1}}}
B(\tilde{p}_{n_{T}},\epsilon/2).
\]

\end{enumerate}

It is readily seen that all three steps can be executed in finite time, and
thus the computation of \texttt{TimeEvolution} with input $(D,\epsilon,T)$
halts in finite time. The second property is satisfied. It remains to show
that the function \texttt{TimeEvolution} also has the first property. Since a
$\rho/2$-grid is used, every point $q\in D$ is within distance $\leq\rho/2$ of
a $\rho/2$-grid point $p$, which implies that $p\in D_{1}$. Furthermore, it
follows from the definition of the function \texttt{ChooseParameters}%
$(\epsilon/2,T)$ that%
\begin{equation}
\left\Vert \tilde{p}_{n_{T}}-\phi_{T}(p)\right\Vert \leq\epsilon/4
\quad(\mbox{since $\phi_{0}(p)=p$}). \label{eq:bound1Euler}%
\end{equation}
On the other hand, it follows from Lemma \ref{Lemma:Perturbation} and
(\ref{Eq:rho}) that%
\[
\left\Vert \phi_{T}(p)-\phi_{T}(q)\right\Vert \leq\frac{\rho}{2}e^{LT},
\quad\left\Vert \phi_{T}(p)-\phi_{T}(q)\right\Vert \leq\frac{(\epsilon
/4)h}{4e^{TL}}e^{LT}\leq\frac{\epsilon}{16}.
\]
Combine the last inequalities together with (\ref{eq:bound1Euler}) yields%
\begin{equation}
\left\Vert \tilde{p}_{n_{T}}-\phi_{T}(q)\right\Vert \leq\left\Vert \tilde
{p}_{n_{T}}-\phi_{T}(p)\right\Vert +\left\Vert \phi_{T}(p)-\phi_{T}%
(q)\right\Vert \leq\frac{\epsilon}{4} + \frac{\epsilon}{16}=5\epsilon/16.
\label{eq:bound2Euler}%
\end{equation}
Hence $\phi_{T}(D)\subseteq$ \texttt{TimeEvolution}$(D,\epsilon,T)$. To verify
that \texttt{TimeEvolution}$(D,\epsilon,T)\subseteq D+\epsilon B(0,1)$, it
suffices to show that if $r\in\cup_{p\in D_{1}}B(\tilde{p}_{n_{T}}%
,\epsilon/2)$, then there exists a point $q\in D$ such that $\left\Vert
r-\phi_{T}(q)\right\Vert \leq\epsilon$. Since $r\in\cup_{p\in D_{1}}%
B(\tilde{p}_{n_{T}},\epsilon/2)$, it follows that there is a point $p\in
D_{1}$ such that%
\[
\left\Vert r-\tilde{p}_{n_{T}}\right\Vert \leq\epsilon/2.
\]
Therefore, there must be a point $q\in D$ within distance $\leq\rho/2$ from
the point $p$. In particular, (\ref{eq:bound2Euler}) must hold. This, together
with the above inequality, shows that%
\[
\left\Vert r-\phi_{T}(q)\right\Vert \leq\left\Vert r-\tilde{p}_{n_{T}%
}\right\Vert +\left\Vert \tilde{p}_{n_{T}}-\phi_{T}(q)\right\Vert \leq
\epsilon/2+5\epsilon/16<\epsilon.
\]

We end this section by defining a computable function
\texttt{HasInvariantSubset} that receives as input some compact set
$D\subseteq\mathbb{D}$, a function $f$, and some rational numbers
$0<\epsilon<1$ and $T>0$. If this function returns 1, then the set $D$ is
guaranteed to be invariant in the sense that $\phi_{t}(D)\subseteq D$ for all
$t\geq T$. If it returns 0, then $\mathbb{D}$ may or may not contain invariant
subsets. We can compute \texttt{HasInvariantSubset}$(f,D,\epsilon,T)$ (which
we will denote simply as \texttt{HasInvariantSubset}$(D,\epsilon,T)$ if $f$ is
clear from the context) as follows:

\begin{enumerate}
\item Compute rationals $t_{0},t_{1},\ldots,t_{N}$ satisfying $0=t_{0}%
<t_{1}<t_{2}<\ldots<t_{N}=T$ and $\left\vert t_{i+1}-t_{i}\right\vert
\leq\epsilon/(4M)$. (Recall that $M=\max_{x\in\mathbb{D}}(\left\Vert
f(x)\right\Vert ,\left\Vert Df(x)\right\Vert )$.)

\item Let $D_{i}=$ \texttt{TimeEvolution}$(D,\epsilon/4,T+t_{i})$ be a union
of finitely many balls with rational centers and radii for each $i=0,1,\ldots
,N$.

\item Compute an over-approximation $A$, consisting on the union of finitely
many balls with rational centers and radii, of $\overline{\mathbb{D}-D}$ with
accuracy bounded by $\epsilon/4.$

\item For every $i=0,1,\ldots, N$, test if $D_{i}\cap(A+B(0,\epsilon
/4))=\varnothing$. If the test succeeds, return 1; otherwise return 0.
\end{enumerate}

Steps 1--3 can be easily computed, as well as $A+B(0,\epsilon/4)$ (it suffices
to increase the radius of each ball defining $A$ by $\epsilon/4$). We can also
determine in finite time whether or not $D_{i}\cap(A+B(0,\epsilon/4))$ is
empty since $D_{i}$ and $A+B(0,\epsilon/4)$ are formed by finitely many
rational balls. Now suppose step 4 returns 1. We show that $\phi_{t}(D)\subset
D$ for all $t\geq T$. We first note that $\phi_{T+t_{i}}(D)\subseteq D_{i}$ by
definition of \texttt{TimeEvolution} and $D_{i}+B(0, \epsilon/4)\subset D$ by
assumption. Furthermore, for any $0\leq t\leq T$, there exists an i, $0\leq
i\leq N-1$, such that $t\in\lbrack t_{i}, t_{i+1}]$. Since $\left\vert
t-t_{i}\right\vert \leq\epsilon/(4M)$, it follows that $\left\Vert \phi
_{T+t}(x)-\phi_{T+t_{i}}(x)\right\Vert \leq\epsilon/4$ for every $x\in D$.
Hence, $\phi_{T+t}(D) \subseteq\phi_{T+t_{i}}(D) + B(0, \epsilon/4) \subseteq
D_{i} + B(0, \epsilon/4) \subseteq D$ for every $0\leq t\leq T$. In other
words, $\phi_{t}(D)\subseteq D$ for all $T\leq t\leq2T$. Now it follows from
$\phi_{T}(D)\subseteq D$ that $\phi_{2T+t}(D) \subseteq\phi_{T+t}(D)\subseteq
D$ for all $0\leq t \leq T$, or $\phi_{t}(D)\subseteq D$ for all $2T\leq
t\leq3T$. Continuing the process inductively, we conclude that $\phi
_{t}(D)\subseteq D$ for all $t\geq T$. A final note. If $D$ includes in its
interior an attracting hyperbolic point or attracting periodic orbit and $D$
is inside the basin of attraction of this attractor, then
\texttt{HasInvariantSubset}$(D,\epsilon,T)$ will return 1 for sufficiently
large $T$ and small enough $\epsilon$ according to Propositions
\ref{Prop:ConvEquilibrium} and \ref{Prop:ConvPeriodic}.

\section{Proof that the non-wandering set is computable\label{Sec:Proof}}

In this section, we prove Theorem \ref{Th:NW(f)}: The operator $\Psi: SS_{2}
\to\mathcal{K}(\mathbb{D})$, $f\mapsto NW(f)$ of (\ref{ODE_Main}), is computable.

To show $NW(f)$ is computable, it suffices to construct an algorithm that
takes as input $(k, f)$ and returns as output a compact non-empty set
$NW_{k}(f)$ such that

\begin{description}
\item[(I)] $NW(f)\subseteq NW_{k}(f)\subseteq\mathbb{D}$.

\item[(II)] $d_{H}(NW(f), NW_{k}(f))\leq1/k$.
\end{description}

The idea underlying the construction is to simulate (\ref{ODE_Main})
numerically using Euler's method with rigorous bounds on the error generated
from approximation as described in the previous section; we will show how to
establish such rigorous error bounds. By using increasingly accurate
approximations, we are able to uniformly approach $NW(f)$ (cf.~Theorem
\ref{Cor:uniform_convergence_full}) with arbitrary precision.

Recall that since $f\in SS_{2}$, $NW(f)$ consists only of $Zero(f)$ and
$Per(f)$, a finite set of equilibrium points and a finite set of periodic
orbits, respectively; moreover, $Zero(f)\neq\emptyset$ but $Per(f)$ might be
empty. An algorithm for computing $Zero(f)$ has been outlined in section
\ref{Sec:Zeros}; the details of that algorithm are presented in \cite{GZ21a}.

\subsection{A first approach to the problem: no saddle
point\label{Sec:1stApproach}}

The case where (\ref{ODE_Main})\ has no saddle point is less complicated
because in this case, we can make use of Theorem
\ref{Cor:uniform_convergence_full} to find approximations to $NW(f)$ by
computing $\phi_{t}(\mathbb{D})$. But how can we actually find such a time
$T_{k}$ given in Theorem \ref{Cor:uniform_convergence_full}? And how do we
decide if a point $x$ is already in $\mathcal{N}_{1/k}(NW(f))$ without knowing
$NW(f)$? (After all our goal is to locate $NW(f)$.) One simple but essential
fact about $NW(f)$ is that it is the \textquotedblleft
minimal\textquotedblright\ time invariant set of the system (\ref{ODE_Main}).
This observation leads to the following tactics for developing a working
algorithm: numerically simulate $\phi_{T}(\mathbb{D})$ using
$\mathtt{TimeEvolution}(\mathbb{D}, \epsilon, T)$; note that
$\mathtt{TimeEvolution}(\mathbb{D}, \epsilon, T)$ is valid for all values of
$T$ as explained at the end of section \ref{Sec:Simulation}, $NW(f)\subset
\mathtt{TimeEvolution}(\mathbb{D}, \epsilon, T)$, and by its definition,
\begin{equation}
\mathtt{TimeEvolution}(\mathbb{D}, \epsilon,T)=%
{\displaystyle\bigcup\limits_{i\in I}}
B(p_{i},\epsilon/2) \label{Eq:Time_evolution}%
\end{equation}
where $I$ is a finite set and $p_{i}\in\mathbb{Q}$. By making use of the
pixels $B(p_{i},\epsilon/2)$ it is possible to compute the invariant sets of
the simulation and select those invariant sets which are minimal, and thus
find approximations to $NW(f)$ with arbitrary precision by using increasingly
accurate simulations -- smaller $\epsilon$ and larger $T$.

In the following we present the algorithm first, and then we show that the
algorithm returns the desired output in finite time.

The algorithm to compute $NW(f)$ runs as follows: taking as input $(k, f)$,
where $1/k$ sets the error bound for the output approximation, $k\in
\mathbb{N}\setminus\{ 0\}$, do:

\begin{enumerate}
\item Let $\epsilon=1/(2k)$ and $T=1.$

\item Compute the set $Zero(f)$ with precision $\epsilon$ using the algorithm
described in \cite{GZ21a}, obtaining an over-approximation $Zero_{\epsilon}(f)$
of $Zero(f)$ satisfying $d_{H}(Zero_{\epsilon}(f), Zero(f))\leq\epsilon$.

\item Compute $\mathtt{TimeEvolution}(\mathbb{D},\epsilon,T)$ and obtain a
finite set $I$ for which (\ref{Eq:Time_evolution}) holds true.

\item Take $P=\varnothing.$

\item For each $J\subseteq I$ check whether $\bigcup_{j\in J}B(p_{j}%
,\epsilon/2)$ is a connected set and whether \texttt{HasInvariantSubset}%
$(D,\epsilon,T)=1$. If both conditions hold, take $P=P\cup\{J\}$.

\item (The minimum principle) For each $J_{1},J_{2}\in P$, test whether
$J_{1}\subsetneqq J_{2}$. If the test succeeds, then take $P=P\backslash
\{J_{2}\}.$

\item For each $J\in P$, test whether $\bigcup_{j\in J}B(p_{j},\epsilon)$ is
contained in a connected component of $Zero_{\epsilon}(f)$. If the test
succeeds, then take $P=P\backslash\{J\}.$ If the test fails, but
$\bigcup_{j\in J}B(p_{j},\epsilon)\bigcap Zero_{\epsilon}(f)\neq\emptyset$,
then take $\epsilon:=\epsilon/2$, $T:=2T$, and go to step 2. For each $J\in P$, test if $\cup_{j\in J}B(p_{j},\epsilon
)\cap\partial\mathbb{D}\neq\varnothing$. If the test succeeds, then take
$\epsilon:=\epsilon/2$, $T:=2T$.

\item Do steps 2--7 to the ODE $x^{\prime}=-f(x)$ defined over $\mathbb{D},$
i.e.\ by taking the transformation $t\mapsto-t$ which reverses time in
(\ref{ODE_Main}), obtaining a set $\tilde{P}$ similar to the set $P$ of step 7.

\item For each set $J\in P\cup\tilde{P}$, use the algorithm described in
Section \ref{Sec:Limit_cycles} to check whether $\bigcup_{j\in J}%
B(p_{j},\epsilon)$ has a doughnut shape with cross-sections bounded by $1/k$.
If this check fails, then take $\epsilon:=\epsilon/2$, $T:=2T$, and go to step 2.

\item Check whether $\mathtt{TimeEvolution}(\mathbb{D}-(Zero_{\epsilon
}(f)\bigcup~\bigcup_{J\in P\bigcup\tilde{P}}(\bigcup_{j\in J}B(p_{j}%
,\epsilon)),\epsilon/3,T))\subseteq Zero_{\epsilon}(f)\bigcup~\bigcup_{J\in
P\bigcup\tilde{P}}(\cup_{j\in J}B(p_{j},\epsilon))$. If this test fails, then
take $\epsilon:=\epsilon/2$, $T:=2T$, and go to step 2.

\item Switch the dynamics to the ODE $x^{\prime}=-f(x)$ and test whether
$\mathtt{TimeEvolution}(\mathbb{D}-(Zero_{\epsilon}(f)\bigcup~\bigcup_{J\in
P\bigcup\tilde{P}}(\bigcup_{j\in J}B(p_{j},\epsilon)),\epsilon/3,T))\subseteq
Zero_{\epsilon}(f)\bigcup~\bigcup_{J\in P\bigcup\tilde{P}}(\bigcup_{j\in
J}B(p_{j},\epsilon))$ for this ODE. If this test fails, then take
$\epsilon:=\epsilon/2$, $T:=2T$, and go to step 2.

\item Output $Zero_{\epsilon}(f)\cup\bigcup_{J\in P\cup\tilde{P}}\cup_{j\in
J}B(p_{j},\epsilon)$ ($=NW_{k}(f)$).
\end{enumerate}

It is perhaps time to explain a bit of roles played by the steps in the
algorithm. Step 2 supplies a subprogram for locating the set of equilibrium
point(s) with arbitrary precision whenever the need arises; in particular, it
detects whether a (time) invariant set contains an equilibrium.
Step 5 identifies the connected invariant sets of the simulation
$\mathtt{TimeEvolution}(\mathbb{D},\epsilon,T)$, which serve as possible
candidates for approximations to $NW(f)$ (see discussion at the end of Section
\ref{Sec:Simulation}). In order to get \textquotedblleft good" candidates for
$Per(f)$, those who contain the equilibrium point(s) or who are not
\textit{minimum} are discarded using steps 6 and 7. After step 8, every
remaining candidate -- if there is any -- is \textit{minimum}, and each of
them contains at least one periodic orbit (see corollary
\ref{corollary:Nonzero_zero}). Afterwards, step 9 checks whether each
remaining candidate is a good enough approximation to those periodic orbit(s)
it contains. Finally, steps 10 and 11 check if there are any periodic orbits
which have not been counted in the current round of simulation.

\begin{figure}
\begin{center}
\includegraphics[width=5cm]{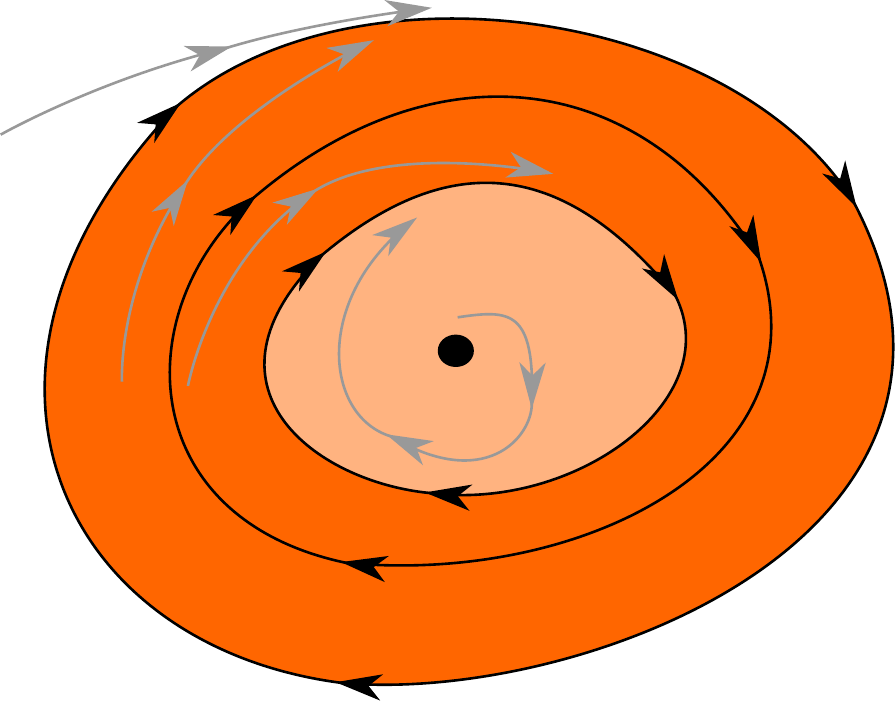}
\caption{Dynamical systems with nested periodic orbits.}
\label{fig:nested_orbits}
\end{center}
\end{figure}

Here is an example demonstrating why only the \textit{minimum} invariant sets
are qualified to serve as candidates: Consider Fig.~\ref{fig:nested_orbits}.
In this figure, we have a repelling equilibrium point which is surrounded by
(in order) an attracting, a repelling, and an attracting periodic orbit which
are nested together. It is not hard to see that if $U$ denotes the region
\textquotedblleft inside\textquotedblright\ the outer periodic orbit
(including this orbit), then $\phi_{t}(U)=U$ for all $t\geq0$. Indeed, it is
readily seen that $\phi_{t}(U)\subseteq U$. For the reverse inclusion, let
$x\in U$ and $t\geq0$ be given. Then $\phi_{-t}(x)\in U$, which implies that
$x=\phi_{t}(\phi_{-t}(x))\in\phi_{t}(U)$, i.e. $\phi_{t}(U)=U$. It is also
clear that $U$ has ``big" regions not contained in $NW(f)$ and it is not a
\textit{minimum} invariant set. Thus if $U$ is not eliminated as an candidate,
then step 9 will run into an infinite loop. The example reveals the reasoning
behind step 6.

Now we show that the algorithm returns the desired output in finite time. It
is clear that the algorithm runs through steps 1 -- 6 in finite time. Since
each equilibrium is either a sink or a source, for $\epsilon$ sufficiently
small, $Zero_{\epsilon}(f)$ becomes a set of disjoint $\epsilon$-balls with
each containing a unique equilibrium that attracts or repels all trajectories
inside the ball towards it or away from it. Furthermore, since there are only
finitely many periodic orbits and each periodic orbit is a compact subset of
$\mathbb{D}$, it follows that, after finitely many updates on $\epsilon$ and
$T$, for every $J\in P$ from step 6, either $\cup_{j\in J}B(p_{j},\epsilon)$
is contained in a connected component of $Zero_{\epsilon}(f)$ or $\cup_{j\in
J}B(p_{j},\epsilon)\cap Zero_{\epsilon}(f)=\emptyset$. In other words, step 7
as well step 8 completes its task in finite time. As output of steps 7 and 8,
if $P\cup\tilde{P}= \emptyset$, move to steps 10 and 11; otherwise, for every
$J\in P\cup\tilde{P}$, since $\phi_{T}(\cup_{j\in J}B(p_{j},\epsilon
))\subseteq\mathtt{TimeEvolution}(\cup_{j\in J}B(p_{j},\epsilon),\epsilon
/3,T)$ by the definition of the simulation $\mathtt{TimeEvolution}$, it
follows from step 5 that $\phi_{T}(\cup_{j\in J}B(p_{j},\epsilon))=\cup_{j\in
J}B(p_{j},\epsilon)$, which implies that $\cup_{j\in J}B(p_{j},\epsilon)$
contains at least one periodic orbit, since it does not contains any
equilibrium point. Now if $\cup_{j\in J}B(p_{j},\epsilon)$ passes the test in
step 9 for every $J\in P\cup\tilde{P}$ and, afterwards, $Zero_{\epsilon
}(f)\bigcup_{J\in P\cup\tilde{P}}\cup_{j\in J}B(p_{j}, \epsilon)$ passes the
tests in steps 10 and 11, then it is not hard to see that the output
$Zero_{\epsilon}(f)\cup\cup_{J\in P\cup\tilde{P}}\cup_{j\in J}B(p_{j},
\epsilon)$ satisfies requirements (I) and (II). But, the minimum nature of the
remaining invariant sets after step 6
together with Theorem \ref{Cor:uniform_convergence_full} ensure that all
invariant sets output by step 7 or step 8 will pass the three tests in steps
9, 10, and 11 after finitely many updates on $\epsilon$ and $T$.

We note that the precision $\epsilon$ and time $T$ are independent of each
other in Theorem \ref{Cor:uniform_convergence_full} but dependent in the
algorithm -- an update doubling $T$ and cutting $\epsilon$ in half. This
technicality is dealt with as follows: suppose that $n$ iterations have been
performed (i.e.~after $n$ updates of the $T$ variable and of $\epsilon$). Then
we have $T=2^{n}$ and $\epsilon=1/(2^{n+1}k)$. According to Propositions
\ref{Prop:ConvEquilibrium} and \ref{Prop:ConvPeriodic}, it is sufficient to
show that there is some $n_{0}\in\mathbb{N}$ such that for all $n\geq n_{0}$
one has%
\begin{align*}
K_{1}e^{-\alpha2^{n}}  &  <\frac{1/(2^{n+1}k)}{6}=\frac{1}{2^{n+2}3k}\text{
\ \ }\Rightarrow\\
-\alpha2^{n}  &  <-\ln\left(  2^{n+2}3kK_{1}\right)  \text{ \ \ }\Rightarrow\\
2^{n}  &  >\frac{\ln\left(  2^{n+2}3kK_{1}\right)  }{\alpha}\text{
\ \ }\Rightarrow\\
2^{n}  &  >\frac{n+2+\ln\left(  3kK_{1}\right)  }{\alpha}%
\end{align*}
Thus there is indeed an $n_{0}\in\mathbb{N}$, which depends on $k,K_{1}%
,\alpha>0$, such that the last inequality is true for all $n\geq n_{0}$.
Therefore, if we simultaneously double $T$ and halve $\epsilon$, we can be
sure that the condition $\phi_{T}(\cup_{j\in J}B(p_{j},\epsilon))\subseteq
\mathcal{N}_{\epsilon/6}(NW(f))$ will eventually hold after $n_{0}$ iterations
of the method.

We mention in passing that, in step 1, we could just have required that one
should compute a finite set $I$ and rationals $p_{i}$, $i\in I$ such that
(\ref{Eq:Time_evolution}) holds. However, computing $\mathtt{TimeEvolution}%
(K,\epsilon,T)$ automatically does it and may provide a computational gain.

\textit{A final note.} The algorithm uses internally uniform time bounds for
performing simulations $\mathtt{TimeEvolution}(\mathbb{D}, \epsilon, T)$,
computations, and tests; if a test fails, it restarts a new round by doubling
the time bound in the previous round (as well halving the error bound).


\subsection{The full picture\label{Sec:FullPicture}}

In Section \ref{Sec:1stApproach} we have assumed that the flow of
(\ref{ODE_Main}) does not have any saddle point. We now drop that assumption.
The major difference from the previous case is that Theorem
\ref{Cor:uniform_convergence_full} is no longer valid, as explained in Remark
\ref{Re:no_saddle}.
A fundamental implication of Theorem \ref{Cor:uniform_convergence_full}
from the view point of computation is that it provides a \textquotedblleft
uniform\textquotedblright\ time bound $T_{k}$ for approximating $NW(f)$, using
$\phi_{T_{k}}(\mathbb{D})$, globally on the entire phase space $\mathbb{D}$
with any precision $1/k$. The algorithm constructed in Section
\ref{Sec:1stApproach} relies on this uniform time bound coupled with the
minimum principle for performing numerical simulations $\mathtt{TimeEvolution}%
(\mathbb{D}, \epsilon, T)$ and computing invariant sets -- approximations to
$NW(f)$ -- of the simulations. The problem with a saddle point $x_{0}$ is that
when a trajectory passes near the stable manifold of $x_{0}$, the flow is to
approach $x_{0}$ but may move very slowly as it gets closer and closer to
$x_{0}$ (since $f$ has a zero at $x_{0}$) before it finally leaves the
vicinity of $x_{0}$ (or eventually converges to $x_{0}$, if the trajectory is
part of the stable manifold of $x_{0}$). In this case, the hope for having
uniform time bounds is slim, if not impossible.

We solve the problem by mending the previous algorithm so that only one time
unit is counted from the moment a trajectory entering a small neighborhood of
a saddle point until the moment it finally leaving the neighborhood, provided
that it indeed leaves the vicinity.
This is done by placing a \textquotedblleft black box\textquotedblright\ at
each saddle. Outside the black box, the algorithm runs as before; upon
entering the box the algorithm switches from simulating trajectories of
(\ref{ODE_Main}) to computing the trajectories of a linear system conjugated
to (\ref{ODE_Main}). One time unit is allotted to the work done in the black
box. The idea is motivated by Hartman-Grobman's Theorem and a computable
version of it (see e.g. \cite[p. 127]{Per01} for its classical statement; and
\cite{GZD12} for its computable version). For completeness, we state the
computable version here: let $\mathcal{F}$ be the set of all functions $f\in
C^{1}(\mathbb{R}^{n};\mathbb{R}^{n})$ such that $0$ is an hyperbolic zero of
$f$ i.e. $f(0)=0$ and $Df(0)$ only has eigenvalues with nonzero real part; let
$\mathcal{O}$ be the set of all open subsets of $\mathbb{R}^{n}$ containing
the origin of $\mathbb{R}^{n}$; and let $\mathcal{I}$ be the set of all open
intervals of $\mathbb{R}$ containing zero.

\begin{theorem}
\label{Prop:HartmanGrobman} There is a computable map $\Theta:\mathcal{F}%
\rightarrow\mathcal{O}\times\mathcal{O}\times C(\mathbb{R}^{n};\mathbb{R}%
)\times C(\mathbb{R}^{n};\mathbb{R}^{n})$ such that for any $f\in\mathcal{F}$,
$f\mapsto(U,V,\mu,H)$, where

\begin{itemize}
\item[(a)] $H:U\to V$ is a homeomorphism\,;

\item[(b)] the unique solution $x(t,\tilde{x})=x(\tilde{x})(t)$ to the initial
value problem $\dot{x}=f(x)$ and $x(0)=\tilde{x}$ is defined on $(-\mu
(\tilde{x}),\mu(\tilde{x}))\times U$; moreover, $x(t,\tilde{x})\in U$ for all
$\tilde{x}\in U$ and $-\mu(\tilde{x})<t<\mu(\tilde{x})$\thinspace;

\item[(c)] $H(x(t,\tilde{x}))=e^{Df(0)t}H(\tilde{x})$ for all $\tilde{x}\in U$
and $-\mu(\tilde{x})<t<\mu(\tilde{x})$\thinspace.
\end{itemize}
\end{theorem}

Recall that for any $\tilde{x}\in\mathbb{R}^{n}$, $e^{Df(0)t}\tilde{x}$ is the
solution to the linear problem $\dot{x}=Df(0)x$, $x(0)=\tilde{x}$. So the
theorem shows that the homeomorphism $H$, computable from $f$, maps
trajectories of the nonlinear problem $\dot{x}=f(x)$ onto trajectories of the
linear problem $\dot{x}=Df(0)x$, near the origin, which is a hyperbolic
equilibrium point. In other words, $H$ is a conjugacy between the linear and
nonlinear trajectories near the origin. Note that the theorem holds true for
any hyperbolic equilibrium point, and the origin is used just for convenience
(see \cite{GZD12} for details).

The linear system in a neighborhood of the origin, $\dot{x}=Df(0)x$ and
$x(0)=\tilde{x}$, can be solved explicitly with the solution $e^{Df(0)t}%
\tilde{x}$. Moreover, the solution is computable from (a $C^{1}$-name of) $f$
because $Df(0)$ is computable from $f$ and $e^{Df(0)}$ is computable from
$Df(0)$ (see, for example, \cite{WZ07}).

Thus Theorem \ref{Prop:HartmanGrobman} provides us a computational tool to
deal with the problem of lacking uniform time bounds near saddles for
simulations $\mathtt{TimeEvolution}$ as outlined as follows: for a saddle of
(\ref{ODE_Main}), use a neighborhood provided by Theorem
\ref{Prop:HartmanGrobman} as an oracle or a black box; once a trajectory
generated by the simulation $\mathtt{TimeEvolution}$ enters the black box,
$\mathtt{TimeEvolution}$ sits in idle and waits for an answer from the black
box; inside the black box, the trajectory of the linear system is computed
starting from a point on the simulated trajectory that enters the black box;
then the black box supplies an answer to $\mathtt{TimeEvolution}$ with a
desirable precision and $\mathtt{TimeEvolution}$ restarts working after
receiving the answer. The waiting time is counted as one time unit. It is not
hard to see that the strategy solves the problem of computing the flow of
(\ref{ODE_Main}) near a saddle point using the simulation
$\mathtt{TimeEvolution}$. Of course, for the strategy to work, we need to show
that the black box is able to supply an answer with required accuracy to
$\mathtt{TimeEvolution}$ for any simulated trajectory that enters it. The
proof is given below.

Apply Theorem \ref{Prop:HartmanGrobman} to each saddle point $x_{i}$ of
(\ref{ODE_Main}) with the homeomorphism $H_{i}: U_{i} \rightarrow V_{i}$,
where $x_{i}$ is the unique equilibrium contained in $U_{i}$. Compute a
rational number $\varepsilon_{i}$ such that $B(H_{i}(x_{i}), \varepsilon
_{i})\subseteq V_{i}$ for each $i$ and set $\varepsilon=\min_{i}%
\varepsilon_{i}$ (recall that $B(H_{i}(x_{i}), \varepsilon_{i})$ is the closed
ball having center $H_{i}(x_{i})$ and radius $\varepsilon_{i}$ with the
max-norm, and there are finitely many saddles). From now on we will
concentrate on the saddle point $x_{i}$; the other saddles can be dealt with
similarly. For simplicity, we assume that $x_{i}=0$.

\begin{figure}[ptb]
\begin{center}
\includegraphics{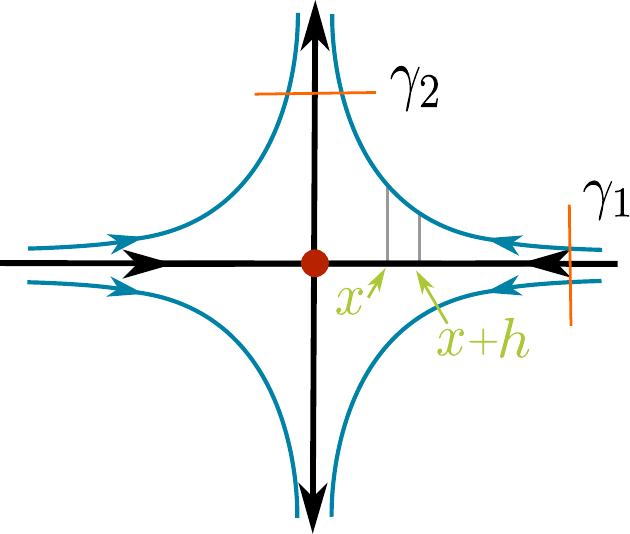}
\end{center}
\caption{Central idea used to determine the time needed to pass by a saddle.}%
\label{fig:saddle_time}%
\end{figure}

Referring to Theorem \ref{Prop:HartmanGrobman}, assume that $\lambda$ and
$\mu$ are eigenvalues of $Df(0)$ with $\lambda< 0 < \mu$. Since the
eigenspaces associated to these eigenvalues have known dimension (=1), and are
the kernel of $Df(0)-\lambda I$ and $Df(0)-\mu I$, we can compute (see
\cite[Theorem 11 - a), c)]{ZB04}) eigenvectors $v_{\lambda},v_{\mu}$
associated to the eigenvalues $\lambda, \mu$, respectively. Let $Q =
[v_{\lambda}\ v_{\mu}]$ be the $2\times2$ matrix formed by the eigenvectors.
Then $Q$ is invertible. We can change the coordinates by using the
(invertible) linear transformation $Q^{-1}: \left[
\begin{array}
[c]{r}%
x\\
y
\end{array}
\right]  \mapsto\left[
\begin{array}
[c]{r}%
u\\
v
\end{array}
\right]  = Q^{-1}\left[
\begin{array}
[c]{r}%
x\\
y
\end{array}
\right]  $. Then the old linear problem $\left[
\begin{array}
[c]{r}%
x^{\prime}\\
y^{\prime}%
\end{array}
\right]  =Df(0)\left[
\begin{array}
[c]{r}%
x\\
y
\end{array}
\right]  $, $\left[
\begin{array}
[c]{r}%
x(0)\\
y(0)
\end{array}
\right]  = \left[
\begin{array}
[c]{r}%
x_{0}\\
y_{0}%
\end{array}
\right]  $, is reduced to a decoupled simpler problem which has the solution
$u(t)=u_{0}e^{\lambda t}$, $v(t)=v_{0}e^{\mu t}$.
The solution to the original linear problem can then be written as follows:
$\left[
\begin{array}
[c]{r}%
x\\
y
\end{array}
\right]  = Q\left[
\begin{array}
[c]{rr}%
e^{\lambda t} & 0\\
0 & e^{\mu t}%
\end{array}
\right]  Q^{-1}\left[
\begin{array}
[c]{r}%
x_{0}\\
y_{0}%
\end{array}
\right]  $. Since the matrix $Q$ is computable from $Df(0)$; in other words,
$Q$ is computable from (a $C^{1}$-name of) $f$, we may assume, without loss of
generality, at the outset that the standard basis of $\mathbb{R}^{2}$ is an
eigenbasis with $x$-axis being the stable separatrix and $y$-axis the unstable
separatrix. In other words, we assume that $Df(0)=\left[
\begin{array}
[c]{rr}%
\lambda & 0\\
0 & \mu
\end{array}
\right]  $ (see Fig.~\ref{fig:saddle_time}). \newline

In order to construct an algorithm that is able to output an approximation, we
need to estimate the time for a trajectory to pass by a saddle. We now turn to
address this problem.
Referring to Fig.~\ref{fig:saddle_time}, it is clear that any trajectory
starting in the interior of a quadrant will stay inside the same quadrant for
all $t$. Thus it suffices to consider the trajectories in the first quadrant.
Since the trajectory starting at a point in the first quadrant $z_{0}=(x_{0},
y_{0})$ is the graph of $z(t, z_{0}) = (e^{\lambda t}x_{0}, e^{\mu t}y_{0})$,
it is readily seen that $t$ can be solved as a function of $x$: $t=t(x)$. Let
$\gamma_{1}$ be a vertical line segment passing the point $(x_{1}, 0)$ and
$\gamma_{2}$ a horizontal line segment passing the point $(0, y_{2})$, where
$x_{1}, y_{2}>0$, and let $z(t,z_{\ast})$ be the trajectory of the linear
system starting at $z_{\ast}$, where $z_{\ast}=(x_{\ast}, y_{\ast})$,
$x_{\ast}>0$, is a point on the upper portion of $\gamma_{1}$.Then $z(t,
z_{\ast})$ will stay in the first quadrant for all $t$ and cross $\gamma_{2}$
at some time instant. Assume that $z(t, z_{\ast})$ crosses $\gamma_{2}$ for
the first time at $z_{\ast\ast}=(x_{\ast\ast}, y_{\ast\ast})$, $x_{\ast\ast
}>0$. Then this \textit{first time} can be computed by the formula below: Let
$T(z_{\ast})=\inf\{t>0: z(t, z_{\ast})\in\gamma_{2}\}$ be the time needed for
the trajectory $z(t, z_{\ast})$ to go from $z_{\ast}$ on $\gamma_{1}$ to a
point $z_{\ast\ast}=(x_{\ast\ast}, y_{\ast\ast})$ on $\gamma_{2}$. Then
\begin{equation}
\label{Eq:Time_saddle}T(z_{\ast})=\int_{x_{\ast\ast}}^{x_{\ast}}\frac
{dx}{\left\vert \lambda x\right\vert }%
\end{equation}
(see Theorem 8*.3.3 from \cite[p. 221]{HW95}). \newline

\begin{figure}
\begin{center}
\includegraphics{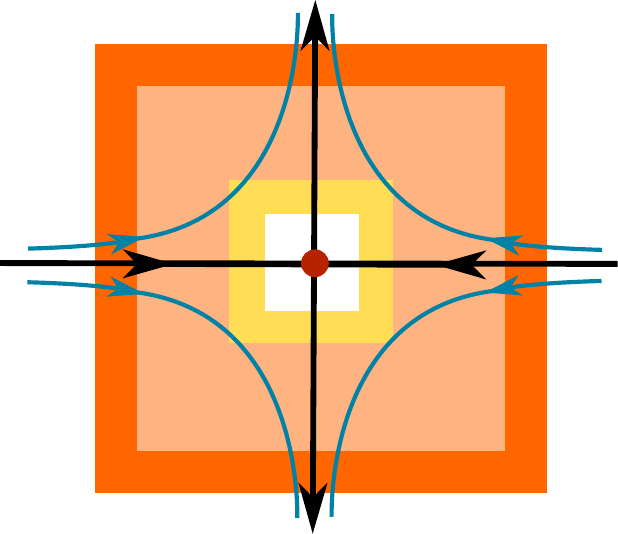}
\caption{Vicinity of a saddle point.}
\label{fig:saddle_vicinity}
\end{center}
\end{figure}

Now we turn to the construction of the black box at the saddle point $0$.
First we partition the ball $B(0,\varepsilon)$, which is the outer square in
Fig.~\ref{fig:saddle_vicinity} and it will remain fixed for any application of
$\mathtt{TimeEvolution}(\mathbb{D},\epsilon,T)$, into four regions:
$A=\{x\in\mathbb{R}^{2}:\left\Vert x\right\Vert \leq1/T\}$, $B=\{x\in
\mathbb{R}^{2}:1/T\leq\left\Vert x\right\Vert \leq2/T\}$, $C=\{x\in
\mathbb{R}^{2}:2/T\leq\left\Vert x\right\Vert \leq3\varepsilon/4\}$, and
$D=\{x\in\mathbb{R}^{2}:3\varepsilon/4\leq\left\Vert x\right\Vert
\leq\varepsilon\}$, where we assume that $2/T<3\varepsilon/4$. The region $A$
is depicted in white in Fig.~\ref{fig:saddle_vicinity}, $B$ in yellow, $C$ in
light orange, and $D$ in darker orange. Since the solution operator,
$e^{Df(0)}: (t, x) \mapsto e^{Df(0)t}x$, of the linear system is computable,
we are able to compute the image of (any subset or point of) $A,B,C$, or $D$
at any given time $t>0$. The neighborhood $H^{-1}_{i}(A\cup B\cup C)$ of $0$
serves as the black box.

Here is how the black box is programmed and incorporated with the simulation
$\mathtt{TimeEvolution}$. The simulation $\mathtt{TimeEvolution}$ performs
normally outside $U_{i}$ as well as in $U_{i}-H^{-1}_{i}(A\cup B\cup C)$ but
stops whenever a simulated trajectory enters the black box. The linear system
will then pick up a point on this simulated trajectory via $H_{i}$ as the
initial point, compute the trajectory until it reaches the region $D$, and
then returns a point in $D$ sufficiently close to the linear trajectory via
$H^{-1}_{i}$ to $\mathtt{TimeEvolution}$; the point supplied by the black box
will be $\rho/2$-close to some point on a trajectory of (\ref{ODE_Main}); upon
receiving the answer from the black box, $\mathtt{TimeEvolution}$ resumes its
normal working routine.

For a bit more details. Suppose that a trajectory computed by
$\mathtt{TimeEvolution}$ enters $H_{i}^{-1}(A\cup B\cup C)$. Then we transform
the system into its linearized version via the map $H_{i}$, obtaining a point
$z_{j}$ in $A\cup B\cup C$. Now we compute the solution of the linear system
with $z_{j}$ as the initial point, until it eventually reaches $D$. Assume
that $e^{Df(0)t_{1}}z_{j}\in D$ for some $t_{1}>0$. At this moment, we pick a
rational point $z_{j+1}\in D$ such that $|z_{j+1}-e^{Df(0)t_{1}}z_{j}|
\leq2^{-m(l)}$, where $m$ is a modulus of continuity for the homeomorphism
$H^{-1}_{i}$ and $2^{-l}\leq\rho/2$. Then $H^{-1}_{i}(z_{j+1})\in H^{-1}%
_{i}(D)\subseteq U_{i}$ and $|H^{-1}_{i}(z_{j+1}) - H^{-1}_{i}(e^{Df(0)t_{1}%
}z_{j})|\leq\rho/2$. In other words, $H^{-1}_{i}(z_{j+1})$ is an approximation
with error $\rho/2$ of a point on the trajectory $\phi_{t}(\bar{x})=H^{-1}%
_{i}(e^{Df(0)t}z_{j})$ of (\ref{ODE_Main}) starting at $\bar{x}=H^{-1}%
_{i}(z_{j})$. Now $\mathtt{TimeEvolution}$ resumes its normal work.

There is one problem with the argument above; that is, when the point $z_{j}$
lands on the stable separatrix -- the $x$ axis -- of the saddle $0$. In this
case, $e^{Df(0)t}z_{j}$ will stay on the stable separatrix for all $t>0$ and
move towards the origin; thus it will never reach $D$. In other words, the
black box won't be able to supply an answer to $\mathtt{TimeEvolution}$ in
this case. Even worse, if $z_{j}$ is on the stable separatrix or,
equivalently, the $x$-coordinate of $z_{j}$ equals zero, there is in general
no effective way of verifying this. The black box has to be re-programmed.

The remaining proof is devoted to re-programming the black box so that it is
able to return an answer, no matter where $z_{j}$ lands in $A\cup B\cup C$. We
begin by constructing, algorithmically, a finite set $J$ contained in $D$ to
be used as an answer to $\mathtt{TimeEvolution}$ from the black box when
$z_{j}$ is \textquotedblleft likely\textquotedblright\ to land on the stable
separatrix. Since the system is linear, it follows that, for any trajectory,
if it enters $A$ along a direction other than the $x$-axis, then it will leave
$A$, then $B$ and $C$, and enter $D$; in particular, once it leaves a region,
it won't re-enter that region ever again. Furthermore, for any trajectory
$z(t, z_{0})=(x_{0}e^{\lambda t}, y_{0}e^{\mu t})$, $y_{0}\neq0$, if it enters
$A$, $B$, $C$, or $D$, then the only way it can leave the region is through
the upper or lower horizontal border of that region, because the trajectory is
the graph of $x(t)=x_{0}e^{\lambda t}$ and $y(t)=y_{0}e^{\mu t}$, where
$\lambda< 0 < \mu$. Let $m=\min_{z\in B\cup C\cup D}\| Df(0)z\|$ and
$M=\max_{z\in B(0, \varepsilon)}\| Df(0)z\|$. Note that $m>0$ since $Df(0)$ is
invertible and $Df(0)z=0$ has a unique solution $z=0$. Both $m$ and $M$ are
computable from $f$, $\epsilon$, and $T$. Let region$_{H}$ denote the union of
the upper and the lower horizontal border of the region. For any point
$z_{\ast}=(x_{\ast}, y_{\ast})$ on $B_{H}$, the trajectory $z(t, z_{\ast}) =
e^{Df(0)t}z_{\ast}$ will not re-enter $B$ but leave $C$ and $D$, and
eventually reach $D_{H}$. The time it takes for $z(t, z_{\ast})$ to go from
$z_{\ast}$ to a point on $D_{H}$, say $z_{\ast\ast}=(x_{\ast\ast}, y_{\ast
\ast})$, is bounded by $3\varepsilon/m$ following (\ref{Eq:Time_saddle}):
\begin{equation}
0<\int_{x_{\ast\ast}}^{x_{\ast}}\frac{dx}{\left\vert \lambda x\right\vert
}\leq\int_{0}^{3\varepsilon}\frac{dx}{m}=\frac{3\varepsilon}{m}.
\label{Eq:TimeBound}%
\end{equation}
(note that any two points in $B(0,\varepsilon)$ are within Euclidean distance
$<3\varepsilon$). This time bound allows us to pick, effectively, a finite set
of points with rational coordinates on $B_{H}$, $W=\{ w_{1}, w_{2}, \ldots,
w_{q}\}$, which has the following property: for every $z_{\ast}\in B_{H}$,
there is some $w_{k}$ such that
\begin{equation}
\left\Vert w_{k}-z_{\ast}\right\Vert \leq\min\left(  \frac{\varepsilon
}{32e^{3\varepsilon M/m}},\frac{1}{T}\right)  \label{Eq:Approx}%
\end{equation}
(recall that $T$ is a rational number). Then it follows from Lemma
\ref{Lemma:Perturbation} that for any $t\in[0, 3\varepsilon/m]$,
\begin{equation}
\left\Vert z(t, w_{k})-z(t, z_{\ast})\right\Vert \leq\left\Vert w_{k}-z_{\ast
}\right\Vert e^{3\varepsilon M/m} \leq\varepsilon/32 \label{Eq:TimeBound1}%
\end{equation}
Without loss of generality we may assume that $\| H^{-1}_{i}(z_{1})-H^{-1}%
_{i}(z_{2})\|<\rho/2$ whenever $\| z_{1}-z_{2}\| \leq\varepsilon/32$.
Otherwise we simply increase the precision of (\ref{Eq:Approx}). (Recall that
the modulus of continuity of $H^{-1}_{i}$ is computable.) Next pick a rational
number $0<\tau< \varepsilon/(32M)$, and numerically compute $z(l\tau, w)$ for
every $w\in W$ and $l=0,1,\ldots,r=\left\lfloor 3\varepsilon/(\tau
m)\right\rfloor $ as long as $z((l+1)\tau, w)\in B(0,\varepsilon)$, where
$\left\lfloor x\right\rfloor =\max\{n\in\mathbb{Z}:n\leq x\}$. It is clear
that the following inequality holds true for any trajectory $z(t, w)$, $w\in
W$:
\begin{equation}
\label{Eq:TimeBound2}\left\Vert z((l+1)\tau, w)-z(l\tau, w)\right\Vert
\leq\varepsilon/32
\end{equation}
Inequality (\ref{Eq:TimeBound2}) further implies that for each $w\in W$, there
exists some $0\leq l_{w}\leq r$ such that $3\varepsilon/4 + 3\varepsilon
/32\leq\| z(l_{w}\tau, w)\| \leq3\varepsilon/4 + 5\varepsilon/32$. Finally it
follows from the inequalities (\ref{Eq:Approx}), (\ref{Eq:TimeBound1}), and
(\ref{Eq:TimeBound2}) that for every $z_{\ast}\in B_{H}$, there is a $w\in W$
such that $\| z(l_{w}\tau, z_{\ast}) - z(l_{w}\tau, w)\| \leq\varepsilon/32$;
in particular,
\[
3\varepsilon/4 + \varepsilon/16\leq\| z(l_{w}\tau, z_{\ast})\| \leq
3\varepsilon/4 + 3\varepsilon/16
\]
thus $z(l_{w}\tau, z_{\ast})\in D$. Let $J=\{ z(l_{w}\tau, w): \, w\in W\}$.
It follows from its construction that the finite set $J$ is computable from
$f$, $\varepsilon$, and $T$.

We are now ready to re-program the black box. Proceed as before: assume that
$z_{j}$ is a point in $A\cup B\cup C$ sent by $H_{i}$ and compute the
trajectory $e^{Df(0)t}z_{j}$. If the numerically computed trajectory of
$e^{Df(0)t}z_{j}$, with accuracy $\leq1/T$ and using the computable version of
Hartman Grobman's theorem, does not enter $B$, then for sure the (real)
trajectory $e^{Df(0)t}z_{j}$ will not enter $A$ and it will leave $C$ and then
enter $D$ because the flow is linear; in this case, the black box operates as
before and it will supply an answer to $\mathtt{TimeEvolution}$. On the other
hand, if the numerically computed trajectory of $e^{Df(0)t}z_{j}$ enters $B$
first, the (real) trajectory $e^{Df(0)t}z_{j}$ can enter the problematic
region $A$. To avoid the problem of region $A$, the black box will immediately
return, as the answer to $\mathtt{TimeEvolution}$, the $H^{-1}_{i}$-image of
the set $J$. It is not hard to see that the black box is now able to produce
an answer whenever there is a simulated trajectory entering it because if
$z_{j}$ lands on the stable separatrix of the saddle point $0$, then
$e^{Df(0)t}z_{j}$ will converge to the origin and thus it will enter $A$ at
some time instant and stay in $A$ thereafter.

A final note: the use of the black boxes does not affect the way how $NW(f)$
are approximated by the invariant sets of the simulation
$\mathtt{TimeEvolution}$. On the one hand, the black box does not throw out
any information on possible sets invariant under simulation
$\mathtt{TimeEvolution}$. The only trajectory that ``disappears" in the black
box is the trajectory starting on the stable separatrix of the saddle; but
this trajectory will move into the connected component of $Zero_{\epsilon}(f)$
containing the saddle when $T$ is large enough (or $\epsilon$ is small
enough). On the other hand, for every other trajectory that enters the black
box in a direction not along the stable separatrix, the black box will either
supply one point in $U_{i}$ or an over-approximation $H^{-1}_{i}(J)\subset
U_{i}$; in either case, every point supplied by the black box is within
$\rho/2$ error bound to a point on some trajectory of (\ref{ODE_Main}). Since
$\mathtt{TimeEvolution}$ uses $\rho/2$-grid for simulation, the invariant sets
of the simulation are not going to change when the black boxes are used.

Now we come to the construction of the desired algorithm -- the algorithm that
computes $NW(f)$ for the full case, where saddle points may exist. The
algorithm is constructed as follows: take as input $(k, f)$, where $f\in
SS_{2}$ and $k\in\mathbb{N}\backslash\{0\}$ ($k$ sets the accuracy $1/k$ for
the output approximation of $NW(f)$), do:

\begin{enumerate}
\item Let $\epsilon=1/(2k)$ and $T=1.$

\item Compute the set $Zero(f)$ with precision $1/k$ with the algorithm
described in \cite{GZ21a}, obtaining an over-approximation $Zero_{\epsilon}(f)$
of $Zero(f)$. If $Zero_{\epsilon}(f)\cap\partial K\neq\varnothing$, then take
$\epsilon:=\epsilon/2$, $T:=2T$, and repeat this step.

\item Compute the eigenvalues of each equilibrium in $Zero(f)$ and check
whether they have opposite signs. If yes, label the equilibrium point as a
saddle point.

\item If there are saddle points $z_{i}$, $i=1,\ldots,n$ compute
homeomorphisms $H_{i},H_{i}^{-1}$ and $\varepsilon>0$ as stated above, and use
the adapted version of $\mathtt{TimeEvolution}$.

\item Compute $\mathtt{TimeEvolution}(\mathbb{D},\epsilon,T)$ and obtain a
finite set $I$ for which (\ref{Eq:Time_evolution}) holds.

\item Take $P=\varnothing.$

\item For each $J\subseteq I$ check if $\cup_{j\in J}B(p_{j},\epsilon/2)$ is a
connected set and \texttt{HasInvariantSubset}$(D,\epsilon,T)=1$. If both
conditions hold, take $P=P\cup\{J\}$.

\item (The minimum principle) For each $J_{1},J_{2}\in P$, test if
$J_{1}\subsetneqq J_{2}$. If the test succeeds, then take $P=P\backslash
\{J_{2}\}.$

\item For each $J\in P$, test if $\cup_{j\in J}B(p_{j},\epsilon)$ is contained
in a connected component of $Zero_{\epsilon}(f)$. If the test succeeds, then
take $P=P\backslash\{J\}.$ If the test fails but $\cup_{j\in J}B(p_{j}%
,\epsilon)$ intersects a connected component of $Zero_{\epsilon}(f)$
containing a sink or a source, then take $\epsilon:=\epsilon/2$, $T:=2T$, and
go to 2. For each $J\in P$, test if $\cup_{j\in J}B(p_{j},\epsilon
)\cap\partial\mathbb{D}\neq\varnothing$. If the test succeeds, then take
$\epsilon:=\epsilon/2$, $T:=2T$.

\item For each $J\in P$ and each saddle point $z_{i}$, test if $\cup_{j\in
J}B(p_{j},\epsilon)$ intersects $H_{i}^{-1}(A\cup B)$ by computing the
distance $d(p_{j},H_{i}^{-1}(A\cup B))$, $j\in J$, and check whether
$d(p_{j},H_{i}^{-1}(A\cup B))<\epsilon$.
If this is true, then take $\epsilon:=\epsilon/2$, $T:=2T$, and go to step 5.

\item Do steps 5--10 to the ODE $x^{\prime}=-f(x)$ defined over $\mathbb{D},$
i.e.\ by taking the transformation $t\mapsto-t$ which reverses time in
(\ref{ODE_Main}), obtaining a set $\tilde{P}$ similar to the set $P$ of step 10.

\item For each set $J\in P\cup\tilde{P}$, use the algorithm described in
Section \ref{Sec:Limit_cycles} to check whether $\cup_{j\in J}B(p_{j}%
,\epsilon)$ has a doughnut shape with cross-section bounded by $1/k$. If this
check fails, then take $\epsilon:=\epsilon/2$, $T:=2T$, and go to step 5.

\item Check if $\mathtt{TimeEvolution}(\mathbb{D}-(Zero_{\epsilon}(f)\cup
~\cup_{J\in P\cup\tilde{P}}(\cup_{j\in J}B(p_{j},\epsilon)),\epsilon
/3,T)\subseteq Zero_{\epsilon}(f)\cup\cup_{J\in P\cup\tilde{P}}(\cup_{j\in
J}B(p_{j},\epsilon))$. If this test fails, then take $\epsilon:=\epsilon/2$,
$T:=2T$, and go to step 2.

\item Switch the dynamics to the ODE $x^{\prime}=-f(x)$ and test if
$\mathtt{TimeEvolution}(\mathbb{D}-(Zero_{\epsilon}(f)\cup~\cup_{J\in
P\cup\tilde{P}}(\cup_{j\in J}B(p_{j},\epsilon)),\epsilon/3,T)\subseteq
Zero_{\epsilon}(f)\cup\cup_{J\in P\cup\tilde{P}}(\cup_{j\in J}%
B(p_{j},\epsilon/2))$ for this ODE. If this test fails, then take
$\epsilon:=\epsilon/2$, $T:=2T$, and go to step 5.

\item Output $Zero_{\epsilon}(f)\bigcup\cup_{J\in P\cup\tilde{P}}%
\cup_{j\in J}B(p_{j},\epsilon)$.
\end{enumerate}

This algorithm, with the adaptations explained above, is essentially similar
to the one presented in the preceding section. A novel step that has not yet
been explained is step 10. The test is designed to eliminate possible
\textquotedblleft fake\textquotedblright\ cycles which might occur when the
unstable separatrix of a saddle point comes back very near to the stable
separatrix. The invariant sets output by step 9 will eventually pass step 10;
for the system is structurally stable and there is no saddle connection by
Theorem \ref{Th:GeneralPeixoto}. We note that since $H^{-1}_{i}: V_{i} \to
U_{i}$ is a homeomorphism and $A\cup B\subset V_{i}$ is a compact set, it
follows that $H^{-1}_{i}(A\cup B)$ is computable from $H^{-1}_{i}$ and $A\cup
B$, and thus so is the distance $d(p_{j}, H^{-1}_{i}(A\cup B))$.


\subsection{Notes about the computation of limit
cycles\label{Sec:Limit_cycles}}

In this section, we construct the algorithm needed for completing step 12 in
the preceding algorithm. Assume that $C_{i}$ is a connected component returned
from step 10 or step 11. Then $C_{i}\cap\partial\mathbb{D}=\emptyset$ and
$C_{i}$ contains at least one periodic orbit, say $\gamma_{i}$. According to
the Jordan Curve Theorem, $\gamma_{i}$ separates $\mathbb{R}^{2}$ into two
disjoint regions, the interior (a bounded region) and the exterior (an
unbounded region). The Poincar\'{e}-Bendixson Theorem further implies that
there is a square $s_{i}$ (of side length $\epsilon$) containing an
equilibrium point such that $s_{i}$ and $\partial\mathbb{D}$ lie in different
connected components of $\mathbb{D}- C_{i}$, $s_{i}$ is in the interior while
$\partial\mathbb{D}$ is in the exterior. We show that there exists a
sub-algorithm that, taking as input the sets returned by steps 10 or 11, halts
and outputs an approximation with an error bounded by $1/k$ of the set of all
periodic orbits contained in the output approximation.
The idea is to use a color-scheme algorithm to check whether $C_{i}$ is in the
shape of a \textquotedblleft doughnut\textquotedblright\thinspace\. Once $C_{i}$ is confirmed to have such a
shape, then the error in the approximation can be determined by measuring the
width of cross sections.

As the output of step 10 or step 11 of the main algorithm presented in Section
\ref{Sec:FullPicture}, $C_{i}$ has the form $C_{i}=\cup_{j\in J_{i}}%
B(p_{i,j},\epsilon)$; hence the closure of its complement, $\overline
{\mathbb{D}-C_{i}}$, can be written in the form of $\overline{\mathbb{D}%
-C_{i}}=\cup_{j\in\bar{J}_{i}}B(\bar{p}_{i,j},\epsilon_{i,j})$, where $\bar
{J}_{i}$ is a finite set of indices, $\epsilon_{i,j}$ are rational numbers
satisfying $0<\epsilon_{i,j}\leq\epsilon$, and $\bar{p}_{i,j}\in\mathbb{D}$
has rational coordinates. In the following, we call $B(\bar{p}_{i,j}%
,\epsilon_{i,j})$ a square (it can be viewed as a pixel).

The sub-algorithm embedded in step 12 of the main algorithm of Section
\ref{Sec:FullPicture} is defined as follows:

\begin{enumerate}
\item Choose a square $s$ contained in $\overline{\mathbb{D}-C_{i}}$ such that
$s\cap\partial\mathbb{D}\neq\varnothing$. Paint this square $s$ blue.

\item If $s^{\prime}$ is a square contained in $\overline{\mathbb{D}-C_{i}}$
adjacent to a blue square, then paint $s^{\prime}$ blue. Repeat this procedure
until there are no more squares which can be painted blue. Let $\widehat
{C}_{i,blue}$ be the union of all blue squares.

\item Pick an unpainted square $\tilde{s}\subseteq Zero_{\epsilon}%
(f)\cap\overline{\mathbb{D}-C_{i}}$, if there is any, and paint it red.

\item If $s^{\prime}$ is a square contained in $\overline{\mathbb{D}-C_{i}}$
adjacent to a red square, then paint $s^{\prime}$ red. Repeat this procedure
until there are no more squares which can be painted red. Let $\widehat
{C}_{i,red}$ be the union of all red squares.

\item Let $\widehat{C}_{i}=\overline{\mathbb{D}-\widehat{C}_{i,blue}%
-\widehat{C}_{i,red}}$. If $\widehat{C}_{i}\cap Zero_{\epsilon}(f)\neq
\varnothing$, return \texttt{False}\textbf{ }(the sub-algorithm has failed).

\item Compute an overapproximation $w_{i}$ with accuracy $\epsilon$ of the
quantity
\begin{equation}
\max\left(  \max_{x\in\widehat{C}_{i}}\min_{y\in\widehat{C}_{i,blue}
}\left\Vert x-y\right\Vert ,\max_{x\in\widehat{C}_{i}}\min_{y\in\widehat
{C}_{i,red}}\left\Vert x-y\right\Vert \right)  >0. \label{Eq:wi}%
\end{equation}
If $2w_{i}\leq\frac{1}{k}$, then return \texttt{True}\textbf{ }(the algorithm
succeeded), else return \texttt{False}.
\end{enumerate}

Note that all sets used in this sub-algorithm, $Zero_{\epsilon}(f)$,
$\widehat{C}_{i,blue}$, $\partial\mathbb{D}$, $\widehat{C}_{i}$, etc., are the
unions of finitely many polytopes with rational vertices. Thus, whether or not
their intersections are empty can be computed in finite time.


We note that steps 1--4 are rather straightforward and can be computed in a
finite amount of time. In addition, neither $\widehat{C}_{i,blue}$ nor
$\widehat{C}_{i,red}$ returned by steps 1--4 is empty as shown below.
Recall that $C_{i}\cap\partial\mathbb{D}=\emptyset$, $C_{i}\cap Zero_{\epsilon
}(f)=\emptyset$, $C_{i}$ contains at least one periodic orbit, say $\gamma
_{i}$, and $\partial\mathbb{D}$ is contained in the exterior region delimited
by $\gamma_{i}$. Subsequently, $\widehat{C}_{i,blue}$ is nonempty and it does
not intersect the interior of $\gamma_{i}$. If $\widehat{C}_{i,red}$ is empty,
then either $Zero_{\epsilon}(f)\subseteq\widehat{C}_{i,blue}$ or there is a
square in $Zero_{\epsilon}(f)$ that cannot be colored blue. If $Zero_{\epsilon
}(f)\subseteq\widehat{C}_{i,blue}$, then $Zero_{\epsilon}(f)$ is contained in
the exterior of $\gamma_{i}$, which is a contradiction to the
Poincar\'{e}-Bendixson Theorem. Hence there is at least one square in
$Zero_{\epsilon}(f)$ that cannot be colored blue. Since $C_{i}\cap
Zero_{\epsilon}(f)=\emptyset$, this square is available for being picked up by
step 3, which confirms that $\widehat{C}_{i,red}\neq\emptyset$.

We mention in passing that the reasoning that ensures steps 1--4 output
non-empty sets $\widehat{C}_{i,blue}$ and $\widehat{C}_{i,red}$ in finite time
is classical, and $\gamma_{i}$ is used in its classical capacity -- its
existence in $C_{i}$. In other words, the algorithm dictates computations; but
the correctness of the algorithm -- halting in finitely many steps with the
intended output -- is proved in a classical mathematical way. In the remaining
of this subsection, $\gamma_{i}$ is to be used repeatedly in this capacity to
show that steps 5 and 6 are guaranteed to halt with output True for
sufficiently small $\epsilon$.

Next we show that if the sub-algorithm returns a \texttt{True} answer on input
set $C_{i}$, then $C_{i}$ is ensured to be an over-approximation to the set of
all periodic orbits it contains with accuracy $\leq1/k$. Afterwards, we prove
that the sub-algorithm will return \texttt{True} for sufficiently large $T$
and small enough $\epsilon$.


Assume that the sub-algorithm returns a \texttt{True} answer on input set
$C_{i}$. Then every point in $\widehat{C}_{i}$ is guaranteed to be, at most,
at a distance of $\leq w_{i}$ from $\widehat{C}_{i,blue}$ as well from
$\widehat{C}_{i,red}$. In particular, this ensures that every point in
$\widehat{C}_{i}$ will be, at most, at a distance of $\leq2w_{i}$ from a point
of $\gamma$, where $\gamma$ is an arbitrary periodic orbit contained in
$\widehat{C}_{i}$. Indeed, let $x\in\widehat{C}_{i}$. Then there is some
$x_{red}\in\widehat{C}_{i,red}$ and some $x_{blue}\in\widehat{C}_{i,blue}$
such that
\begin{equation}
\left\Vert x-x_{red}\right\Vert \leq w_{i}\text{ and }\left\Vert
x-x_{blue}\right\Vert \leq w_{i}. \label{Eq:red_blue}%
\end{equation}
It follows from the triangular inequality that $\left\Vert x_{red}%
-x_{blue}\right\Vert \leq2w_{i}$. Since $x_{red}\in\widehat{C}_{i,red}$ and
$x_{blue}\in\widehat{C}_{i,blue}$, the line segment $\overline{x_{red}%
x_{blue}}$ will have to cross $\gamma$ at some point $y$, as to be shown
momentarily. Since the line segment has length bounded by $2w_{i}$, it follows
that $\left\Vert y-x_{red}\right\Vert \leq w_{i}$ or $\left\Vert
y-x_{blue}\right\Vert \leq w_{i}$, which in turn implies that $\left\Vert
x-y\right\Vert \leq2w_{i}$. Consequently, we arrive at the conclusion that the
Hausdorff distance between $\widehat{C}_{i}$ and the set of periodic orbits
contained inside $\widehat{C}_{i}$ is at most $2w_{i}$, provided that the
steps 5 \& 6 of the sub-algorithm both return True answers on input $C_{i}$.
Since $C_{i}\subset\widehat{C}_{i}$ by definition, the same conclusion holds
true for $C_{i}$. To show that any line segment $\Sigma$ that goes from a
point on $\widehat{C}_{i, blue}$ to a point on $\widehat{C}_{i, red}$ must
cross each and every periodic orbit contained in $\widehat{C}_{i}$, we argue
by way of a contradiction. Suppose $\gamma$ was a periodic orbit contained in
$\widehat{C}_{i}$ and $\gamma\cap\Sigma=\emptyset$. Then there is a square $s$
in $Zero_{\epsilon}(f)$ lying in the interior of $\gamma$ but not colored red
because $\widehat{C}_{i}\cap Zero_{\epsilon}(f)=\emptyset$ and the red region
is path connected by construction (see step 4 of the sub-algorithm). Since
$\widehat{C}_{i}$ does not contain any equilibrium point, it follows that each
square in $Zero_{\epsilon}(f)$ is colored either blue or red and thus $s$ has
color blue, which in turn implies that $s$ is in the same connected component
as of $\partial\mathbb{D}$. We arrive at a contradiction, for $s$ is in the
interior of $\gamma$ while $\partial\mathbb{D}$ in the exterior of $\gamma$.
We have now proved that if the sub-algorithm returns a \texttt{True} answer on
input set $C_{i}$, then $C_{i}$ is ensured to be an over-approximation to the
set of all periodic orbits contained in $C_{i}$ with accuracy $\leq1/k$.


It remains to show that the sub-algorithm will return \texttt{True} for
sufficiently large $T$ and small enough $\epsilon$. Let $C_{i}$ be an output
of step 10 or step 11 (of the main algorithm presented in section 8.2) that
fails either step 5 or step 6 of the sub-algorithm. Assume that $C_{i}$ is an
output of step 10, then $C_{i}$ contains at least one attractive periodic
orbit named $\gamma_{i}$ as in the previous paragraphs. Our strategy is to
show that, after finitely many updates on $\epsilon$, this periodic orbit
$\gamma_{i}$ will be over-approximated with an error bound $1/k$ in the sense
that there is a set, also called $C_{i}$ for simplicity, output by step 10
such that $\gamma_{i}$ is contained in $C_{i}$ and the sub-algorithm will
output True on input $C_{i}$. Since there are only finitely many periodic
orbits, the strategy ensures that, as an input to the sub-algorithm, every
output of step 10 or step 11 will return True, after finitely many updates on
$\epsilon$. (Recall that every output of step 10 or step 11 contains at least
one periodic orbit.) The strategy is executed as follows: first, we show that
$C_{i}$ will be inside a basin of attraction of $\gamma_{i}$ and $2w_{i}%
\leq\frac{1}{k}$ after finitely many updates on $\epsilon$ and $T$,
$\epsilon:=\epsilon/2$ and $T:=2T$; second, we prove that if $C_{i}$ is inside
a basin of attraction of $\gamma_{i}$, then $\widehat{C}_{i}\cap
Zero_{\epsilon}(f)=\emptyset$. Hence steps 5 \& 6 are guaranteed to halt with
outputs True.
We mention in passing that there is no need to find the exact number of
updates; it suffices to show $\mathit{classically}$ that there exists a
rational number such that the two mentioned conditions will be satisfied
whenever $\epsilon$ is updated to be less than this rational number.

Now for the details. Assume that $C_{i}$ is an output of step 10 for some
$0<\epsilon<1/(8k)$ and $T>0$. A similar argument applies to the case where
$C_{i}$ is an output of step 11. Then $C_{i}$ contains at least one attractive
periodic orbit, say $\gamma_{i}$.



We begin by showing that $C_{i}$ will be inside a basin of attraction of
$\gamma_{i}$ and $2w_{i} < 1/k$ after finitely many updates on $\epsilon$. Let
$\delta_{1}>0$ be a rational number such that $\mathcal{N}_{\delta_{1}}%
(\gamma_{i})$ is inside a basin of attraction of $\gamma_{i}$. Pick a point
$x_{0}$ on $\gamma_{i}$. Let $\Sigma\subset\mathbb{D}$ be a line segment
orthogonal to $\gamma_{i}$ at $x_{0}$; i.e., let $\Sigma= \{ x\in\mathbb{D}:
\, (x-x_{0})\cdot f(x_{0})=0\}$. (Note that $f(x_{0})\neq0$ since $\gamma_{i}$
cannot have any equilibrium point of (\ref{ODE_Main}) on it.) Then there is a
rational number $0<\delta_{2}<\delta_{1}$ and a unique function $\tau(x)$,
defined and continuously differentiable for $x\in\mathcal{N}_{\delta_{2}%
}(x_{0})$, such that $\mathcal{N}_{\delta_{2}}(x_{0}) \subset\mathcal{N}%
_{\delta_{1}}(\gamma_{i})$, $\tau(x_{0})=T_{\gamma_{i}}=$ the period of
$\gamma_{i}$, and $\phi_{\tau(x)}(x)\in\Sigma$. The restriction of $\tau$ on
$\mathcal{N}_{\delta_{2}}(x_{0})\cap\Sigma$ is called the first return map or
the Poincar\'{e} map. Note that $\mathcal{N}_{\delta_{2}}(x_{0})\cap\Sigma$
intersects both the interior and exterior delimited by $\gamma_{i}$. Pick two
points, $y_{0}$ and $z_{0}$, on $\mathcal{N}_{\delta_{2}}(x_{0})\cap\Sigma$
such that $y_{0}$ is in the interior while $z_{0}$ is in the exterior.
Let $y_{n+1}=\phi_{\tau(y_{n})}(y_{n})$ and $z_{n+1}=\phi_{\tau(z_{n})}%
(z_{n})$. Then the trajectory $\phi_{t}(y_{0})$ ($\phi_{t}(z_{0})$,
respectively) stays in the interior (exterior, respectively) for all $t>0$ and
$y_{n} \to x_{0}$ ($z_{n}\to x_{0}$, respectively) monotonically as
$n\to\infty$ (see, e.g., \cite[Lemma 8*.5.10,on p. 242]{HW95}).

\begin{figure}[ptb]
    \begin{center}
        \includegraphics[width=5cm]{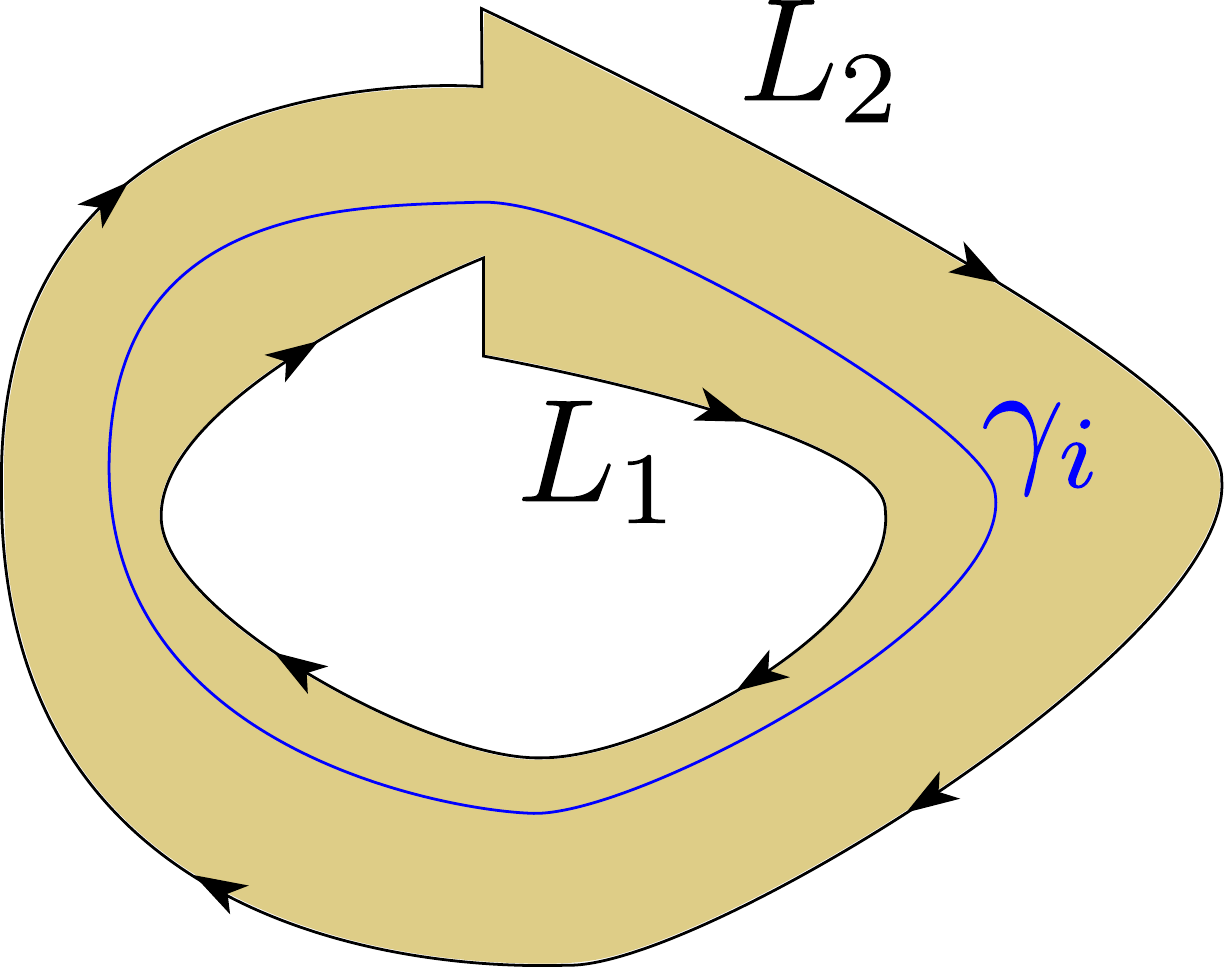}
    \end{center}
    \caption{Curves enclosing a periodic orbit $\gamma_{i}$.}%
    \label{Fig:curves}%
\end{figure}

Let $t_{0}>T_{\gamma_{i}}$ be a positive number such that $Ke^{-\alpha
(t_{0}-T_{\gamma_{i}})}\leq1/(8k)$, where $T_{\gamma_{i}}$ is the period of
$\gamma_{i}$. Pick an $\tilde{n}\in\mathbb{N}$ such that $\Vert y_{\tilde{n}%
}-x_{0}\Vert\leq1/(8k)$, $\Vert z_{\tilde{n}}-x_{0}\Vert\leq1/(8k)$, and
\[
\min\left\{  \sum_{n=0}^{\tilde{n}-1}\tau(y_{n}),\sum_{n=0}^{\tilde{n}-1}%
\tau(z_{n})\right\}  >t_{0}%
\]
where
\begin{equation}
K=\max_{1\leq j\leq q}K_{j}\quad\mbox{and}\quad\alpha=\min_{1\leq j\leq
q}\alpha_{j}/T_{j} \label{Eq:K_alpha}%
\end{equation}
$K_{j}$, $\alpha_{j}$ and $T_{j}$ are the corresponding values of $K,\alpha,T$
provided by Proposition \ref{Prop:ConvPeriodic} for the periodic orbit
$\gamma_{j}$. Then it follows from Proposition \ref{Prop:ConvPeriodic} that
\[
d(\phi(t_{0}-T_{\gamma_{i}}+t,y_{0}),\gamma_{i})<1/(8k),d(\phi(t_{0}%
-T_{\gamma_{i}}+t,z_{0}),\gamma_{i})<1/(8k)
\]
for every $t\geq0$. Let $L_{1}$ be the trajectory of (\ref{ODE_Main}) from
$y_{\tilde{n}}$ to $y_{\tilde{n}+1}$, $L_{2}$ the trajectory of
(\ref{ODE_Main}) from $z_{\tilde{n}}$ to $z_{\tilde{n}+1}$, $l_{1}$ the line
segment from $y_{\tilde{n}}$ to $y_{\tilde{n}+1}$, and $l_{2}$ the line
segment from $z_{\tilde{n}}$ to $z_{\tilde{n}+1}$ (see Fig.~\ref{Fig:curves}). It is clear that $L_{1}%
\cup\,l_{1}$ is a simple closed curve in the interior of $\gamma_{i}$ and
$L_{2}\cup l_{2}$ a simple closed curve in the exterior of $\gamma_{i}$. Since
$y_{0}, z_{0}\in\mathcal{N}_{\delta_{2}}(\gamma_{2})\cap\Sigma\subset
\mathcal{N}_{\delta_{1}}(\gamma_{i})$ and $\mathcal{N}_{\delta_{1}}(\gamma
_{i})$ is a basin of attraction of $\gamma_{i}$, it follows that the closed
curves $L_{1}\cup\, l_{1}$ and $L_{2}\cup\, l_{2}$ as well the region bounded
by them are contained inside $\mathcal{N}_{\delta_{1}}(\gamma_{i})$.

We show next that the Hausdorff distance between $L_{1}\cup l_{1}$ and
$\gamma_{i}$ as well as between $L_{2}\cup l_{2}$ and $\gamma_{i}$ is bounded
by $1/(8k)$. It suffices to show that $d(L_{1}\cup l_{1},\gamma_{i}%
(t-t_{y_{0}}))\leq1/(8k)$ for every $T^{\prime}\leq t\leq T^{\prime}%
+T_{\gamma_{i}}$, where $t_{y_{0}}$ is the asymptotic phase corresponding to
$y_{0}$, for some $T^{\prime}$. Two cases are to be considered:

Case (1): $\tau(y_{\tilde{n}})\geq T_{\gamma_{i}}$. In this case, by taking $T^{\prime}=\sum_{n=0}^{\tilde{n}-1}\tau(y_{n})$, for any
\[
\sum_{n=0}^{\tilde{n}-1}\tau(y_{n}) \leq t \leq\sum_{n=0}^{\tilde{n}-1}%
\tau(y_{n}) + T_{\gamma_{i}} \leq\sum_{n=0}^{\tilde{n}}\tau(y_{n})
\]
we have either $\phi(t, y_{0})=y_{\tilde{n}}$, or $\phi(t, y_{0})=y_{\tilde
{n}+1}$, or
\[
\phi(t - \sum_{n=0}^{\tilde{n}-1}\tau(y_{n}), y_{\tilde{n}}) = \phi(t, y_{0})
\in L_{1}
\]
and
\[
\| \phi(t, y_{0}) - \gamma_{i}(t - t_{y_{0}})\| \leq1/(8k)
\]

Case (2): $\tau(y_{\tilde{n}})<T_{\gamma_{i}}$. Let $T_{d}$ be the positive
number
\[
0<T_{d}=T_{\gamma_{i}}-\tau(y_{\tilde{n}})<T_{\gamma_{i}}%
\]
and take $T^{\prime}=\sum_{n=0}^{\tilde{n}-1}\tau(y_{n})-T_{d}$. Then for any
\[
\sum_{n=0}^{\tilde{n}-1}\tau(y_{n})-T_{d}\leq t\leq\sum_{n=0}^{\tilde{n}%
-1}\tau(y_{n})+T_{\gamma_{i}}-T_{d}=\sum_{n=0}^{\tilde{n}}\tau(y_{n})
\]
if
\[
\sum_{n=0}^{\tilde{n}-1}\tau(y_{n})\leq t\leq\sum_{n=0}^{\tilde{n}}\tau
(y_{n})
\]
then the proof is similar as that for case (1). On the other hand, if
\[
\sum_{n=0}^{\tilde{n}-1}\tau(y_{n})-T_{d}\leq t\leq\sum_{n=0}^{\tilde{n}%
-1}\tau(y_{n})
\]
then there exists some $1\leq m\leq\tilde{n}-1$ with
\[
\sum_{n=0}^{m-1}\tau(y_{n})\leq t\leq\sum_{n=0}^{m}\tau(y_{n})
\]
which implies that $\phi(t,y_{0})$ is on the trajectory somewhere between
$y_{m}$ and $y_{m+1}$. Since $t\geq t_{0}-T_{\gamma_{i}}$, we also have
\[
\Vert\phi(t,y_{0})-\gamma_{i}(t-t_{y_{0}})\Vert\leq1/(8k)
\]
Note that $\phi(t,y_{0})$ is not on $L_{1}\cup l_{1}$. However, since
$L_{1}\cup l_{1}$ is a closed curve surrounding $\gamma_{i}$ and
$y_{n}\rightarrow x_{0}$ monotonically, it follows that $L_{1}\cup l_{1}$ lies
between the trajectory from $y_{m}$ to $y_{m+1}$ and $\gamma_{i}$; thus the
line goes from the point $\phi(t,y_{0})$ to the point $\gamma_{i}(t-t_{y_{0}%
})$ must cross a point on $L_{1}$. Consequently,
\begin{equation}
d(L_{1}\cup l_{1},\gamma_{i}(t-t_{y_{0}}))\leq1/(8k) \label{Eq:closed_curve1}%
\end{equation}


Define
\[
\delta_{k}=\min\left\{  \min_{x\in L_{1}\cup l_{1}}d(x,\gamma_{i}),\min_{x\in
L_{2}\cup l_{2}}d(x,\gamma_{i}), \delta_{1} \right\}
\]
Then $\delta_{k}>0$ since the compact sets $\gamma_{i}$, $L_{1}\cup l_{1}$,
and $L_{2}\cup l_{2}$ are mutually disjoint. It is clear that $\mathcal{N}%
_{\delta_{k}}(\gamma_{i})$ is contained in $\mathcal{N}_{\delta_{1}}%
(\gamma_{i})$ as well in the region bounded by $L_{1}\cup l_{1}$ and
$L_{2}\cup l_{2}$.

We now show that the main algorithm will output an updated $C_{i}$ that is
contained in $\mathcal{N}_{\delta_{k}}(\gamma_{i})$ after finitely many
updates on $\epsilon$, with the property that $2w_{i}<1/k$ and $C_{i}$ is in a
basin of attraction of $\gamma_{i}$. This is indeed the case as shown as
follows. Recall that the main algorithm starts with $\epsilon=1/k$ and $T=1$.
With each iteration of the algorithm, the updates $\epsilon:=\epsilon/2$ and
$T:=2T$ are performed. Thus, after $n$ iterations, $T=2^{n}$ and
$\epsilon=1/(2^{n}k)$. Let $0<\zeta<\delta_{k}$. Then it is readily seen that
the condition $2\epsilon<\zeta$ will be met after $n_{0}$ iterations with
$n_{0}>1+|\log_{2}(k\zeta)|$. Similarly, for any $x\in\mathcal{N}_{\delta_{1}%
}(\gamma_{i})$ but $x\notin\gamma_{i}$, the condition $0<d(\phi_{t}%
(x),\gamma_{i})<\zeta/4$ is reached after $n_{1}$ iterations as long as
$Ke^{-\alpha2^{n_{1}}}<\zeta/4$ or, equivalently, $n_{1}>\left\vert \log
_{2}\left(  \frac{1}{\alpha}\ln\left(  \frac{4K}{\zeta}\right)  \right)
\right\vert $, where $K$ and $\alpha$ are defined in (\ref{Eq:K_alpha}). Note
that after $n_{1}$ iterations, $T\geq2^{n_{1}}$. Set
\[
N=\max\left\{  1+|\log_{2}(k\zeta)|,\left\vert \log_{2}\left(  \frac{1}%
{\alpha}\ln\left(  \frac{4K}{\zeta}\right)  \right)  \right\vert \right\}
\]
and assume that the main algorithm is on the $n$th iteration for some $n\geq
N$, and $C_{i}=\cup_{j\in J_{i}}B(p_{j},\epsilon/2)$ is an output of step 10
that contains $\gamma_{i}$. Let $\widetilde{C}_{i}=\cup\{B(p_{j}%
,\epsilon/2):\,j\in J_{i},\,d(p_{j},\gamma_{i})\leq3\zeta/4\}$. Then
$\widetilde{C}_{i}$ covers $\gamma_{i}$. Moreover,the following are true: (i)
if $x\in\widetilde{C}_{i}$, then $d(x,\gamma_{i})<\zeta$ because
$\epsilon/2<\zeta/4$, which implies that $\widetilde{C}_{i}\subset
\mathcal{N}_{\zeta}(\gamma_{i})\subset\mathcal{N}_{\delta_{1}}(\gamma_{i})$;
(ii) if $y\in\phi_{T}(\widetilde{C}_{i})$, then $d(y,\gamma_{i})<\zeta/4$
because $T\geq2^{n_{1}}$; and (iii) $\mathtt{TimeEvolution}(\widetilde{C}%
_{i},\epsilon,T)\subseteq\widetilde{C}_{i}$ (invariance of $\widetilde{C}_{i}$
by the operator $\mathtt{TimeEvolution}$). Since, by construction, $C_{i}$ is
the minimal invariant set by $\mathtt{TimeEvolution}$ that includes
$\gamma_{i}$, this implies that $C_{i}\subseteq\widetilde{C}_{i}$. Thus, for
every $x\in C_{i}$, $d(x,\gamma_{i})<\zeta$. Therefore, $C_{i}\subseteq
\mathcal{N}_{\zeta}(\gamma_{i})\subset\mathcal{N}_{\delta_{k}}(\gamma
_{i})\subset\mathcal{N}_{\delta_{1}}(\gamma_{i})$. This shows that $C_{i}$ is
inside the region enclosed by the closed simple curves $L_{1}\cup l_{1}$ and
$L_{2}\cup l_{2}$ (we shall call it the $L$-region for simplicity), which in
turn is contained in $\mathcal{N}_{\delta_{1}}(\gamma_{i})$, a basin of
attraction of $\gamma_{i}$.

This suggests that if we can show that $\widehat{C}_{i}$ is contained in the
$L$-region and every square in $Zero_{\epsilon}(f)$ is disjoint from the
$L$-region, then $\widehat{C}_{1}\cap Zero_{\epsilon}(f)=\emptyset$. Hence
step 5 returns True.

To show that $\widehat{C}_{i}$ is contained in the $L$-region, we begin with
the assumption that $C_{i}\subseteq\mathcal{N}_{\delta_{k}/3}(\gamma_{i})$ and
$\epsilon\leq\delta_{k}/3$. Clearly this can be done as proved above. Recall
that $\mathcal{N}_{\delta_{k}}(\gamma_{i})$ is inside the $L$-region and
$C_{i}$ is in the exterior of $L_{1}\cup\, l_{1}$. Then for every $x\in
L_{1}\cup\, l_{1}$ or in its interior,
\begin{equation}
\label{L-estimate}d(x, C_{i})\geq d(x, \mathcal{N}_{\delta_{k}/3}(\gamma
_{i})\geq2\delta_{k}/3 \geq2\epsilon> \epsilon.
\end{equation}
We already know that there is at least one red square $s$ (of side length
$\epsilon$ containing an equilibrium point $x_{e}$) in the interior of
$\gamma_{i}$ and all equilibrium points in the interior of $\gamma_{i}$ must
be in the interior of $L_{1}\cup\, l_{1}$ due to the fact that the $L$-region
is in a basin of attraction of $\gamma_{i}$. Thus this red square covers a
portion of the interior of $L_{1}\cup\, l_{1}$ and is disjoint from $C_{i}$ by
(\ref{L-estimate}). We shall show that $L_{1} \cup\, l_{1}$ together with its
interior has color red, which implies that $\widehat{C}_{i, red}$ covers
$L_{1} \cup\, l_{1}$ and its interior. It can be shown similarly that
$\widehat{C}_{i, blue}$ covers $L_{2}\cup\, l_{2}$ and its exterior. This
ensures that $\widehat{C}_{i} = \mathbb{D} - \widehat{C}_{i, red} -
\widehat{C}_{i, blue}$ is indeed a subset of the $L$-region. We now look at
the interior of $L_{1}\cup\, l_{1}$. By the Jordan curve theorem, the interior
is a path-connected region and, by (\ref{L-estimate}), every square with
center on $L_{1}\cup\, l_{1}$ or in its interior having side length $\epsilon$
is disjoint from $C_{i}$. Then step 4 of the subalgorithm colors every such
square red starting with the red square containing $x_{e}$. Since each square
in $Zero_{\epsilon}(f)$ has side-length $\epsilon$ and contains one
equilibrium point either in the interior of $L_{1}\cup\, l_{1}$ or in the
exterior of $L_{1}\cup\, l_{2}$, it follows that $Zero_{\epsilon}%
(f)\subset\widehat{C}_{i, red}\cup\, \widehat{C}_{i, blue}$. Hence
$\widehat{C}_{i}\cap Zero_{\epsilon}(f)=\emptyset$.

Finally, it follows from (\ref{Eq:closed_curve1}) that $d(\widehat{C}%
_{i},\gamma_{i})\leq1/(8k)$. Hence the quantity (\ref{Eq:wi}) is bounded by
$1/(4k)$. We recall that $w_{i}$ is an over-approximation of the quantity
(\ref{Eq:wi}) with an error bound $\epsilon\leq1/(8k)$. This implies that
$w_{i}<1/(2k)$ or, equivalently, $2w_{i}<1/k$. Thus step 6 also returns True.
The proof is complete.

We conclude this section with the computation of a cross-section for each
$\widehat{C}_{i}$, where $\widehat{C}_{i}$ is the output of a successful run
of the subalgorithm, $i$ is in some finite index set $I$, and $\bigcup_{i\in
I} C_{i}$ is an over-approximation of $Per(f)$ with error bound $1/k$. Recall
that a cross-section of $\widehat{C}_{i}$ is a line segment that lies in
$\widehat{C}_{i}$, is transversal to all trajectories across it, and
intersects with all periodic orbits contained in $\widehat{C}_{i}$. These
cross-sections will be used in the next section for computing the exact number
of periodic orbits contained in $C_{i}$. The same technique used to write the
subalgorithm can be extended to compute the desired cross-sections.

Given $C_{i}$, let $\widehat{C}_{i}$, $\widehat{C}_{i,red}$, $\widehat
{C}_{i,blue}$ be the corresponding sets obtained by a successful run of the
preceding subalgorithm, $i\in I$. Then to compute a cross-section for
$\widehat{C}_{i}$, proceed as follows:

\begin{enumerate}
\item Compute a rational point $p_{i}$ on $\widehat{C}_{i}\cap\widehat
{C}_{i,red}$ (it suffices to look at the vertices of the finitely many squares
which form $\widehat{C}_{i}$ and pick one such vertex which is also painted
red. Note that all the vertices of the squares defining $\widehat{C}_{i}$ and
$\widehat{C}_{i,red}$ have rational coordinates, and therefore such a $p_{i}$
can be computed in finite time).

\item Consider the vector $n(p_{i})=(-f_{2}(p_{i}),f_{1}(p_{i}))$ which is
orthogonal fo $f(p_{i})\neq0$. Compute a rational approximation $v_{i}$ to
$n(p_{i})$ with accuracy bounded by $\epsilon$.

\item Test whether $\measuredangle(v_{i},n(p_{i}))\leq\pi/10$ or
$\measuredangle(v_{i},n(p_{i}))\geq\pi/11$, where $\measuredangle(x, y)$
denotes the positive angle between vectors $x$ and $y$. If $\measuredangle
(v_{i},n(p_{i}))\geq\pi/11$, then update $\epsilon:=\epsilon/2$ and go to step 2.

\item Let $s_{1},\ldots,s_{l}$ be the squares (with rational vertices)
computed by the preceding subalgorithm such that $\widehat{C}_{i,blue}%
=\cup_{j=1,\ldots,l}s_{j}$; let $l_{v_{i}}$ be the ray starting at $p_{i}$ and
parallel to $v_{i}$ such that $p_{i}$ is the only red point on $l_{v_{i}}$.
Starting from $j=1$, for each $1\leq j\leq l$, decide whether $s_{j}\cap
l_{v_{i}}\mathbb{=\varnothing}$. If the condition holds true, take $r_{j}$
being the point $(3, 0)$ in $\mathbb{R}^{2}$ and move to $s_{j+1}$. If
$s_{j}\cap l_{v_{i}} \neq\emptyset$, compute the point $r_{j}\in\mathbb{D}$
(which has rational coordinates) satisfying $\left\Vert r_{j}-p_{i}\right\Vert
=\min_{y\in s_{j}\cap l_{v_{i}} }\left\Vert y-p_{i}\right\Vert $, and then
move to $s_{j+1}$. (We note that, in the latter case, $r_{j}$ is a blue point.)

\item Take $q_{i}$ be some $r_{j}$ satisfying $\left\Vert q_{i} -
p_{i}\right\Vert =\min_{j=1,\ldots,l}\left\Vert r_{j}-p_{i}\right\Vert $.

\item Compute
\[
\theta=\max_{z\in\overline{p_{i}q_{i}}}\measuredangle(f(p_{i}),f(z))=\left|
\arccos\left(  \frac{f(p_{i})\cdot f(z)}{\left\Vert f(p_{i})\right\Vert
\left\Vert f(z)\right\Vert }\right)  \right|
\]
and test whether $\theta\leq\pi/10$ or $\theta\geq\pi/11$. If $\theta\geq
\pi/11$, then update $\epsilon:=\epsilon/2$ and $T:=2T$ in the main algorithm
of Section \ref{Sec:FullPicture}, obtaining new sets $C_{i}$, $\widehat{C}%
_{i},$ $\widehat{C}_{i,red},$ $\widehat{C}_{i,blue}$ and repeat the current
algorithm starting from step 1.

\item Output the line segment $\overline{p_{i}q_{i}}$ as a cross-section for
$\widehat{C}_{i}$.
\end{enumerate}

We need to show that the algorithm halts after finitely many updates on
$\epsilon$; when it halts, it returns a cross-section of $\widehat{C}_{i}$.
The first five steps are rather straightforward, for $f(p_{i})$ is computable
and every square in $\widehat{C}_{i,red}$ or $\widehat{C}_{i,blue}$ has
rational corners and rational side-length. Since the ray $l_{v_{i}}$ starting
from the red point $p_{i}$ will move out of $\mathbb{D}$, it must cross the
blue -colored $\partial\mathbb{D}$. Hence, there is at least one square, say
$s_{j}$, in $\widehat{C}_{i,blue}$ such that $s_{j}\cap l_{v_{i}}\neq
\emptyset$. This ensures that $q_{i}$ has color blue and the length of
$\overline{p_{i}q_{i}}>0$ (recall that $\widehat{C}_{i,red}\cap\widehat
{C}_{i,blue}=\emptyset$). We now turn to step 6. We begin with the definition
of the function
\[
h:\bigcup_{i\in I}\widehat{C}_{i}\times\bigcup_{i\in I}\widehat{C}%
_{i}\rightarrow\mathbb{R},\quad h(x,y)=\measuredangle(f(x),f(y))=\arccos
\left(  \frac{f(x)\cdot f(y)}{\left\Vert f(x)\right\Vert \left\Vert
f(y)\right\Vert }\right)
\]
For simplicity, we call $\bigcup_{i\in I}\widehat{C}_{i}$ a
hat-over-approximation of $Per(f)$. Since $\bigcup_{i\in I}\widehat{C}_{i}$ is
a compact subset of $\mathbb{D}$ and $\widehat{C}_{i}\cap Zero_{\epsilon
}(f)=\emptyset$ for all $i\in I$, the function $h$ is well-defined and
uniformly continuous on $\bigcup_{i\in I}\widehat{C}_{i}\times\bigcup_{i\in
I}\widehat{C}_{i}$. Thus, there is a rational number $\delta>0$ such that
\begin{equation}
|h(x_{1},y_{1})-h(x_{2},y_{2})|<\pi/10\quad\mbox{whenever}\quad\Vert
(x_{1},y_{1})-(x_{2},y_{2})\Vert\leq\delta\label{function_h}%
\end{equation}
Hence, if the length of $\overline{p_{i}q_{i}}$ is no larger than $\delta$,
then $\theta\leq\pi/10$. On the other hand, suppose the length of
$\overline{p_{i}q_{i}}$ is larger than $\delta$. In this case, we recall
briefly some facts which were worked out in detail in the proof of the
preceding subalgorithm. First, when $\epsilon$ and $T$ are updated, the
updated hat-over-approximation is a subset of the hat-over-approximation
before the update. This indicates that (\ref{function_h}) holds true for any
updated hat-over-approximation. Second, given any positive integer $l$, the
main algorithm of section 8.2 and the preceding subalgorithm can output a
hat-over-approximation such that $w_{i}+\epsilon<1/l$ for all $i$ by updating
$\epsilon$ and $T$ finitely many times, where $w_{i}$ is an approximation with
accuracy $\epsilon$ to the quantity defined in (\ref{Eq:wi}). Hence, by
picking $l$ such that $(1/l)<\delta$ and by updating $\epsilon$ and $T$, the
condition $w_{i}+\epsilon<\delta$ is ensured to be met for all $i$. Now since
the length of $\overline{p_{i}q_{i}}$ cannot be larger than $w_{i}+\epsilon$
by definition of $q_{i}$, it follows that the length of $\overline{p_{i}q_{i}%
}$ will be bounded by $\delta$ for all $i$ after updating $\epsilon$ and $T$
finitely many times. Hence, step 6 will halt and output $\theta\leq\pi/10$. As
the last step, we show that $\overline{p_{i}q_{i}}$ is indeed a cross-section
of $\widehat{C}_{i}$. It is clear that $\overline{p_{i}q_{i}}$ lies in
$\widehat{C}_{i}$. We recall that it has been shown in the proof of the
preceding subalgorithm that any line segment that goes from a point on
$\widehat{C}_{i,red}$ to a point on $\widehat{C}_{i,blue}$ must cross each and
every periodic orbit contained in $\widehat{C}_{i}$. It remains to show that
the trajectories move through $\overline{p_{i}q_{i}}$ transversely at every
point on $\overline{p_{i}q_{i}}$. For any $z\in\overline{p_{i}q_{i}}$, since
$\measuredangle(f(p_{i}),f(z))\leq\theta\leq\pi/10$ and $\measuredangle
(v_{i},n(p_{i}))\leq\pi/10$, it follows that $\measuredangle(f(z),v_{i}%
)\geq\measuredangle(f(p_{i}),n(p_{i}))-\measuredangle(n(p_{i}),v_{i}%
)-\measuredangle(f(p_{i}),f(z))\geq(\pi/2)-(\pi/10)-(\pi/10)=3\pi/10$. Hence,
the trajectories cross $\overline{p_{i}q_{i}}$ transversely at every point on
$\overline{p_{i}q_{i}}$.


\section{Computing Poincar\'{e} maps\label{Sec:Poincare}}

In this section, we make use of Poincar\'{e} maps (or first return maps) to
construct an algorithm for computing the number of periodic orbits contained
in each $C_{i}$. A cross-section is needed in order to define a first return
map. Since we have an algorithm for computing a corss-section of $\widehat
{C}_{i}$, it is natural to work with $\widehat{C}_{i}$ instead of $C_{i}$
provided that $\widehat{C}_{i}$ is invariant for all $t\geq T$ for some $T>0$
and $\widehat{C}_{i} \cap\widehat{C}_{j} = \emptyset$ whenever $i\neq j$. The
first condition guarantees that the first return map can be defined in a
neighborhood of a cross-section when $T$ is sufficiently large, and the second
condition ensures that $\widehat{C}_{i}$ contains the exact number of periodic
orbits as $C_{i}$ does because $C_{i}\subseteq\widehat{C}_{i}$ and
$Per(f)\subset\cup_{i}C_{i}$.

We give a sketch that both conditions shall be met. Suppose $\widehat{C}_{i}$,
$i\in I$, are True outcomes of the subalgorithm with input parameters
$\epsilon$ and $T$. Now compute \texttt{HasInvariantSubset}$(\widehat{C}_{i},
\epsilon, T)$ and test whether $\widehat{C}_{i}\cap\widehat{C}_{j}=\emptyset$.
If \texttt{HasInvariantSubset}$(\widehat{C}_{i}, \epsilon, T) = 1$ for all $i$
and $\widehat{C}_{i}\cap\widehat{C}_{j}=\emptyset$ whenever $i\neq j$, then
return True. Otherwise, update $\epsilon:= \epsilon/2$ and $T:= 2T$ and rerun
the main- and the sub-algorithm. It can be proved that the computation and the
test will return True after finitely many updates on $\epsilon$ and $T$ by an
argument similar to that used to confirm that steps 5 and 6 in the
subalgorithm will return True after finitely many updates on $\epsilon$ and
$T$.

In the remainder of this section, we assume that $\widehat{C}_{i}$ are the
True returns and $Per(f)\subset\cup_{i} \widehat{C}_{i}$. For simplicity, we
further assume that $\widehat{C}_{i}$ is invariant for all $t > 0$ by a
transformation of time. For each $\widehat{C}_{i}$, let $\overline{p_{i}q_{i}%
}$ be the cross-section of $\widehat{C}_{i}$ computed at the end of the
previous section. Then a Poincar\'{e} map $P_{i}$ can be defined on
$\overline{p_{i}q_{i}}$: it assigns to every $p$ on $\overline{p_{i}q_{i}}$
the point on $\overline{p_{i}q_{i}}$ that is first reached by following the
trajectory $\phi_{t}(p)$ for $t>0$. In particular, a point $q$ on
$\overline{p_{i}q_{i}}$ is on a periodic orbit iff $q$ is a fixed point of the
Poincar\'{e} map, $P_{i}(q)=q$. Hence, in order to compute the number of
periodic orbits contained in $\widehat{C}_{i}$, we just need to compute the
number of fixed points of $P_{i}$ or, equivalently, the number of zeros of the
function $Q_{i}:\overline{p_{i}q_{i}}\rightarrow\mathbb{R}^{2}$ defined by
$Q_{i}(x)=P_{i}(x)-x$. The number of zeros of $Q_{i}$ can be computed using
the algorithm from Section \ref{Sec:Zeros} as long as $Q_{i}$ has the
following properties:


\begin{enumerate}
\item If $q$ is a zero of $Q_{i}$, then the jacobian of $Q_{i}$ at point $q$,
$DQ_{i}(q)$, is invertible;

\item $Q_{i}$ and $DQ_{i}$ are computable from $f$ of (\ref{ODE_Main}).
\end{enumerate}

It is well known that if all periodic orbits in $C_{i}$ are hyperbolic, then
condition 1 holds.

Concerning condition 2, it suffices to show that the Poincar\'{e} map $P_{i}$
and its derivative $DP_{i}$ are computable from $f$.

\begin{theorem}
\label{Th:Poincare}Let $\Sigma=\overline{p_{i}q_{i}}\subseteq\mathbb{D}%
\subseteq\mathbb{R}^{2}$ be a computable line segment which defines a
cross-section for a periodic orbit $\gamma$ of (\ref{ODE_Main}), where $f\in
C^{1}(\mathbb{D})$. Suppose also that the Poincar\'{e} map $P:\Sigma
\rightarrow\Sigma$ is defined for all all points $x\in\Sigma$. Then $P$ and
$DP$ are computable from $f$ of (\ref{ODE_Main}).
\end{theorem}

\begin{proof}
The techniques from \cite{GRZ18} and \cite{Tuc02a} are used to prove the computability of $P$ and $DP$, respectively.


We begin by showing that $P$ is computable. Since $\Sigma=\overline{p_{i}%
q_{i}}$ is a cross-section on an approximation $\widehat{C}_{i}$ of some periodic
orbit(s), the flow of (\ref{ODE_Main}) crosses this section transversaly. This
implies that for any point $x\in\Sigma$, the angle $\measuredangle
(f(x),\Sigma)$ between $f(x)$ and $\Sigma$ is nonzero, $\measuredangle
(f(x),\Sigma)>0$. Let $\theta=\frac{\min_{x\in\Sigma}\measuredangle
(f(x),\Sigma)}{2}>0$. Then $\theta$ is computable from $f$. Furthermore, by
continuity of $f$, there exists some $\varepsilon>0$ such that%
\[
\min_{x\in B(\Sigma,\varepsilon)}\measuredangle(f(x),\Sigma)>\theta>0
\]
where $B(\Sigma,\varepsilon)=\{x\in\mathbb{D}:\left\Vert y-x\right\Vert
\leq\varepsilon$ for some $y\in\Sigma\}$ contains no zeros of $f$. Let
\[
\alpha=\min_{x\in B(\Sigma,\varepsilon)}\left\Vert f(x)\right\Vert ,\text{
\ \ }\beta=\max_{x\in B(\Sigma,\varepsilon)}\left\Vert f(x)\right\Vert
\]
Since $B(\Sigma,\varepsilon)$ is compact and
contains no zero of $f$, it follows that $\alpha,\beta>0$. It is convenient to view  $B(\Sigma,\varepsilon)$ as a
rectangle. Note that
$B(\Sigma,\varepsilon)\cap \widehat{C}_{i}$ is divided into two parts\ $B_{1}$ and
$B_{2}$ by the line passing through $p_{i}$ and $q_{i}$. Let us assume,
without loss of generality, that the flow passes from $B_{1}$ through $\Sigma$
and then moves through $B_{2}$ until it leaves $B(\Sigma,\varepsilon)$. A
simple analysis shows that the flow of (\ref{ODE_Main}) cannot take more than
$2\varepsilon/(\alpha\sin\theta)>0$ time units to cross $B(\Sigma
,\varepsilon)$ (the flow will have to cross this rectangle; but since the norm
of the orthogonal component is at least $\alpha\sin\theta$, this will be done
in time $2\varepsilon/(\alpha\sin\theta)$), but requires at least
$2\varepsilon/\beta>0$ time units to cross it (because the norm of the
orthogonal component is bounded by $\beta$). Therefore if $x$ and $y$ are
solutions of (\ref{ODE_Main}) with initial conditions $x(0)=x_{0}$ and
$y(0)=y_{0}$, with $x_{0},y_{0}\in B(\Sigma,\varepsilon)\cap \widehat{C}_{i}$, then
$x(t)$ and $y(t)$ leave $B(\Sigma,\varepsilon)$ at times $t_{x},t_{y}%
\in\lbrack0,2\varepsilon/(\alpha\sin\theta)]$, respectively.

Now take some rational $\varepsilon_{0}>0$ satisfying $\varepsilon_{0}\leq
\min\{\epsilon\alpha\sin\theta/2,\varepsilon\}$, where $\epsilon>0$ is a
rational yet to be defined. In particular, this implies that the time for the
flow in $\widehat{C}_{i}$ to cross $B(\Sigma,\varepsilon_{0})\subseteq B(\Sigma
,\varepsilon)$ is bounded by%
\begin{equation}
2\varepsilon_{0}/(\alpha\sin\theta)\leq\frac{2\epsilon\alpha\sin\theta
}{2\alpha\sin\theta}\leq\epsilon\label{Eq:aux}%
\end{equation}

Let $\bar{B}_{1}=B_{1}\cap B(\Sigma,\varepsilon_{0})$, $\bar{B}_{2}=B_{2}\cap
B(\Sigma,\varepsilon_{0})$, and $\delta=\varepsilon_{0}/(2\beta)$. Next
consider the sequence of iterates $\phi_{t_{i}}(x)$, where $x\in\Sigma$, $0<t_{i+1}-t_{i}%
\leq\delta$, and $\{t_{i}\}_{i\in\mathbb{N}}$ is computable.
Since the flow $\phi_{t}(x)$ takes at least $2\delta$ time units
to cross each band $\bar{B}_{1}$ or $\bar{B}_{2}$, we are certain that
$\phi_{t_{1}}(x),\phi_{t_{2}}(x)\in\bar{B}_{2}$ when the flow first leaves
$\Sigma$ from $x$ and that there is some $k>0$ such that $\phi_{t_{k}}%
(x),\phi_{t_{k+1}}(x)\in\bar{B}_{1}$ with $t_{k},t_{k+1}>t_{1}>0$. Note that,
at most, only one of the iterates $\phi_{t_{k}}(x),\phi_{t_{k+1}}(x)$ are on
$\partial\bar{B}_{1}$, so that at least one of the iterates $\phi_{t_{k}%
}(x),\phi_{t_{k+1}}(x)$ is in the interior of $\bar{B}_{1}$. Since $\bar
{B}_{1}$ is compact and computable, so does the closure of its
complement. Thus, we can decide whether one of the iterates $\phi_{t_{k}}%
(x),\phi_{t_{k+1}}(x)$ is on $\bar{B}_{1}$ or on the closure of its complement
(the problematic case to detect is when one of these iterates is on
$\partial\bar{B}_{1}$, and that's why we always use two iterates, since this
ensures that at least one of the iterates will not be on the boundary of
$\bar{B}_{1}$). If we conclude that one iterate is on the closure of the
complement of $\bar{B}_{1}$, then we skip this iteration and increment $k$
until (and that will eventually happen) it reaches the first $k$ for which we
conclude that $\phi_{t_{k}}(x)$ belongs to $B(\Sigma,\varepsilon_{0})$ and
thus also necessarily to $\bar{B}_{2}$ and then return $\phi_{t_{k}}(x)$.

Now we turn to find some sufficiently small $\epsilon>0$ such that
$\phi_{t_{k}}(x)$ is close enough to the real value $P(x)$.
Assume that we need to compute $P(x)$ with accuracy $2^{-j}$ for some $j>0$.
Recall that the time needed for the flow in $\widehat{C}_{i}$ to cross $B(\Sigma
,\varepsilon_{0})$ is bounded by (\ref{Eq:aux}). Hence, the time it takes to
cross $\bar{B}_{1}$ until reaching $\Sigma$ is bounded by $\epsilon/2$. Let
$M$ be a rational such that $M\geq\sup_{x\in\bar{B}_{1}}\left\Vert
f(x)\right\Vert $, which is computable from $f$. Then $\left\Vert
P(x)-\phi_{t_{k}}(x)\right\Vert \leq M\epsilon/2$. Furthermore, since we
usually cannot compute $\phi_{t_{k}}(x)$ exactly, but only an approximation
$\bar{\phi}_{t_{k}}(x)$ of it, it follows that if we take $\epsilon\leq2^{-j}/M$
and compute $\bar{\phi}_{t_{k}}(x)$ with accuracy bounded by $M\epsilon/2$, then we
have%
\begin{align*}
\left\Vert P(x)-\bar{\phi}_{t_{k}}(x)\right\Vert  &  \leq\left\Vert
P(x)-\phi_{t_{k}}(x)\right\Vert +\left\Vert \phi_{t_{k}}(x)-\bar{\phi}_{t_{k}%
}(x)\right\Vert \\
&  \leq M\epsilon/2+M\epsilon/2\\
&  \leq M\epsilon\\
&  \leq2^{-j}%
\end{align*}
It remains to show that $DP$ is computable. The proof essentially follows
along the lines of \cite[Section 5]{Tuc02a} and uses the fact just shown above
that $P$ is computable. We may assume that the
(computable) cross-section $\Sigma$ is parallel to the vertical axis. If the asumption is not true, a (computable) change of basis will result in the desirable case.


First we note that if $\phi(t,\bar{x})$ denotes the solution of
(\ref{ODE_Main}) with initial condition $x(0)=\bar{x}$ at time $t$, and if
given some $x\in\Sigma$, $\tau(x)>0$ denotes the first time where the
trajectory starting on $x\in\Sigma$ will hit $\Sigma$ again, we have
$P(x)=\phi(\tau(x),x)$. Let us now calculate the partial derivatives of
$P=(P_{1},P_{2})$. We have%
\begin{align}
\frac{\partial P_{i}}{\partial x_{j}}(x)  &  =\frac{\partial}{\partial x_{j}%
}\phi_{i}(\tau(x),x)\nonumber\\
&  =\frac{\partial\phi_{i}(\tau(x),x)}{\partial t}\frac{\partial\tau
(x)}{\partial x_{j}}+\frac{\partial\phi_{i}(\tau(x),x)}{\partial x_{j}%
}\nonumber\\
&  =f_{i}(\phi(\tau(x),x))\frac{\partial\tau(x)}{\partial x_{j}}%
+\frac{\partial\phi_{i}(\tau(x),x)}{\partial x_{j}}\nonumber\\
&  =f_{i}(P(x))\frac{\partial\tau(x)}{\partial x_{j}}+\frac{\partial\phi
_{i}(\tau(x),x)}{\partial x_{j}} \label{Eq:derivatives_Poincare}%
\end{align}
for $i=1,2$. To obtain the partial derivatives of $\tau(x)$, we note that
$P_{1}(x)$ is constant. Hence, by (\ref{Eq:derivatives_Poincare})
\[
0=\frac{\partial P_{1}}{\partial x_{j}}(x)=f_{1}(P(x))\frac{\partial\tau
(x)}{\partial x_{j}}+\frac{\partial\phi_{1}(\tau(x),x)}{\partial x_{j}}.
\]
Solving for $\partial\tau(x)/\partial x_{j}$, we get%
\[
\frac{\partial\tau(x)}{\partial x_{j}}=-\frac{\partial\phi_{1}(\tau
(x),x)}{\partial x_{j}}\frac{1}{f_{1}(P(x))}%
\]
(note that $f_{1}(P(x))\neq0$ as the flow of (\ref{ODE_Main}) is transverse to
$\Sigma$). Replacing this last expression into (\ref{Eq:derivatives_Poincare})
(note that we only need to compute the partial derivatives of $P_{2}$ as the
partial derivatives of $P_{1}$ are zero, as we have seen), we get%
\begin{equation}
\frac{\partial P_{2}}{\partial x_{j}}(x)=-f_{2}(P(x))\frac{\partial\phi
_{1}(\tau(x),x)}{\partial x_{j}}\frac{1}{f_{1}(P(x))}+\frac{\partial\phi
_{2}(\tau(x),x)}{\partial x_{j}} \label{Eq:Poincare_deriv}%
\end{equation}
The only element still missing in order to compute the partial derivatives of
$P_{2}$ is the computation of the partial derivatives of $\phi_{l}$, for $l=1,2$. This can
be achieved as follows. From (\ref{ODE_Main}) we get that $\left(  \phi
_{k}(t,x)\right)  ^{\prime}=f_{k}(\phi(t,x))$ for $k=1,2$. By applying partial
derivatives to both sides and switching the order of differentiation on the
left-hand side ($\phi$ is $C^{2}$ and thus this operation will not affect the
result), we get%
\[
\frac{d}{dt}\frac{\partial\phi_{l}(t,x)}{\partial x_{j}}=\sum_{i=1}^{2}%
\frac{\partial f_{l}}{\partial x_{i}}(\phi(t,x))\frac{\partial\phi_{i}%
}{\partial x_{j}}(t,x)\text{, }l=1,2.
\]
In matrix form this can be written as
\[
\frac{d}{dt}D\phi(t,x)=Df(\phi(t,x))D\phi(t,x)
\]
which is a linear ODE with the initial condition $D\phi(0,x)=I$ (note that
$\phi(0,x)=x$ is the identity map and that the partial derivatives are only
taken in order to $x$). The solution of this initial-value problem can be
computed as the solution of an ODE, for an arbitrary amount of time as in
\cite{GZB07}. Moreover the time $\tau(x)$ used in (\ref{Eq:Poincare_deriv})
can also be computed as in \cite{GRZ18}. Indeed, from the above arguments, one
can conclude that to compute $\tau(x)$ with accuracy $2^{-n}$ it suffices to
take $\epsilon< 2^{-n}$ and return $t_{k}$ as above. This proves the theorem.
\end{proof}

\section{Proof of Theorem B -- Putting it all
together\label{Sec:ThmBProofFinal}}
Let us assume that $Zero_{\epsilon}(f)$ does not include any (very small) periodic orbit. This can be ensured, for example, by using Theorem \ref{Prop:HartmanGrobman}, the computable Hartman-Grobman theorem, to compute some $\epsilon>0$ such that if $x_0$ is an equilibrium point, then $B(x_0,\epsilon)$ is inside the neighborhood computed by the computable version of the Hartman-Grobman Theorem. In this neighborhood there are no periodic orbits, since the flow is conjugated to a linear flow there, and therefore periodic orbits can only exist on the doughnut-like sets computed in previous sections.

These doughnut-like sets thus include all the periodic orbits with arbitrarily high
accuracy. We have also shown that on each of these doughnut-like sets we can
define a cross-section and compute a Poincar\'{e} map, $P$, there. Moreover,
using the technique of Section \ref{Sec:Zeros} we can compute the number of
zeros of $P(x)-x$, i.e. the number of fixed points of $P$. But the number of
fixed points is equal to the number of periodic orbits. Hence Theorem B is proved.

\section{Connections with Hilbert's 16th problem}\label{Sec:Hilbert16}

The algorithm constructed in Theorem B computes the positions and
the exact number of the periodic orbits for every vector field in
$SS_2$. Let us call the algorithm \texttt{Algo}. As a by-product,
\texttt{Algo} produces a computable function $\phi: SS_2 \to
\mathbb{N}$, where $\phi (f) =$ the number of the periodic orbit(s)
of $f$ on $\mathbb{D}$. The pre-images of $\phi$ decompose $SS_2$
into mutually disjoint open connected components $C_j$, $j\in
\mathbb{N}$. Let $\mathcal{P}$ (respectively,
$\mathcal{P}_{\mathbb{Q}}$) denote the set of all polynomials
(respectively, all polynomials with rational coefficients) defined
on $\mathbb{D}$. Then $\mathcal{P}_{\mathbb{Q}}$ is dense in
$C^1(\mathbb{D})$. Since $SS_2$ is an open subset of
$C^1(\mathbb{D})$, it follows that every $C_j$ contains (infinitely
many) polynomials from $\mathcal{P}_{\mathbb{Q}}$. Hence, if there
is an algorithm, say \texttt{A-for-H} (an algorithm for Hilbert's 16th
problem), that computes an upper bound $u(n)$ for the numbers of
periodic orbits of polynomials of degree $n$ in $SS_2 \bigcap
\mathcal{P}_{\mathbb{Q}}$, then $u(n)$ is also an upper bound for
the numbers of periodic orbits of all polynomials of degree $n$ in
$SS_2$. This result would provide an affirmative answer to Hilbert's
16th problem restricted to polynomials in $SS_2$, that is, the
structurally stable polynomials on $\mathbb{D}$. Since $SS_2$ is an
open dense subset in $C^1(\mathbb{D})$, a property true on $SS_2$ is
typical and generic.

Whether an \texttt{A-for-H} algorithm exists is an open problem. Theorem C
indicates that one may not be able to construct algorithms of an
\texttt{A-for-H} nature but for computing sharp upper bounds over certain
classes of polynomial systems. On the other hand, we construct in
the following an \texttt{A-for-H} algorithm over $SS_2\bigcap \mathcal{P}$
that works relative to the Halting problem. In other words, it is
possible to devise a Turing machine (an algorithm) that solves
Hilbert's 16th problem over $SS_2\bigcap \mathcal{P}$, provided that
the Halting problem
\[
HALT=\{(M,i):\text{the Turing machine $M$ halts with input }i\}
\]
 is given as an oracle. As usual, we use the notation $C\leq HALT$ to denote that the problem $C$ is solvable by an algorithm that works relative to the halting problem (or, more generally, relative to another problem).

We begin by effectively listing the polynomial systems in
$\mathcal{P}_{SS_2}=SS_2\bigcap \mathcal{P}_{\mathbb{Q}}$. First we note that
$\mathcal{P}_{\mathbb{Q}}$ can be enumerated as
$\mathcal{P}_{\mathbb{Q}}=\{ P_j\}_{j=1}^{\infty}$. Let $M$ be a
Turing machine that on input $k$ computes the first $k$ steps of \texttt{Algo}$(P_1)$,
\texttt{Algo}$(P_2)$, $\ldots$, \texttt{Algo}$(P_k)$. $M(k)$ outputs $P_j$ if $1\leq j\leq k$, \texttt{Algo}$(P_j)$
halts in $\leq k$ steps, and $P_j$ is not the output of $M(l)$ for
$l < k$; otherwise, $M(k)$ outputs the empty set. It is clear that
$M$ lists recursively all polynomial systems in $SS_2\bigcap
\mathcal{P}_{\mathbb{Q}}$. For simplicity, we use $\{ P_j\}$ to
denote this computable sequence. Then $\{P_j\}$ is a subset of
$\mathcal{P}\bigcap SS_2$.

Let $\mathcal{A}\subseteq\mathcal{P}$ and let
$Hilbert_{16}(\mathcal{A})$ be the problem of solving the Hilbert's
16th problem over the set $\mathcal{A}$, i.e.~  computing an upper
bound for the numbers of periodic orbits of all elements of
$\mathcal{A}$ of degree $n$, where $n$ is given as an input.

\begin{theorem}\label{th:Hilbert-halt}
$Hilbert_{16}(\mathcal{P}\cap SS_2)\leq HALT$.
\end{theorem}

\begin{proof}
Consider the Turing machines $N_{n}$, $n\in \mathbb{N}\setminus \{
0\}$, defined as follows, where $k\in\mathbb{N}$ is the input of
$N_{n}$:
\begin{enumerate}
\item Let $i=1$.

\item Consider the polynomial vector fields with rational coefficients $P_1, P_2,\ldots, P_i$
(in $\mathcal{P}_{SS_2}$) and retain only those which have degree
$n$. For each of the retained polynomial vector fields, use it as an
input to the algorithm \texttt{Algo}, but simulate only $i$ steps of
the algorithm. If in any of these computations the algorithm
\texttt{Algo} stops and returns a number greater than or equal to
$k$, then $N_{n}$ stops its computation. Otherwise, $i$ is
incremented and Step 2 is repeated.
\end{enumerate}

It is clear that $N_n$ halts on input $k$ only if there is a vector
field in $\{ P_j\}$ with $k$ or more periodic orbits. Otherwise
$N_n$ will not halt with input $k$.

Now consider the Turing machine that has oracle access to $HALT$
defined as follows:  on input $n$ (the degree of the polynomials),
\begin{enumerate}
\item Set $k=1$.

\item Using the oracle to decide whether $(N_n,k)\in HALT$. If the answer is positive, then increment $k$ and repeat this step.
Otherwise return $k$.
\end{enumerate}

It is readily seen that the output of this Turing machine on input
$n$ would provide a (sharp) upper bound for the numbers of periodic
orbits of polynomials of degree $n$ in $\{ P_j\}$, provided that
such a bound exists. In the case that the upper bound doesn't exist,
the Turing machine won't halt.

As we mentioned at the beginning of this section, an upper bound for
the numbers of periodic orbits of polynomials of degree $n$ in $\{
P_j\}$ is also an upper bound of the same nature of polynomials of
degree $n$ in $\mathcal{P}\bigcap SS_2$.
\end{proof}

The result above can be generalized from the compact domain
$\mathbb{D}$ to $\mathbb{R}^2$. Now let $\mathcal{P}$ denote the set
of polynomial vector fields defined over the \textit{whole} plane
$\mathbb{R}^2$, and let $\mathcal{P}_m$ denote the set of polynomial
vector fields of degree $m$. Then we have:

\begin{theorem}\label{Th:Hilbert-plane}
There are dense subsets $A,B$ of $\mathcal{P}$, where $B$ is also
open, such that $B\subseteq A\subseteq\mathcal{P}$ and
$Hilbert_{16}(A)\leq HALT$.
\end{theorem}

\begin{proof}
The proof of this result is similar to the previous theorem, but
some adaptations are needed. In \cite[Theorem 4, p.~327]{Per01} it
is mentioned that the set $\mathcal{B}_m$ of polynomial vector
fields of degree $m$ which are structurally stable on $\mathbb{R}^2$
under $C^1$-perturbations is open and dense in $\mathcal{P}_m$.
Furthermore, every system in $\mathcal{B}_m$ only has hyperbolic
equilibrium points/periodic orbits without any saddle connections
even if it has saddles at infinity. Let $\mathcal{A}_m$ be the
systems in $\mathcal{P}_m$ which (i) only have hyperbolic
equilibrium points and periodic orbits, and (ii) have no ``finite''
saddle connections but may have saddle connections involving saddles
at infinity. It is clear that $\mathcal{B}_m\subseteq
\mathcal{A}_m\subseteq \mathcal{P}_m$; thus  $\mathcal{A}_m$ is also
dense in $\mathcal{P}_m$. Let
$\mathcal{A}=\cup_{m\in\mathbb{N}}\mathcal{A}_m$ and
$\mathcal{B}=\cup_{m\in\mathbb{N}}\mathcal{B}_m$. Then $\mathcal{B}$
is open, $\mathcal{A}$ is dense in $\mathcal{P}$, and $\mathcal{B}\subseteq
\mathcal{A}\subseteq \mathcal{P}$.

We now show that $Hilbert_{16}(\mathcal{A})\leq HALT$. We note that
the algorithm \texttt{Algo} can be applied to the systems defined on
$\mathbb{D}_n=\{ x\in \mathbb{R}^2: \| x\|\leq n\}$, $n\in
\tilde{\mathbb{N}}=\mathbb{N}\setminus \{ 0\}$, by rescaling the
systems to $\mathbb{D}$. More precisely, let \texttt{AlgoGen} be the
algorithm that, on input $(n,f)$,  rescales the vector field $f$
defined on $\mathbb{D}_n$ to $\mathbb{D}_1=\mathbb{D}$ and then
applies algorithm \texttt{Algo} to the rescaled system, where $n\in
\tilde{\mathbb{N}}$ and $f\in C^1(\mathbb{D}_n)$.

Note, however, that we have assumed that the flow defined by
\eqref{ODE_Main} points inwards across the boundary of $\mathbb{D}$
when using the algorithm \texttt{Algo} and thus a similar
requirement seems to be in order for \texttt{AlgoGen}. However, that
requirement is only needed because of the possibility that there
might be trajectories (which e.g.~might be part of a periodic orbit)
tangent to the boundary of $\mathbb{D}$ and it is computationally
impossible to detect that (using the slope of the flow, this is in
essence equivalent to determining if two real numbers are equal).
However, for our case of \texttt{AlgoGen} applied to $\mathbb{D}_n$,
we could instead apply it to $\mathbb{D}_{n+1}$ with the following
adaptation: if a trajectory starting in $\mathbb{D}_n$ enters the
region $\mathbb{D}_{n+1}\setminus\mathbb{D}_n$, then we remove that
trajectory (concluding that the flow will leave $\mathbb{D}_n$ and
accordingly it will be dealt with later in a larger region
$\mathbb{D}_l$ for some $l>n$)
and restart the computation starting at another ``pixel'' of
$\mathbb{D}_n$. Afterwards we count the number of periodic orbits of
\eqref{ODE_Main} which are clearly inside $\mathbb{D}_n$ (if they
intersect the boundary of $\mathbb{D}_n$, they are not counted).
Note that if \eqref{ODE_Main} is defined on the plane and has only a
finite number of hyperbolic equilibria and periodic orbits with no
saddle connections, then \texttt{AlgoGen} will always halt.
Moreover, if all periodic orbits of \eqref{ODE_Main} are inside
$\mathbb{D}_n$, then \texttt{AlgoGen}$(l,f)$ will return either the
number $\alpha$ of periodic orbits of \eqref{ODE_Main} when $l\geq
n+1$ or a number $\leq \alpha$ when $l\leq n$.

Let $a(j, m)$, $j , m\in \tilde{\mathbb{N}}$, be a computable
sequence listing all rational polynomials on $\mathbb{R}^2$ of
degree $m$, and let $N_{m}$ be the Turing machine defined as
follows: on input $k,l\in\mathbb{N}$,
\begin{enumerate}
\item set $i=1$;

\item simulate $i$ steps of the algorithm
\texttt{AlgoGen}$(l,\cdot)$ on each of $a(1, m), \ldots, a(i, m)$. If
in any of these computations, the algorithm \texttt{AlgoGen} stops
and returns a number greater than or equal to $k$, then $N_{m}$
stops the computation. Otherwise $i$ is incremented and Step 2 is
repeated.
\end{enumerate}

From the design it is clear that $N_{m}$ halts with input $k,l$ iff
there is a polynomial vector field in $\mathcal{P}_{m}\cap
SS_2(\mathbb{D}_l)$ that has $k$ or more periodic orbits. This fact
follows from Peixoto's theorem applied to $\mathbb{D}_l$ because
Peixoto's theorem implies that $SS_2(\mathbb{D}_l)$ is the union of
mutually disjoint open connected components, and all vector fields
in the same component have the same number of periodic orbits. Thus,
if there is a polynomial vector field $p \in \mathcal{P}_{m}\cap
SS_2(\mathbb{D}_l)$ with $k$ periodic orbits, then the open
component in which $p$ lies must contain a polynomial vector field
$q$ of degree $m$ with rational coefficients, for the set of all
polynomial vector fields with rational coefficients is dense in
$C^1(\mathbb{D}_l)$. Since $q$ is in the same component as of $p$,
it follows that $q$ also has $k$ periodic orbits.
Hence, $N_m$ will halt with input $k,l$ when it lists $q$ and is
allowed for a sufficiently many steps of \texttt{AlgoGen}$(l,\cdot)$
with input $q$. On the other hand, if there is no vector field in
$\mathcal{P}\cap SS_2(\mathbb{D}_l)$ with $k$ or more periodic
orbits, then $N_{m}$ will not halt with input $k,l$. We also note
that if a rational polynomial $a(i,m)$ has $k$ (hyperbolic) periodic
orbits, then there is some $l\in\mathbb{N}$ such that all periodic
orbits of $a(i, m)$ is contained in $\mathbb{D}_l$; hence $N_{m}$
would halt with input $k,l+1$.

Now consider the Turing machine with oracle access to $HALT$, which
operates on input $m$ as follows, where $\left\langle
\cdot,\cdot\right\rangle:\mathbb{N}^2\to\mathbb{N}$ is a computable
bijection, e.g.~like the one defined in \cite[p.~27]{Odi89}:
\begin{enumerate}
\item Let $i=\left\langle k,l\right\rangle =1$.

\item Using the oracle, decide whether $(N_m,(k,l))\in HALT$ (note that $k$ and $l$
can be computed from $i$). If the answer is positive increment $k$
and repeat this step. Otherwise return $k$.
\end{enumerate}

In is readily seen that the output of this Turing machine on input
$m$ returns a (sharp) upper bound (if it exists) on the number of
periodic orbits for elements of $A_m$.
\end{proof}

\section{Conclusion}

In this paper, we have shown that, in general, one cannot compute
the number of periodic orbits that a polynomial ODE (\ref{ODE_poly})
can have. Even sharp upper bounds on the number of periodic orbits
cannot, in general, be computed for subfamilies of polynomial
systems.

On the other hand, we have shown that the exact number of periodic
orbits can be computed uniformly for all structurally stable planar
dynamical systems (\ref{ODE_Main}) defined on the unit ball, as well
as the limit set $NW(f)$. The algorithm computing the exact number
of periodic orbits also portrays them with any precision one wishes
to have.

We conclude the paper with a question: What is the computational
complexity of computing the number of periodic orbits (or the limit
set $NW(f)$) when (\ref{ODE_Main}) is structurally stable?  In other
words, what computational resources (e.g.~in terms of time or
space/memory) are required to compute the number of periodic orbits?
It is known that when solving an ODE (\ref{ODE_Main}) with a
polynomial-time computable vector field $f$, the complexity of
computing the solution at $t=1$ can be arbitrarily high if $f$ does
not satisfy a Lipschitz condition \cite[p. 469]{Mil70}, and
$PSPACE$-complete when $f$ satisfies a Lipschitz condition or is of
class $C^{k}$ \cite{Kaw10}, \cite{KORZ14}. So we might expect that
the complexity of computing the number of period orbits of a general
ODE is at least as high as those bounds. On the other hand,
hyperbolicity often provides some degree of regularity which can be
exploited to lower the complexity upper bounds, such as in
e.g.~\cite{Bra05}, \cite{Ret05}. Therefore, it could also be the
case that the hyperbolicity of the periodic orbits might be
exploited to obtain smaller complexity upper bounds. That would be an
interesting question for further work.\medskip

\textbf{Acknowledgments.}Daniel Gra\c{c}a was partially funded by FCT/MCTES through national funds and when applicable co-funded EU funds under the project UIDB/50008/2020. \includegraphics[width=3.5mm]{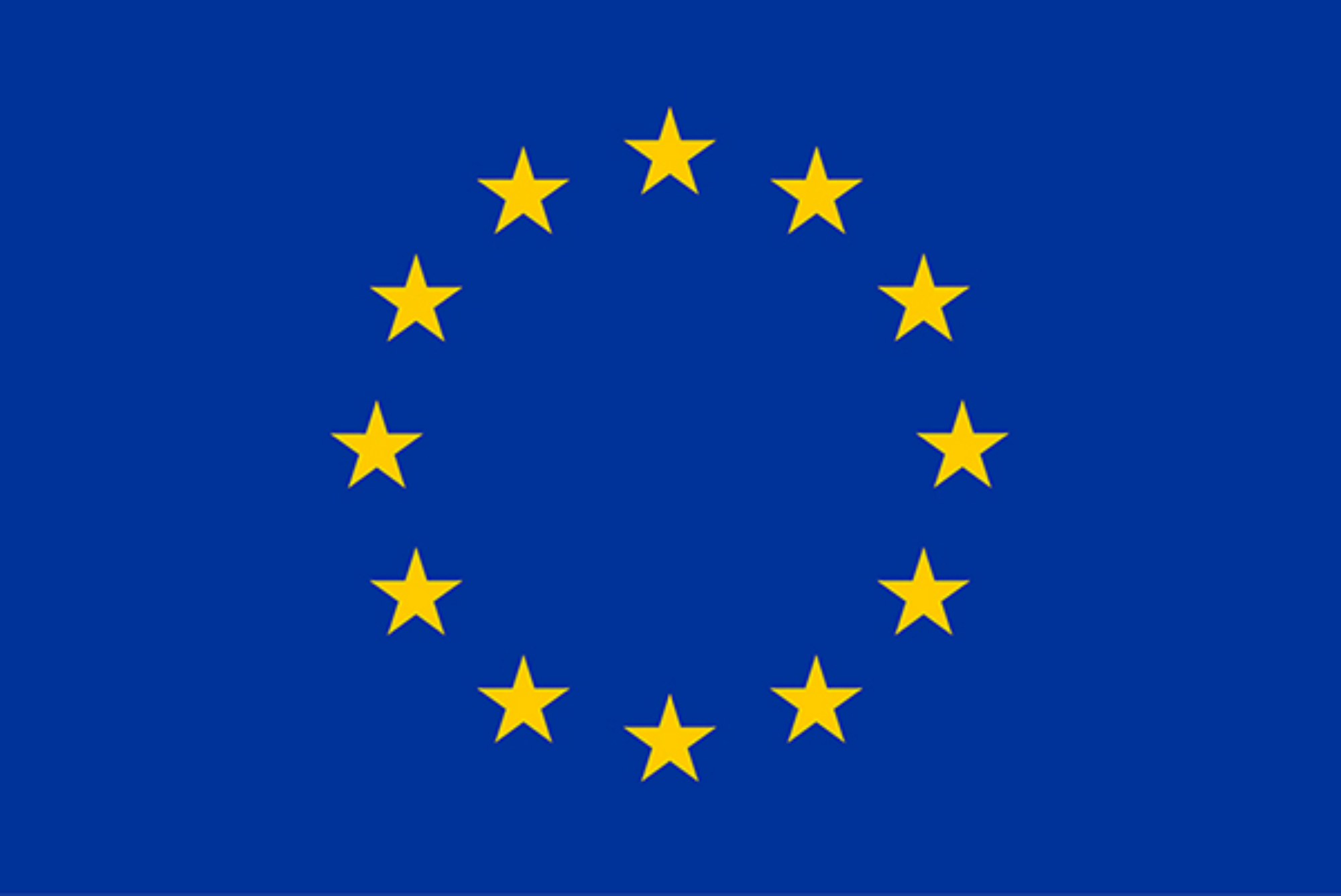} This project has received funding from the European Union's Horizon 2020 research and innovation programme under the Marie Sk\l{}odowska-Curie grant agreement No 731143.

\bibliographystyle{plain}
\bibliography{ContComp}

\appendix

\section{Proofs of results about dynamical systems}\label{Sec:AppendixDynSys}

\textbf{Proof of Lemma \ref{Lemma:equilibrium_point_cycle}. }Notice that the
dynamics inside the region $\mathcal{D}$ delimitated by $\gamma$ never leave
this region for all $t\in\mathbb{R}$. Therefore, if $\gamma$ is repelling, the
trajectories move away from $\gamma$ towards an attractor $\varsigma_{1}$ as
$t\rightarrow+\infty$. If $\gamma$ is attracting, the trajectories move away
from $\gamma$ towards a repeller $\varsigma_{1}$ as $t\rightarrow-\infty$. By
Theorem \ref{Th:GeneralPeixoto}, $\varsigma_{1}$ can only be an equilibrium
point or a periodic orbit. If $\varsigma_{1}$ is an equilibrium point, we are
done. If $\varsigma_{1}$ is a periodic orbit, we can repeat the procedure to
obtain a new limit object $\varsigma_{2}$. Since there is only a finite number
of periodic orbits, repeating this procedure will eventually yield an
equilibrium point $\varsigma_{m}$ inside $\mathcal{D}$, thus showing the result.

\begin{theorem}
\label{Th:Non-wandering}Let (\ref{ODE_Main})\ be defined on a compact set
$K\subseteq\mathbb{R}^{n}$. Then for every $\varepsilon>0$ there exist some
$\delta>0$ and $T>0$ such that for any $x\in K$, if $d(x,NW(f))\geq
\varepsilon$, then $\phi_{t}(B(x,\delta))\cap B(x,\delta)=\varnothing$ for
every $t\geq T$, where $B(x,\delta)=\{y\in K:\left\Vert x-y\right\Vert
<\delta\}.$
\end{theorem}

\begin{proof}
By definition, if $x\notin NW(f)$, this means that there is some neighboorhood
$U_{x}$ and some $T_{x}$ such that $\phi_{t}(U_{x})\cap U_{x}=\varnothing$ for
every $t\geq T_{x}$.

Since $NW(f)$ is closed, the set $A_{\varepsilon}=\{y:d(y,NW(f))<\varepsilon
\}$ is open. This is because $NW(f)$ is compact and therefore there exists a
dense sequence of points $\{x_{j}\}_{j\in\mathbb{N}}$ in $NW(f)$, which
implies that $A_{\varepsilon}=\cup_{j\mathbb{\in N}}\mathring{B}%
(x_{j},\varepsilon)$, where $\mathring{B}(x_{j},\varepsilon)=\{y\in
\mathbb{R}^{n}:\left\Vert y-x_{j}\right\Vert <\varepsilon\}$ is the interior
of $B(x_{j},\varepsilon)$. Since every open ball $\mathring{B}(x_{j}%
,\varepsilon)$ is an open set and a countable union of open sets is also an
open set, then $A_{\varepsilon}$ must be open.

We then conclude that $\{U_{x}\}_{x\in K-NW(f)}\cup A_{\varepsilon}$ defines
an open covering of $K$. Since $K$ is compact, this implies that there are
$x_{1},\ldots,x_{k}$ such that $\{U_{x_{l}}\}_{1\leq l\leq k}\cup
A_{\varepsilon}$ is a cover of $K$. Moreover, by the Lebesgue number Lemma
(see e.g.\ \cite[Lemma 27.5 on p.~175]{Mun00}), there exists some $\delta>0$
such that for every $x\in K$, $B(x,\delta)$ is contained in some element of
the covering $\{U_{x_{l}}\}_{1\leq l\leq k},A_{\varepsilon}$. In particular,
if $d(x,NW(f))\geq\varepsilon$, $B(x,\delta)$ cannot be contained in
$A_{\varepsilon}$. Instead it must be $B(x,\delta)\subseteq U_{x_{l}}$ for
some $l$. Take $T=\max_{1\leq l\leq k}T_{x_{l}}$. Since $\phi_{t}(U_{x_{l}%
})\cap U_{x_{l}}=\varnothing$ for all $t\geq T\geq T_{x_{l}}$, it must be
$\phi_{t}(B(x,\delta))\cap B(x,\delta)=\varnothing$ for all $t\geq T$. This
proves the lemma.
\end{proof}

In the remaining of this paper, if $x\in\mathbb{R}^{n}$ and $A\subseteq
\mathbb{R}^{n}$, we take%
\[
d(x,A)=\inf_{y\in A}\left\Vert x-y\right\Vert .
\]
The following lemma is a consequence of Propositions
\ref{Prop:ConvEquilibrium} and \ref{Prop:ConvPeriodic}.

\begin{lemma}
\label{Lemma:basin}Let (\ref{ODE_Main})\ be a structurally stable system
defined on a compact set $K\subseteq\mathbb{R}^{2}$. Then there is some
$\varepsilon>0$ such that for any attracting periodic orbit $\gamma$ or
equilibrium point $x_{0}$, there exists a neighborhood $U_{A}$ (the
\emph{basin of attraction}) such that $d_{H}(\overline{U}_{A},A)\geq
\varepsilon$, where $A=\gamma$ or $A=\{x_{0}\}$ depending on the case, where
any trajectory starting from a point of $U_{A}$ will converge towards $A$.
Moreover, given $\delta>0$ with $\delta<\varepsilon$, there exists some
$T_{\delta}$, independent of the attractor $A$, such that if $x(t)$ is a
solution of (\ref{ODE_Main}) such that $x(t_{1})\in U_{A}$, then
$d(x(t_{1}+T_{\delta}),A)<\delta$.
\end{lemma}

\begin{proof}
The existence of the neighborhood $U_{A}$ comes from Propositions
\ref{Prop:ConvEquilibrium} and \ref{Prop:ConvPeriodic}. Moreover, $A$ is
always a closed set, which implies that $d_{H}(\overline{U}_{A},A)>0$
($\overline{U}_{A}$ denotes the closure of $U_{A}$). Indeed, since $A\subset
U_{A}$, it follows that $A\bigcap(K-U_{A})=\emptyset$, and thus there exists
some $\eta>0$ such that $d_{H}(A,K-U_{A})\geq\eta$. We claim that
$d_{H}(A,\overline{U}_{A})\geq\eta$. For any $x\in\overline{U}_{A}$, if
$x\not \in U_{A}$, then $d(x,A)\geq d_{H}(K-U_{A},A)\geq\eta$; if
$x\in\overline{U}_{A}-U_{A}$, then there exists a sequence $p_{i}$ of points
in $K-U_{A}$ such that $p_{i}\rightarrow x$; thus $d(x,A)=\lim_{i\rightarrow
\infty}d(p_{i},A)\geq\eta$. From Propositions \ref{Prop:ConvEquilibrium} and
\ref{Prop:ConvPeriodic}, we know that there is some $\varepsilon_{A}>0$ such
that on
\[
A_{\varepsilon_{A}}=\left\{  x\in K:d(x,A)<\varepsilon_{A}\right\}
\]
the convergence to $A$ is exponentially fast, i.e.\ there are constants
$M,\alpha>0$ such that%
\[
d(x,A)<\varepsilon_{A}\text{ implies }d(\phi_{t}(x),A)<Me^{-\alpha t}%
\]
Let $\mathcal{A}$ be the set consisting of all attractors of (\ref{ODE_Main}).
By taking $\varepsilon=\min_{A\in\mathcal{A}}\{d_{H}(\overline{U}%
_{A},A),\varepsilon_{A}\}$, we then immediately conclude the lemma (note that
the number of attractors is finite).
\end{proof}

Using the above Lemma, we can adapt Theorem \ref{Th:Non-wandering} to prove
Theorem \ref{Cor:uniform_convergence_full}.

\textbf{Proof of Theorem \ref{Cor:uniform_convergence_full}. } We note that any
trajectory that starts on a point not in $NW(f)$ will have to converge to an
attracting equilibrium point or to an attracting periodic orbit. Let us call
such an attracting equilibrium point or attracting periodic orbit an
attractor. Due to Lemma \ref{Lemma:basin}, we know that the attractor is
hyperbolic and has a neighborhood with the property that each point of this
neighborhood converges exponentially fast to the attractor. Then given an
attractor $A$, there is some $\bar{\epsilon}>0$ such that if $x_{0}%
\in\mathbb{D}$ is such that $d(x_{0},A)\leq\bar{\epsilon}$, then the
trajectory starting at $x_{0}$ will converge, when $t\rightarrow+\infty$ (or
when $t\rightarrow-\infty$, in the case of repellers), exponentially fast to
$A$. Then we can consider the set%
\begin{equation}
A_{\bar{\epsilon}}=\mathcal{N}_{\bar{\epsilon}}(A)=%
{\displaystyle\bigcup\limits_{x\in A}}
\mathring{B}(x,\bar{\epsilon})=%
{\displaystyle\bigcup\limits_{x\in A}}
\{y\in\mathbb{D}:\left\Vert x-y\right\Vert <\bar{\epsilon}\}
\label{Eq:A_epsilon}%
\end{equation}
which works like an hyperbolic basin of attraction to $A$. Since $NW(f)$ only
has a finite number of connected components (see Theorem
\ref{Th:GeneralPeixoto}), we can take $\bar{\epsilon}$ to be the minimum of
all the $\bar{\epsilon}$'s for each particular $A$. Moreover, let us apply
Theorem \ref{Th:Non-wandering} with $\varepsilon=\min(\bar{\epsilon}%
,\epsilon)/2$, obtaining some values $T^{\ast}$ and $\delta$ such that the
conditions of this theorem hold. Let us also assume, without loss of
generality, that $\delta\leq\varepsilon/2$. Let us also take $T=16T^{\ast
}/\delta^{2}$. Let $Att_{\epsilon}$ be formed by the (finite) union of all
$A_{\varepsilon}$, where $A$ is an attractor. Notice that if a trajectory
enters a connected component of $Att_{2\varepsilon}$, then it will stay there.
Therefore, to prove the theorem it is enough to prove that any trajectory
starting on $\mathbb{D}-A_{\epsilon}$ will reach $A_{\epsilon}$ in time $\leq
T=16T^{\ast}/\delta^{2}$.

First let us find a finite number of rationals $p_{1},\ldots,p_{l}%
\in\mathbb{D}$ ($l\leq(4/\delta)^{2}=16/\delta^{2}$) such that%
\[
\mathbb{D}=B(0,1)\subseteq%
{\displaystyle\bigcup\limits_{i=1}^{l}}
B(p_{i},\delta).
\]
Let $x\in\mathbb{D}$. If $d(x,A)\leq2\varepsilon$, then the result is
obviously true. Let us hence suppose that $d(x,A)>2\varepsilon$. Therefore $x$
will belong to some ball $B(p_{i_{1}},\delta)$ which does not intersect
$A_{\varepsilon}$. Hence we can apply Theorem \ref{Th:Non-wandering} to
conclude that $\phi_{t}(x)\notin\mathbb{D}-B(p_{i_{1}},\delta)$ for all $t\geq T^{\ast}$.
If $\phi_{T^{\ast}}(x)\in A_{2\varepsilon}$, then the result is true.
Otherwise, assuming that $\phi_{T^{\ast}}(x)\in\mathbb{D}-B(p_{i_{1}},\delta
)$, we conclude that $\phi_{T^{\ast}}(x)\in B(p_{i_{2}},\delta)$. Hence, in a
similar manner, we conclude that $\phi_{t}(x)\notin\mathbb{D}-(B(p_{i_{1}},\delta)\cup
B(p_{i_{2}},\delta))$ for all $t\geq 2T^{\ast}$. We can continue this
procedure, \textquotedblleft exhausting\textquotedblright\ one ball
$B(p_{i},\delta)$ at a time. Since there are only $l\leq16/\delta^{2}$ such
balls, in time $16T^{\ast}/\delta^{2}=T$ we are sure to have exhausted all of
them. Hence, for $t\geq T$, $\phi_{t}(x)\in A_{2\varepsilon}$ for all $t\geq T^{\ast}$. This
proves the result.

\textbf{Proof of Lemma \ref{Lemma:unstable_convergence}. }Since $x_{0}$ is an
hyperbolic saddle point, $U_{x_{0}}$ is a 1-dimensional manifold in
$\mathbb{D}$, i.e. it is a curve which contains the point $x_{0}$.
Furthermore, by the stable manifold theorem (see e.g.~\cite[pp. 107--108]%
{Per01}), we also know that $\phi_{-t}(U_{x_{0}})\subseteq U_{x_{0}}$ for all
$t\geq0$. Since $U_{x_{0}}$ is a curve which contains $x_{0}$, that means that
$x_{0}$ divides $U_{x_{0}}$ into two 1-dimensional manifolds $\Gamma_{1}$ and
$\Gamma_{2}$ not containing $x_{0}$. Let us take $y_{1}\in\Gamma_{1}$. It is
well known that $\Gamma_{1}\subseteq\{\phi_{t}(y_{1}):t\in\mathbb{R}\}$. We
thus conclude that that if $\phi_{t}(y_{1})$ converges to an
attractor $\Omega_{1}(x_{0})$ as time increases, the same will happen for any
point in $\Gamma_{1}$. We notice that $\phi_{t}(y_{1})$ cannot converge to a
saddle point (i.e. $\Gamma_{1}$ cannot be part of the stable manifold of a
saddle point) because otherwise we would have a saddle connection, and this
cannot happen on structurally stable systems due to Peixoto's theorem. A
similar reasoning shows that any trajectory starting on a point of $\Gamma
_{2}$ will converge to an attractor $\Omega_{2}(x_{0})$ (note that it might be
$\Omega_{1}(x_{0})=\Omega_{2}(x_{0})$). This proves the result.\bigskip

To prove Theorem \ref{Thm:No_saddle}, we need the following definition (see
e.g. \cite[pp. 511-514]{HW95} for more details) which will be useful to show
the preliminary Lemma \ref{Lemma:Going_through_saddle}.

\begin{definition}
Consider the one-dimensional ODE $x^{\prime}=h(t,y)$, where $h$ is $C^{1}$. Then:

\begin{enumerate}
\item A $C^{1}$\emph{ function }$\alpha:I\subseteq\mathbb{R}\rightarrow
\mathbb{R}$, where $I$ is some interval, with the property that $\alpha^{\prime}(t)\leq h(t,\alpha(t))$ for all $t\in I$, is called a
\emph{lower fence};

\item A $C^{1}$\emph{ function }$\beta:I\subseteq\mathbb{R}\rightarrow
\mathbb{R}$, where $I$ is some interval, with the property that  $\beta^{\prime}(t)\geq h(t,\beta(t))$ for all $t\in I$, is called an
\emph{upper fence};

\item Consider a pair of functions $\alpha,\beta$, such that $\alpha$ is a
lower fence and $\beta$ is an upper fence on a common interval $I\subseteq
\mathbb{R}$, such that $\alpha(t)\leq\beta(t)$ for all $t\in I$. In these
conditions a \emph{funnel} is the set $\{(t,x)\in\mathbb{R}^2:t\in I$ and
$\alpha(t)\leq x\leq\beta(t)\}$.
\end{enumerate}
\end{definition}

\begin{figure}
\begin{center}

\includegraphics[width=6cm]{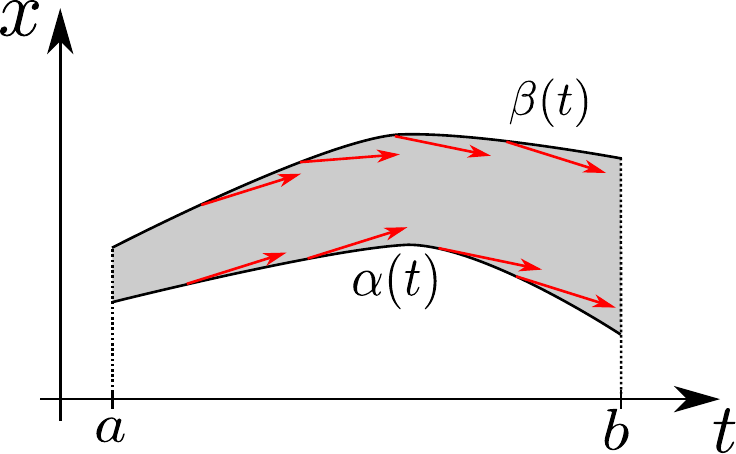}
\caption{A funnel for the ODE $x^{\prime}%
=h(t,x)$. The red arrows indicate the vector field defined by $h$ and
the gray area defines a funnel. Any solution which starts inside the funnel will stay there until $t=b$.}
\label{fig:funnel}
\end{center}
\end{figure}

The following result can be found in \cite[p. 514]{HW95} and can be depicted
graphically as in Fig.~\ref{fig:funnel}.

\begin{proposition}
\label{Prop:funnel}Let $\alpha,\beta:I=[a,b)\subseteq\mathbb{R}\rightarrow
\mathbb{R}$, where $b$ might be infinite, be a lower and upper fence,
respectively, which define a funnel for the ODE $x^{\prime}=h(t,y)$, where
$h$ is $C^{1}$. Assume also that $h$ satisfies a Lipschitz condition on the funnel. Then any solution which starts in the
funnel at $t=a$ remains in the funnel for all $t\in\lbrack a,b)$.
\end{proposition}

\begin{lemma}
\label{Lemma:Going_through_saddle}Let (\ref{ODE_Main}) define a structurally
stable system over the compact set $\mathbb{D}\subseteq\mathbb{R}^{2}$. Then
there exists $\varepsilon>0$ with the following properties. Let $x_{0}$ be an
hyperbolic saddle point and let $\Omega_{1}(x_{0}),\Omega_{2}(x_{0})$ be
defined as in Lemma \ref{Lemma:unstable_convergence}. Let also $U_{x_{0}}$ be
a local unstable manifold of $x_{0}$. Then $B(x_{0},\varepsilon)\subseteq
\mathbb{D}$ and any trajectory starting in $B(x_{0},\varepsilon)$ will either:
(i) converge to $x_{0}$ or (ii) converge to one of
the attractors $\Omega_{1}(x_{0}),\Omega_{2}(x_{0})$.
\end{lemma}

\begin{proof}
Let $x_{0}$ be some saddle point of (\ref{ODE_Main}). First let us pick some
$\varepsilon^{\ast}>0$ small enough so that $B(x_{0},\varepsilon^{\ast
})\subseteq\mathbb{D}$. This is always possible since $x_{0}$ cannot be at the
boundary of $\mathbb{D}$ or otherwise the system defined by (\ref{ODE_Main})
would not be structurally \ stable. Moreover, by Peixoto's theorem, the number
of saddle points is finite and thus we can take $\varepsilon^{\ast}>0$ to be
the minimum of all $\varepsilon^{\ast}$ for each saddle points. Since
$d(x_{0},\partial\mathbb{D})>0$ it suffices to pick some $\varepsilon^{\ast
}>0$ satisfying $\varepsilon^{\ast}<d(x_{0},\partial\mathbb{D})$ for all
saddle points $x_{0}$. We also know that, by the Hartman-Grobman theorem,
there is a homeomorphism $H$ from an open $U\subseteq\mathbb{D}$ containing
$x_{0}$ to an open subset $V\subseteq\mathbb{R}^{2}$ containing the origin
such that it maps trajectories of (\ref{ODE_Main}) in $U$ to trajectories of
the following linearized ODE%
\begin{equation}
x^{\prime}=A_{0}x \label{Eq:Linearized_system}%
\end{equation}
where $A_{0}=Df(x_{0})$ and $H(x_{0})=0$. Let us also suppose without loss of
generality that $\varepsilon^{\ast}$ is small enough so that $B(x_{0}%
,\varepsilon^{\ast})\subseteq U.$%
\begin{figure}
\begin{center}
\includegraphics{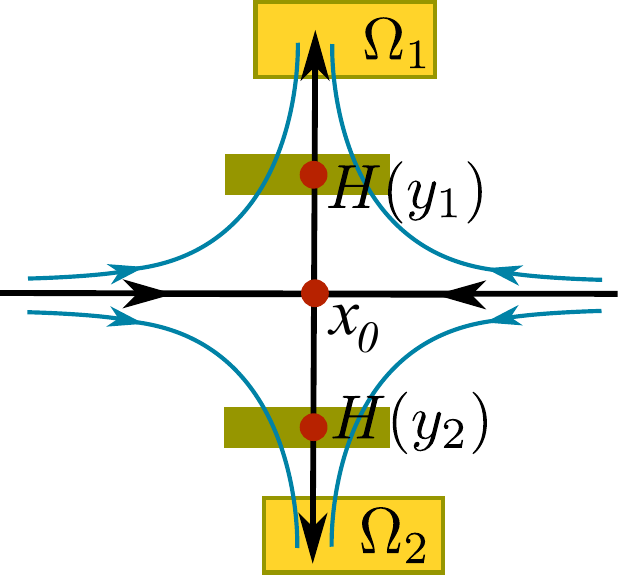}
\caption{Flow near a saddle point, converging to two attractors $\Omega
_1$ and $\Omega_2$.}
\label{fig:saddle_refined}
\end{center}
\end{figure}%

Let $\Gamma_{1},\Gamma_{2}$ be defined as in the proof of Lemma
\ref{Lemma:unstable_convergence} (i.e. they form, together with the point
$x_{0}$, a local unstable manifold of $x_{0}$). Let $y_{1}\in\Gamma_{1}\cap
B(x_{0},\varepsilon^{\ast})$. Then the trajectory $\phi_{t}(y_{1})$ will
converge to some
attractor $\Omega_{1}(x_{0})$. This means that the trajectory starting on $y_{1}$
will eventually reach a basin of attraction $\mathcal{B}_{\Omega_{1}(x_{0})}$
of the atractor $\Omega_{1}(x_{0})$ as defined in Lemma \ref{Lemma:basin} in
some time $T_{1}$ and then stay there, i.e. $\phi_{t}(y_{1})\in\mathcal{B}%
_{\Omega_{1}(x_{0})}$ for all $t\geq T_{1}$. Since $\mathcal{B}_{\Omega
_{1}(x_{0})}$ is open, and due to Lemma \ref{Lemma:unstable_convergence},
there is some $\delta_{1}>0$ such that $B(\phi_{t}(y_{1}),\delta_{1}%
)\subseteq\mathcal{B}_{\Omega_{1}(x_{0})}$ for all $t\geq T_{1}$. Because the
operator $\phi_{T_{1}}$ is continuous (Corollary \ref{Cor:Continuity}) and
because $B(\phi_{t}(y_{1}),\delta_1)$ is an open set, this implies that
\[
\phi_{T_{1}}^{-1}(B(\phi_{T_{1}}(y_{1}),\delta_{1}))=\phi_{-T_{1}}%
(B(\phi_{T_{1}}(y_{1}),\delta_{1}))
\]
is an open neighborhood of $y_{1}$ with the property that any trajectory
starting in this neighborhood will have reached the basin of attraction
$\mathcal{B}_{\Omega_{1}(x_{0})}$ in time $T_{1}.$

By the Hartman-Grobman theorem, $H\circ\phi_{-T_{1}}(B(\phi_{T_{1}}%
(y_{1}),\delta_{1}))$ will be a neighborhood of $H(y_{1})$ for $\delta_{1}$
small enough. Moreover, $H(y_{1})$ will also belong to the unstable manifold
of the origin in the linearized system (\ref{Eq:Linearized_system}).
Similarly, we can pick a point $y_{2}\in\Gamma_{2}\cap B(x_{0},\varepsilon
^{\ast})$ and $\delta_{2},T_{2}>0$ such that $\phi_{-T_{2}}(B(\phi_{T_{2}%
}(y_{2}),\delta_{2}))$ is an open neighborhood of $y_{2}$ with the property
that any trajectory starting in this neighborhood will have reached the basin
of attraction $\mathcal{B}_{\Omega_{2}(x_{0})}$ in time $T_{2}$. Furthermore,
$H\circ\phi_{-T_{2}}(B(\phi_{T_{2}}(y_{2}),\delta_{2}))$ will be a
neighborhood of $H(y_{2})$, and $H(y_{2})$ will also belong to the unstable
manifold of the origin in the linearized system (\ref{Eq:Linearized_system}).

The classical theory of ODE gives us a way of finding the solutions of
(\ref{Eq:Linearized_system}) explictly. Namely, since $x_{0}$ is an hyperbolic
saddle point, $Df(x_{0})$ will have two eigenvalues $\lambda_{1}<0$ and
$\lambda_{2}>0$. By picking appropriate coordinates (namely by picking as a
basis for the coordinate system non-zero eigenvectors $v_{1},v_{2}$ associated
to the eigenvalues $\lambda_{1},\lambda_{2}$, respectively. Without loss of
generality this change of coordinates can be assumed to be comprised in the
homeomorphism $H$), we see that this system is given by%
\begin{equation}
\left[
\begin{array}
[c]{c}%
x_{1}^{\prime}\\
x_{2}^{\prime}%
\end{array}
\right]  =\left[
\begin{array}
[c]{cc}%
\lambda_{1} & 0\\
0 & \lambda_{2}%
\end{array}
\right]  \left[
\begin{array}
[c]{c}%
x_{1}\\
x_{2}%
\end{array}
\right]  =\left[
\begin{array}
[c]{c}%
\lambda_{1}x_{1}\\
\lambda_{2}x_{2}%
\end{array}
\right]  \label{Eq:Linearized_eigenvectors}%
\end{equation}
which has as solution curves (assuming that they start at $(x_{1}%
(0),x_{2}(0))$)%
\begin{equation}
\left\{
\begin{array}
[c]{c}%
x_{1}=x_{1}(0)e^{\lambda_{1}t}\\
x_{2}=x_{2}(0)e^{\lambda_{2}t}%
\end{array}
\right.  \label{Eq:Solution_Linear}%
\end{equation}
Graphically the flow is as depicted in Fig.\ \ref{fig:saddle_refined}. Without
loss of generality, we assume that the homeomorphism $H$ maps the trajectories
of (\ref{ODE_Main}) around $x_{0}$ to trajectories of
(\ref{Eq:Linearized_eigenvectors}) around the origin. Let $\delta>0$ be such
that $B(H(y_{i}),\delta)\subseteq H\circ\phi_{-T_{i}}(B(\phi_{T_{i}}%
(y_{i}),\delta_{i}))$ for $i=1,2.$ From the expression of the solution
(\ref{Eq:Solution_Linear}) for the ODE (\ref{Eq:Linearized_eigenvectors}), we
conclude that any trajectory of (\ref{Eq:Solution_Linear}) starting on a point
$z=(z_{1},z_{2})\in B(0,\delta)$ will either: (i) converge to the origin (when
$z_{2}=0$, i.e. when $z$ lies on the stable manifold) or (ii) enter the ball
$B(H(y_{1}),\delta)$ (when $z_{2}>0$) or (iii) enter the ball $B(H(y_{2}%
),\delta)$ (when $z_{2}<0$). This implies that any trajectory starting on
$H^{-1}(B(0,\delta))$ will either converge to $x_{0}$ or enter the open set
$\cup_{i=1,2}\phi_{-T_{i}}(B(\phi_{T_{i}}(y_{i}),\delta_{i}))$. But since once
a trajectory enters $\phi_{-T_{i}}(B(\phi_{T_{i}}(y_{i}),\delta_{i}))$ it will
reach the basin of attraction $\mathcal{B}_{\Omega_{i}(x_{0})}$, we conclude
the desired result for this case.
\end{proof}

\textbf{Proof of Theorem \ref{Thm:No_saddle}.} By Lemma
\ref{Lemma:Going_through_saddle}, and the arguments used in its proof, we know
that there is some $\lambda>0$ (independent of $i)$ such that any trajectory
starting in $B(x_{i},\lambda)$ will either converge to $x_{i}$ (if the
trajectory starts on the stable manifold of $x_{i}$) without leaving
$B(x_{i},\lambda)$ or it will converge to some
attractor $\Omega_{i}$.

Since there are no saddle connections on structurally stable systems
$W^{u}(x_{i})$ cannot intersect $W^{s}(x_{j})$. Let%
\[
\delta=\min_{\substack{i,j\in\{1,\ldots,n\}\\i\neq j}}d(\overline{W^{u}%
(x_{i})},x_{j})>0.
\]
Note also that the time it takes to go from a point $y_{i}\in W^{u}(x_{i}%
)\cap\partial B(x_{i},\lambda)=\{y_{i},z_{i}\}$ to $\mathcal{B}_{\Omega_{i}}$
is finite. Moreover, because
$\mathcal{B}_{\Omega}$ is open, we conclude that there is
some $\epsilon>0$ such that $B(\phi_{t}(y_{i}),\epsilon)\subseteq\mathcal{B}_{\Omega_{i}}$. Hence, by
Lemma \ref{Lemma:Perturbation}, Corollary \ref{Cor:Continuity}, and using
techniques similar to those used in the proof of Lemma
\ref{Lemma:Going_through_saddle}, we conclude that there is some $\delta
^{\ast}>0$\ and some $T^{\ast}\geq 0$, independent of $i$ and of
$\{y_{i},z_{i}\}$, where $y_{i}$ and $z_{i}$ belong to \textquotedblleft
different sides\textquotedblright\ of $W^{u}(x_{i})\cap\partial B(x_{i}%
,\lambda)$, such that any point in $B(w,\delta^{\ast}),$ $w\in\{y_{i}%
,z_{i}\},$ will be on
$\mathcal{B}_{\Omega}$ for some attractor $\Omega$ at time $T^{\ast}$ and,
moreover, $d(\phi_{t}(B(w,\delta^{\ast})),\phi_{t}(w))<\delta/2$ for any
$0\leq t\leq T^{\ast}.$ Since $d(\phi_{t}(w),x_{j})\geq\delta$ for all $0\leq
t\leq T^{\ast}$, this implies that $d(\phi_{t}(B(w,\delta^{\ast}%
)),x_{j})>\delta/2$. Therefore, it suffices to take $\varepsilon=\min
(\delta^{\ast},\delta/2)$ to conclude that no trajectory leaving from
$B(x_{i},\varepsilon)$ will enter $B(x_{j},\varepsilon)$ for $i\neq j.$

\section{Error bounds on the plane for Euler's method}\label{Sec:AppendixEuler}

In this appendix we derive the error bound (\ref{Eq:ErrorEuler}) for the Euler
method used in Section \ref{Sec:Simulation}.

We recall that in Euler's method we start from a point $x_{0}$ and then obtain
several iterates $x_{1},\ldots,x_{N}$ which approach the solution of the
initial-value problem $z^{\prime}=f(z),$ $x(0)=x_{0}$ at times $0,h,\ldots
,Nh\subseteq\lbrack0,b]$, where $b>0$, $h>0$ is the time step, and $x_{i}$
approximates $z(ih)$ for $i=1,\ldots,N$. Here we follow the arguments used in
\cite[pp. 346--350]{Atk89}\ for obtaining error bounds for the one-dimensional case to get error bounds for the two-dimensional case.

If $f=(f_{1},f_{2})$ and $z=(z_{1},z_{2})\in\mathbb{R}^{2}$, Taylor's formula
gives us%
\[
z_{1}(t_{0}+h)=z_{1}(t_{0})+hz_{1}^{\prime}(t_{0})+\frac{h^{2}}{2}%
z_{1}^{\prime\prime}(\xi)
\]
for some $t_{0}<\xi<t$. Note that $z_{1}^{\prime}(\xi)=f_{1}(z(\xi
))=f(z_{1}(\xi),z_{2}(\xi))$. Hence
\begin{align*}
z_{1}^{\prime\prime}(\xi)  &  =\frac{\partial f_{1}}{\partial z_{1}}%
(z(\xi))z_{1}^{\prime}(\xi)+\frac{\partial f_{1}}{\partial z_{2}}(z(\xi
))z_{2}^{\prime}(\xi)\\
&  =\frac{\partial f_{1}}{\partial z_{1}}(z(\xi))f_{1}(z(\xi))+\frac{\partial
f_{1}}{\partial z_{2}}(z(\xi))f_{2}(z(\xi)).
\end{align*}
If we take $M=\max_{x\in\mathbb{D}}(\left\Vert f(x)\right\Vert ,\left\Vert
Df(x)\right\Vert )$, then we get%
\[
\left\Vert z_{1}^{\prime\prime}(\xi)\right\Vert \leq2M^{2}.
\]
A similar bound holds for $\left\Vert z_{2}^{\prime\prime}(\xi)\right\Vert $.
Let us now assume that the rounding error is bounded by $\rho>0$ when using
Euler's method (\ref{Eq:Euler}) and take $X_{i}=z(hi)$ for $i=0,1,\ldots,N$.
Then%
\begin{align*}
X_{i+1}  &  =X_{i}+hf(X_{i})+\frac{h^{2}}{2}(z_{1}^{\prime\prime}(\xi
_{i,1}),z_{1}^{\prime\prime}(\xi_{i,2}))\\
x_{i+1}  &  =x_{i}+hf(x_{i})+\rho_{i}%
\end{align*}
where $\rho_{i}$ is the rounding error on each step, with $\left\vert \rho
_{i}\right\vert \leq\rho$ and $i=1,\ldots,N$. Subtracting both equations, we
get%
\[
e_{i+1}=e_{i}+h\left(  f(X_{i})-f(x_{i})\right)  +\frac{h^{2}}{2}%
(z_{1}^{\prime\prime}(\xi_{i,1}),z_{1}^{\prime\prime}(\xi_{i,2}))-\rho_{i}%
\]
where $e_{i}=X_{i}-x_{i}$. Let $\tau_{i}=\frac{h}{2}(z_{1}^{\prime\prime}%
(\xi_{i,1}),z_{1}^{\prime\prime}(\xi_{i,2}))-\frac{\rho_{i}}{h}$. Then the
above identity yields (note that since $M\geq\left\Vert Df(x)\right\Vert $ on
$\mathbb{D}$, it works as a Lipschitz constant there, but let us just consider
a Lipschitz constant $L>0$ there)%
\begin{align*}
\left\Vert e_{i+1}\right\Vert  &  \leq\left\Vert e_{i}\right\Vert +h\left(
L\left\Vert e_{i}\right\Vert \right)  +h\left\Vert\tau_{i}\right\Vert\\
&  \leq(1+hL)\left\Vert e_{i}\right\Vert +h\left\Vert\tau_{i}\right\Vert%
\end{align*}
Applying this last formula recursively, we get%
\begin{align*}
\left\Vert e_{i}\right\Vert  &  \leq(1+hL)^{i}\left\Vert e_{0}\right\Vert
+\left(  1+(1+hL)+\ldots+(1+hL)^{i-1}\right)  h\left\Vert\tau_{i}\right\Vert\\
&  \leq(1+hL)^{i}\left\Vert e_{0}\right\Vert +\left(  1+(1+hL)+\ldots
+(1+hL)^{i-1}\right)  h\tau
\end{align*}
where $\tau=hM^{2}+\rho/h$. Using the formula for the sum of a finite
geometric series $(r\neq1)$%
\[
1+r+r^{2}+\ldots+r^{n-1}=\frac{r^{n}-1}{r-1}%
\]
we get%
\[
\left\Vert e_{i}\right\Vert \leq(1+hL)^{i}\left\Vert e_{0}\right\Vert +\left(
\frac{(1+hL)^{i}-1}{L}\right)  \tau.
\]
Since%
\[
(1+hL)^{i}\leq e^{ihL}\leq e^{bL}%
\]
The last inequality then yields%
\[
\left\Vert e_{i}\right\Vert \leq e^{bL}\left\Vert e_{0}\right\Vert +\left(
\frac{e^{bL}-1}{L}\right)  \left(  hM^{2}+\frac{\rho}{h}\right)
\]
which is essentially (\ref{Eq:ErrorEuler}).
\end{document}